\documentclass[12pt, a4paper]{article}

\usepackage[english]{babel}
\usepackage{amsmath}
\usepackage{amsthm}
\usepackage{amsfonts}
\usepackage{amssymb}
\usepackage[title,toc,titletoc,page]{appendix}
\usepackage{float}
\usepackage{hyperref}
\usepackage{algorithm}
\usepackage{algpseudocode}
\usepackage{enumitem}
\usepackage{titlesec}
\usepackage{secdot}
\sectiondot{subsection}
\sectiondot{subsubsection}
\titleformat*{\section}{\Large\bfseries\fontfamily{cmss}\selectfont}

\titleformat*{\subsection}{\large\bfseries\fontfamily{cmss}\selectfont}
\titleformat*{\subsubsection}{\bfseries\fontfamily{cmss}\selectfont}
\setlist[enumerate]{topsep=0pt, itemsep=.5ex, parsep=0ex}
\usepackage{tkz-tab}
\usepackage{fancyhdr}
\usepackage{color}
\usepackage[colorinlistoftodos]{todonotes}

\hypersetup{colorlinks=true, allcolors = blue}
\usepackage{tcolorbox}
\usepackage{subcaption}
\usepackage{titlesec}
\usepackage{relsize}
\usepackage{scalerel}
\titleformat*{\section}{\large\bfseries}

\tcbuselibrary{theorems}

\newtcbtheorem{mytheorem}{\emph{Định lý}}{colback = white!30, colframe = black!90, fonttitle = \bfseries}{th}
\usepackage{tabularx}
\usepackage[bitstream-charter]{mathdesign}

\usepackage[left=2.5cm,right=2.5cm,top=1.54cm,bottom=1.54cm,headheight=1cm]{geometry}
\usepackage{xcolor}
\hypersetup{
    colorlinks,
    linkcolor={red!70!black},
    citecolor={blue!70!black},
    urlcolor={blue!70!black}
}
\usepackage{mathtools}
\usepackage{multicol}
\usepackage{amsthm}
\usepackage{mathrsfs}
\usepackage{multirow}
\usepackage{relsize}
\usepackage{fancyhdr}

\theoremstyle{plain}
\newtheorem{definition}{Definition}
\theoremstyle{plain}

\usepackage[mathcal]{euscript}
\numberwithin{hq}{section}

\usepackage{pgf}
\usepackage{setspace}
\usepackage{mathrsfs}
\usepackage{tikz}
\usetikzlibrary{arrows.meta,positioning}
\newcommand{\norm}[1]{\left\lVert#1\right\rVert}

\makeatletter
\newcommand{\leqnomode}{\tagsleft@true}
\newcommand{\reqnomode}{\tagsleft@false}
\makeatother
\usepackage{titling}
\usepackage{multirow}
\thanksmarkseries{arabic}
\makeatletter
\renewcommand{\ALG@name}{Algorithm}
\makeatother

\usepackage{makecell}
\usepackage{array, booktabs}

\begin{document}
\title{Finite Expression Method with TranNet-based Function Learning for High-Dimensional Partial Differential Equations}
\author{Toan Huynh, Feng Bao, Haizhao Yang, Ahmed Zytoon}

\maketitle
\section*{Abstract}
In this paper, we study a machine-learning-based solver for high-dimensional partial differential equations (PDEs). Computing accurate solutions efficiently for such problems remains challenging because of the curse of dimensionality, which severely limits the scalability of classical numerical methods. Our approach builds on the recently developed finite expression method (FEX), which approximates PDE solutions in a function space generated by finitely many analytic expressions. This framework has been shown to achieve high, and in some cases machine-level, accuracy with polynomial memory complexity and favorable computational cost. We propose an extension of FEX in which the functional pool is generated by shallow neural network operators whose parameters are initialized using the transferable neural network method TransNet. Numerical experiments suggest that the proposed extension is an effective alternative for solving several high-dimensional PDEs.
\section{Introduction}
Partial differential equations (PDEs) play a fundamental role in many scientific fields because they are effective in modeling various physical phenomena, such as diffusion~\cite{Crank1975}, fluid dynamics~\cite{Temam2001}, and ocean modeling~\cite{Miller2007}. Numerous classical numerical methods, including finite difference, finite element, and finite volume methods, have been developed to compute approximate solutions of PDEs with great success~\cite{Adler2024}. However, for high-dimensional problems, the computational cost of these methods often grows exponentially with the dimension of the computational domain~\cite{Weinan2022}. As a result, there has been a surge of interest in using neural networks (NNs) to develop mesh-free methods that leverage data and physical knowledge to solve PDEs~\cite{gin2021deepgreen, Han2018, Khoo2021, Li2021, Lu2021, Raissi2019, Sirignano2018,  Weinan2017, Weinan2018, Zang2020, Zhang2022, Zhang2024TransNet}. Along the same lines, extensive research has also been devoted to applying NNs to SDEs, SPDEs, and stochastic dynamical systems, particularly for simulation, prediction, and uncertainty quantification~\cite{Richard2022, Richard2024, baranek2026stpinnsdeeplearning, Bahmani2025, liu2019neuralsde, Kong2020sdenet, tzen2019neural, jia2019neural, Ryan2026, Tang2025}. From a theoretical perspective, NNs possess strong approximation capabilities and can, in certain settings, help mitigate the curse of dimensionality~\cite{Zhang2022, Jiao2023relu, Shen2021a, Shen2021b, Yarotsky2021elementary}. However, in practice, NN-based solvers often struggle to achieve highly accurate solutions, especially for very high-dimensional problems, even when the target solution has relatively simple structure~\cite{Weinan2018, Liu2020mscalednn}. In addition, such methods may still incur substantial memory usage and computational cost~\cite{Bianco2018benchmark}.

In this work, we employ a recently developed methodology, namely the finite expression method (FEX), which seeks solutions in a function space of mathematical expressions constructed from finitely many operators~\cite{Liang2025finite}. Analytical and numerical results have shown that FEX can reproduce the true solution with high, and in some cases machine-level, accuracy, while mitigating or even avoiding the curse of dimensionality in certain settings. In FEX, the search for a PDE solution is formulated as a combinatorial optimization problem over mathematical expressions built from unary and binary operators chosen from a candidate pool. The goal is to determine both the discrete operator choices that define the expression structure and the continuous parameters associated with them. In principle, such combinatorial optimization (CO) problems can be addressed by traditional approaches such as genetic programming and simulated annealing~\cite{Bello2016neural, Cheung2019thompson, Mazyavkina2021reinforcement, Murraysmith2012modelling}. More recently, reinforcement learning (RL) has emerged as an effective alternative for solving CO problems by learning solution-construction strategies directly from reward feedback, without relying heavily on hand-crafted heuristics. Motivated by the recent success of RL in combinatorial optimization, automatic algorithm design, and symbolic regression~\cite{Bello2017neural, Coreyes2021evolving, Landajuela2021discovering, Petersen2021deep, Ramachandran2018searching}, FEX employs an RL-based search loop to identify promising operator configurations, after which the associated continuous parameters are optimized.

The original FEX framework uses a candidate pool built from explicit, hand-crafted formulas such as \((\cdot)^2\) and \((\cdot)^3\). While effective, this design restricts the search space to operators specified \emph{a priori}. In this work, we investigate an extension in which some operators are no longer given by explicit formulas, but are instead represented by trainable single-hidden-layer neural networks. Once trained, these networks define nonlinear operators that can be incorporated into the FEX search procedure alongside the original hand-crafted candidates. To train these neural operators, we adopt the training methodology introduced in the transferable neural network (TransNet) framework~\cite{Zhang2024TransNet, Zhang2024}.

TransNet is a transfer learning method for PDEs that constructs a pre-trained neural feature space designed to be reused across PDEs with different coefficients, domains, and equation types. Unlike several existing transfer learning approaches for PDEs, such as~\cite{Li2021, Lu2021, Souvik2020, desai2022oneshottransferlearningphysicsinformed}, which require prior knowledge of the target PDE family during pre-training, TransNet builds this feature space without using PDE-specific information. As a result, the learned feature space can be transferred across a broad class of PDEs with different domains and boundary conditions. The feature space is initialized from a random feature model~\cite{Chen2022, Liu2023, Sun2018} and is subsequently fine-tuned using least-squares solvers that minimize the standard PDE residual loss. As shown in~\cite{Zhang2024TransNet}, TransNet achieves competitive accuracy and efficiency on a variety of linear and nonlinear PDEs. However, because TransNet was not designed for high-dimensional PDEs, we use only its training methodology here to optimize neural operators in the FEX candidate pool.

The rest of the paper is organized as follows. In Section~\ref{Sec2_FEX}, we review the methodology of FEX and its implementation. In Section~\ref{Sec3_TransNet}, we introduce the TransNet method together with its training procedure and explain how it can be used to approximate a given function. We then present numerical results in Section~\ref{Sec4_Results} to demonstrate the effectiveness of the proposed algorithm. Finally, we conclude in Section~\ref{Conclusion}.
\section{Finite Expression Method (FEX)}
\label{Sec2_FEX}
\subsection{Methodology}
The goal of FEX is to find a mathematical expression that approximates the solution to a given PDE. Numerical results in~\cite{Liang2025finite} show that FEX can recover the expression of the true solution to very high precision. We first recall the notation introduced in~\cite{Liang2025finite} to define the function space of mathematical expressions and to formulate the FEX problem as a combinatorial optimization (CO) problem. 

\begin{definition}[\textbf{Mathematical expression}]
A mathematical expression is a combination of symbols that is well formed according to syntax rules and defines a valid function. These symbols include operands (variables and numbers), operators (e.g.\ ``$+$'', ``$\sin$'', integral, derivative), brackets, and punctuation.
\end{definition}

For example, ``$\sin(x \times y)+1$'' is a valid mathematical expression. On the other hand, expressions such as ``$5>x$'' and ``$\sin(x \times y)+$'' are not included, since they do not define valid functions.

\begin{definition}[\textbf{$k$-finite expression}]
A mathematical expression is called a $k$-finite expression if the number of operators appearing in the expression is $k$.
\end{definition}

For instance, the expression ``$\sin(x \times y)+1$'' is a $3$-finite expression since it contains the three operators ``$\times$'', ``$\sin$'', and ``$+$''. Infinite mathematical expressions, such as series expansions, are not included in this framework.

\begin{definition}[\textbf{Finite expression method}]
The finite expression method is a numerical methodology for solving a PDE by seeking a finite expression whose associated function approximately satisfies the PDE.
\end{definition}

We denote by $\mathbb{S}_k$ the set of functions represented by finite expressions containing at most $k$ operators. This defines the function space used in FEX. The idea of FEX is to approximate the solution of a given PDE by solving a combinatorial optimization (CO) problem over the space $\mathbb{S}_{k}$.

More precisely, let $\mathcal{L}: \mathbb{S} \to \mathbb{R}$ be a functional associated with a given PDE, where $\mathbb{S}$ is a function space. The minimizer of $\mathcal{L}$ is then viewed as the best approximation to the PDE solution within $\mathbb{S}$. Accordingly, the target CO problem over 
$\mathbb{S}_k$ is defined by:
\vspace{-0.2cm}
\begin{equation}
\label{CO_FEX}
\min\{ \mathcal{L}(u) \mid u \in \mathbb{S}_k \}.
\vspace{-0.3cm}
\end{equation}

The choice of the functional $\mathcal{L}$ depends on the problem under consideration, and one may design a more suitable functional for a specific PDE, constraint, or domain. Some popular choices include least-square methods~\cite{Sirignano2018, Lagaris1998,  Dissanayake1994} and variational formulations~\cite{Weinan2018, Zang2020, Yulei2021, Chen2023} For a detailed discussion on the choice of $\mathcal{L}$, we refer the reader to~\cite{Liang2025finite}.

It was shown in~\cite{Liang2025finite} that the function space $\mathbb{S}_k$ can avoid the curse of dimensionality in the approximation of high-dimensional continuous functions. In particular, the space $\mathbb{S}_k$, generated by operators such as ``$+$'', ``$-$'', ``$\times$'', ``$/$'', ``$\lvert\cdot\rvert$'', ``$\mathrm{sign}(\cdot)$'', and ``$\lfloor \cdot \rfloor$'', was proved to be dense in $C([a,b]^d)$ for arbitrary $a,b \in \mathbb{R}$ and $d \in \mathbb{N}$ whenever $k > \mathcal{O}(d^4)$; see \cite[Theorem 4]{Liang2025finite}. However, since ``$\mathrm{sign}(\cdot)$'' and ``$\lfloor \cdot \rfloor$'' are not commonly used in mathematical expressions, it is more practical to restrict the operator list to more standard operators such as ``$+$'', ``$\times$'', ``$\sin(\cdot)$'', and exponential functions. For spaces $\mathbb{S}_k$ generated by such operators, it was also proved that $\mathbb{S}_k$ can still approximate the H\"older class $\mathcal{H}^{\alpha}_{\mu}\left([0,1]^d\right)$ while avoiding the curse of dimensionality, provided that $k \geq \mathcal{O}\left(d^2\left(\log d + \log \dfrac{1}{\epsilon}\right)^2\right)$ for arbitrarily small $\epsilon > 0$; see \cite[Theorem 5]{Liang2025finite}.

In the next section, we review the numerical implementation of FEX through the introduction of binary trees and a policy-gradient-based reinforcement learning approach for solving the CO problem.

\subsection{An Implementation of FEX}
We initiate FEX by constructing binary trees that represent various finite expressions that can be used to provide an approximation to the solution of a given PDE. To identify the expression that achieves the best accuracy, we formulate the combinatorial optimization (CO) problem~\eqref{CO_FEX} as a joint parameter and operator selection problem. To solve this CO problem, a search loop is introduced to identify effective operators that have the potential to recover the true solution when incorporated into the expression. For the reader's convenience, we summarize the key notations used in this section in Table~\ref{table_notations}.
\begin{table*}[h!]
\centering
\begin{tabular}{ll}
\hline
\textbf{Notation} & \textbf{Explanation} \\
\hline
$\pmb{\mathcal{T}}$      & A binary tree \\
$\pmb{e}$ & Operator sequence \\
$\pmb{\theta}$ & Trainable scaling and bias parameters \\
$\mathcal{L}$ & The functional associated with the PDE solution \\
$S$ & The scoring function that maps an operator sequence to $[0, 1]$\\
$\pmb{\chi}_{\Phi}$ & The controller parameterized by $\Phi$ \\
$\pmb{\mathcal{J}}$ & The objective function for the policy-gradient approach\\
\hline
\end{tabular}
\caption{A summary of notations in the FEX implementation.}
\label{table_notations}
\end{table*}
\vspace{-0.9cm}
\subsection*{Finite Expressions with Binary Trees}
Finite expression method (FEX) represents a candidate solution by using a \emph{binary tree}. The main idea is to encode an analytic expression as a tree whose nodes are filled with operators chosen from prescribed unary and binary operator sets.

Let $\pmb{\mathcal{T}}$ be a binary tree. Each node of the tree is assigned either a unary operator or a binary operator, and all these node values collectively
form an operator sequence $\pmb{e}$, following a predefined order to traverse the tree (e.g., inorder
traversal). Within each node featuring a unary operator, we incorporate two additional parameters, a scaling parameter $\alpha$ and a bias parameter $\beta$, to enhance expressiveness. All these parameters are denoted by $\pmb{\theta}$. Once the operators are assigned, the tree defines a function $u(x;\pmb{\mathcal{T}},\pmb{e},\pmb{\theta})$ as a finite expression. For a fixed tree structure $\pmb{\mathcal{T}}$, the number of operators is bounded above by a constant denoted by $k_{\mathcal{T}}$, and the collection $\{u(\cdot;\pmb{\mathcal{T}},\pmb{e},\pmb{\theta}) \mid \pmb{e}, \pmb{\theta}\}$ forms the function space used in the CO to solve the PDE. 

\begin{figure*}[t]
\centering
\resizebox{0.95\textwidth}{!}{
\begin{tikzpicture}[
>=Stealth,
unary/.style={
    draw,
    rectangle,
    rounded corners=1pt,
    minimum width=7mm,
    minimum height=5mm,
    inner sep=1pt,
    font=\small
},
binary/.style={
    draw,
    circle,
    minimum size=5.5mm,
    inner sep=0pt,
    font=\small
},
every node/.style={font=\small},
lab/.style={font=\scriptsize},
expr/.style={font=\small, align=center},
depthlab/.style={font=\small}
]

% -------------------------------------------------
% Left labels
% -------------------------------------------------
\node[align=center] at (-1.8,  0.8) { \fontsize{8pt}{8pt}\selectfont Binary\\ \fontsize{8pt}{8pt}\selectfont Tree};
\node[align=center] at (-1.8, -0.9) {\fontsize{8pt}{8pt}\selectfont Math\\ \fontsize{8pt}{8pt}\selectfont Expression};
\node[align=center] at (-1.8, -1.8) {\fontsize{8pt}{8pt}\selectfont Depth};

% -------------------------------------------------
% Legend
% -------------------------------------------------
\node[unary]  at (0.2, 2.8) {\fontsize{8pt}{8pt}\selectfont $u$};
\node[right=2mm] at (0.35, 2.8) {\fontsize{8pt}{8pt}\selectfont Unary operator};

\node[binary] at (3.1, 2.8) {\fontsize{8pt}{8pt}\selectfont $b$};
\node[right=2mm] at (3.25, 2.8) {\fontsize{8pt}{8pt}\selectfont Binary operator};

\node at (6.2, 2.8) {\fontsize{8pt}{8pt}\selectfont $i$: input};
\node at (0.55, 2.2) {\fontsize{8pt}{8pt}\selectfont $o$: output};
\node at (3.7, 2.2) {\fontsize{8pt}{8pt}\selectfont $\alpha$: coefficient};
\node at (6.4, 2.2) {\fontsize{8pt}{8pt}\selectfont $\beta$: constant};

% =================================================
% Tree 1: unary, depth 1
% =================================================
\begin{scope}[xshift=0cm]
\node[unary] (u1) at (0,0.8) {\fontsize{8pt}{8pt}\selectfont $u_1$};
\draw[->] (0,0.15) -- node[pos=-0.5,lab] {$i_0$} (u1.south); 
\draw[->] (u1.north) -- node[pos=1.2,lab] {$o_1$} (0,1.4);

\node[expr] at (0,-0.9) {\fontsize{8pt}{8pt}\selectfont $\alpha_1 u_1(i_0)+\beta_1$};
\node[depthlab] at (0,-1.8) {\fontsize{8pt}{8pt}\selectfont$L=1$};
\end{scope}

% =================================================
% Tree 2: binary, depth 1
% =================================================
\begin{scope}[xshift=2.4cm]
\node[binary] (b1) at (0,0.8) {\fontsize{8pt}{8pt}\selectfont $b_1$};
\draw[->] (-0.45,0.3) -- node[pos=-0.45,lab] {$i_{01}$} (b1.south west);
\draw[->] ( 0.45,0.3) -- node[pos=-0.47,lab] {$i_{02}$} (b1.south east);
\draw[->] (b1.north) -- node[pos=1.2,lab] {$o_1$} (0,1.4);

\node[expr] at (0,-0.9) {\fontsize{8pt}{8pt}\selectfont $b_1(i_{01},i_{02})$};
\node[depthlab] at (0,-1.8) {\fontsize{8pt}{8pt}\selectfont$L=1$};
\end{scope}

% =================================================
% Tree 3: depth 2
% =================================================
\begin{scope}[xshift=6.0cm]
\node[binary] (b2) at (0,1.6) {\fontsize{8pt}{8pt}\selectfont $b_2$};
\node[unary]  (u21) at (-0.55,0.9) {\fontsize{8pt}{8pt}\selectfont $u_1$};
\node[unary]  (u22) at ( 0.55,0.9) {\fontsize{8pt}{8pt}\selectfont $u_2$};

\draw[->] (u21.north) -- (b2.south west);
\draw[->] (u22.north) -- (b2.south east);

\draw[->] (-0.55, 0.3) -- node[pos=-0.45,lab] {$i_{01}$} (u21.south);
\draw[->] ( 0.55,0.3) -- node[pos=-0.45,lab] {$i_{02}$} (u22.south);
% \draw[->] (b2.north) -- node[left,lab] {$o_2$} (0,1.9);

\node[expr] at (0,-0.9)
{\fontsize{8pt}{8pt}\selectfont $b_2\!\bigl(\alpha_1u_1(i_{01})+\beta_1,\alpha_2u_2(i_{02})+\beta_2\bigr)$};
\node[depthlab] at (0,-1.8) {\fontsize{8pt}{8pt}\selectfont $L=2$};
\end{scope}

% =================================================
% Tree 4: depth 3
% =================================================
\begin{scope}[xshift=11.1cm]
\node[unary]  (u3)  at (0, 2.1) {\fontsize{8pt}{8pt}\selectfont $u_3$};
\node[binary] (b3)  at (0,1.25) {\fontsize{8pt}{8pt}\selectfont $b_2$};
\node[unary]  (u31) at (-0.55, 0.5) {\fontsize{8pt}{8pt}\selectfont $u_1$};
\node[unary]  (u32) at ( 0.55, 0.5) {\fontsize{8pt}{8pt}\selectfont $u_2$};

\draw[->] (b3.north) -- (u3.south);
\draw[->] (u31.north) -- (b3.south west);
\draw[->] (u32.north) -- (b3.south east);

\draw[->] (-0.55, -0.1) -- node[pos=-0.45,lab] {$i_{01}$} (u31.south);
\draw[->] ( 0.55, -0.1) -- node[pos=-0.45,lab] {$i_{02}$} (u32.south);

\node[expr] at (0,-0.9)
{\fontsize{8pt}{8pt}\selectfont $\alpha_3u_3\!\left(b_2\!\bigl(\alpha_1u_1(i_{01})+\beta_1,\alpha_2u_2(i_{02})+\beta_2\bigr)\right)+\beta_3$};
\node[depthlab] at (0,-1.8) {\fontsize{8pt}{8pt}\selectfont $L=3$};
\end{scope}

\end{tikzpicture}
}
\caption{\small Examples of binary trees of depths 1, 2, and 3, respectively. In each binary tree, every node contains either a unary or a binary operator. For trees with depth greater than 1, the computation is carried out recursively. This figure is adapted from~\cite[Figure 3]{Liang2025finite}.}
\label{tree_schematic}
\vspace{-0.4cm}
\end{figure*}
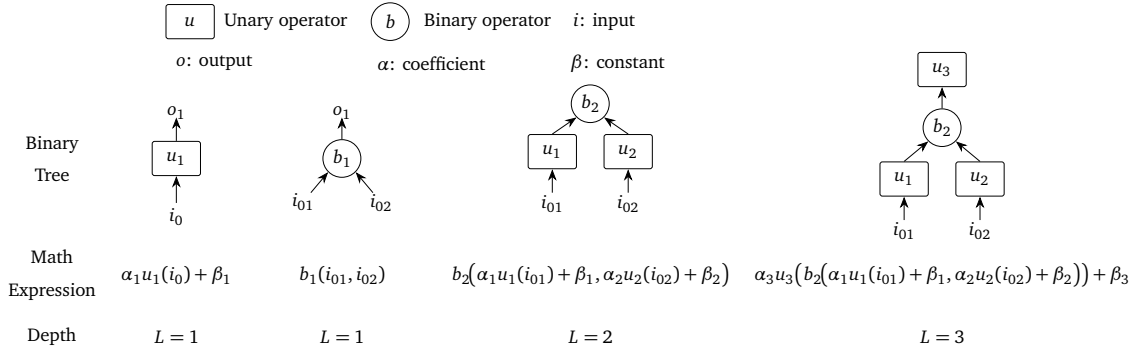

The expression is built recursively from the leaves to the root. If a node is assigned a \emph{unary operator}, then the output at that node is obtained by applying the unary operator to the output of its child, together with additional trainable scaling and bias terms. If a node is assigned a \emph{binary operator}, then the output at that node is obtained by combining the outputs of its two children. For example, for a depth-1 unary node, the output has the form $o_1 = \alpha\, u_1(i_0) + \beta$, where $u_1$ is the unary operator, $i_0$ is the input to that node, and $\alpha,\beta$ are trainable parameters. For a depth-1 binary node, the output is $o_1 = b_1(i_{01}, i_{02})$, where $b_1$ is the binary operator and $i_{01}, i_{02}$ are the two inputs. For deeper trees, this rule is applied recursively: the outputs of lower-level nodes become the inputs of higher-level nodes until the root produces the final expression value. The configuration of the tree of depths 1, 2 and 3 are displayed in Figure~\ref{tree_schematic}.

The binary operator set can consist of $\{+, -, \times, \div, \hdots\}$, and the unary operator set may include elementary and differential-type operators such as $\left\{\sin,\exp,\log,\mathrm{Id},(\cdot)^k,\frac{d\,\cdot}{dx_i},\hdots\right\}$. Because of this design, the method can represent not only standard algebraic expressions but also more structured formulas involving derivatives or integrals.

% An important feature is that unary nodes are not purely symbolic. Each unary node carries continuous parameters, typically a scaling factor $\alpha$ and a bias term $\beta$. Therefore, the search is performed not only over operator choices but also over continuous coefficients. This makes the representation more flexible than a purely discrete symbolic tree.

When the input is $x \in \mathbb{R}^d$, the unary operators attached directly to $x$ act componentwise. In that case, the scaling parameter can be vector-valued, so the tree can naturally process high-dimensional variables while still producing a scalar output. This is one of the ways the binary-tree framework is adapted to high-dimensional PDE problems.
\subsection*{Solving a CO in FEX}
We next focus on solving the CO problem~\eqref{CO_FEX} over the finite-expression function class $\{u(\cdot; \pmb{\mathcal{T}}, \pmb{e}, \pmb{\theta}) \mid \pmb{e}, \pmb{\theta}\} \subset \mathbb{S}_{k_{\mathcal{T}}}$. The target mathematical expression is then identified by minimizing the functional $\mathcal{L}$ associated with the given PDE, that is,
\begin{equation}
\label{true_CO_FEX}
\min \left\{ \mathcal{L}(u(\cdot; \pmb{\mathcal{T}}, \pmb{e}, \pmb{\theta})) \,\middle|\, \pmb{e}, \pmb{\theta} \right\}.
\end{equation}
A search loop is first executed to find a good operator sequence $\pmb{e}$ that may uncover the structure of the true solution. Once such a sequence is identified, the parameter $\pmb{\theta}$ is optimized to further minimize the functional in~\eqref{true_CO_FEX}. The search loop consists of four parts: $(1)$ score computation, which evaluates the score of the operator sequence $\pmb{e}$; $(2)$ operator sequence generation, which generates high-score operator sequences using a pre-chosen controller; $(3)$ controller update, which updates the controller using score feedback from the generated sequences to increase the probability of producing good future operator sequences; and $(4)$ candidate optimization. During the search, we maintain a pool with a fixed number of candidate to store operator sequences with high scores.

\subsubsection*{Score Computation}
The \emph{score} of the operator sequence $\pmb{e}$, which guides the controller toward generating good operator sequences and helps to maintain a candidate pool of high scores, is defined by
\begin{equation}
S(\pmb{e}):=(1+\mathcal{L}(\pmb{e}))^{-1} \in [0, 1].
\end{equation}
where $\mathcal{L}(\pmb{e}):=\min\limits_{\pmb{\theta}} \mathcal{L}\bigl(u(\cdot;\pmb{\mathcal{T}}, \pmb{e}, \pmb{\theta})\bigr)$ is the best loss achievable by the operator sequence $\pmb{e}$ after optimizing the continuous parameters $\theta$. If $\mathcal{L}(\pmb{e})$ is small, then the corresponding expression can fit the PDE well after tuning $\pmb{\theta}$, so $S(\pmb{e})$ is close to $1$. If $\mathcal{L}(\pmb{e})$ is large, then $S(\pmb{e})$ is closer to $0$. So, the score transforms a loss-minimization problem into a reward-like quantity that is easier to use in the later reinforcement-learning framework.

In principle, computing $\mathcal{L}(\pmb{e})$ requires solving an inner optimization problem over $\theta$, which is expensive and nonconvex. Therefore, the authors in~\cite{Liang2025finite} does not compute the exact minimizer every time. Instead, the score is approximated by using a two-stage optimization for $\pmb{\theta}$: a first-order optimization method (such as stochastic gradient descent~\cite{Rumelhart1986learning} or Adam~\cite{Kingma2014adam}) for $T_1$ steps, then a second-order optimization method (e.g., the Newton method~\cite{Avriel2003nonlinear} and the Broyden-Fletcher-Goldfarb-Shanno method (BFGS)~\cite{Fletcher2000practical}) for $T_2$ additional steps. The resulting parameter is denoted by $\pmb{\theta}^{\pmb{e}}_{T_1+T_2}$, and the practical score approximation is
\begin{equation}
\label{practice_score}
S(\pmb{e})\approx \left(1+\mathcal{L}\bigl(u(\cdot;\pmb{\mathcal{T}},\pmb{e},\pmb{\theta}^{\pmb{e}}_{T_1+T_2})\bigr)\right)^{-1}.
\end{equation}

% The important idea is that the algorithm does not need the exact best score for each operator sequence. It only needs a good enough estimate to compare different sequences. In this sense, $S(\pmb{e})$ acts as an approximate reward for symbolic structures.
Because the initialization of $\pmb{\theta}$ is random, the estimated score may exhibit significant variation from run to run.

\subsubsection*{Operator Sequence Generation}
Let $\pmb{\chi}_{\Phi}$ be a controller with model parameter $\Phi$. The role of the controller $\pmb{\chi}_{\Phi}$ is to guide the generation of operator sequences. The parameter $\Phi$ is updated throughout the search loop to increase the probability of generating good operator sequences.

The controller does not directly output one fixed expression. Instead, it outputs a probability distribution over the possible operators at each node of the binary tree. More precisely, if we treat the operator assigned to node $i$ of $\pmb{\mathcal{T}}$ as a random variable, then the controller $\pmb{\chi}_{\Phi}$ outputs a probability mass function $\pmb{p}^{i}_{\Phi}$ to characterize its distribution for $i=1,\hdots,s$, where $s$ is the total number of nodes in the tree. Each operator $e_i$ is then sampled according to the distribution $\pmb{p}^{i}_{\Phi}$, and these sampled values are combined to form the full operator sequence $\pmb{e}=(e_1,e_2,\dots,e_s)$.

To avoid getting stuck too early and to enhance exploration of potentially high-score operator sequences, an $\varepsilon$-greedy strategy~\cite{Sutton2018reinforcement} is employed. More precisely, with probability $\varepsilon$, the operator at a node is sampled uniformly from the operator set, while with probability $1-\varepsilon$, it is sampled from the learned controller distribution.

A larger $\varepsilon$ encourages more exploration, while a smaller $\varepsilon$ means that the algorithm relies more on what the controller has already learned.

\subsubsection*{Controller Update}
The goal of the controller update is to guide the controller toward generating high-score operator sequences $\pmb{e}$. The controller $\pmb{\chi}_{\Phi}$ is modeled as a neural network parameterized by $\Phi$. Its training objective is to maximize the expected score of sampled operator sequences, i.e.,
\begin{equation}
\label{score_FEX}
J(\Phi):=\mathbb{E}_{\pmb{e}\sim \pmb{\chi}_\Phi}[S(\pmb{e})].
\end{equation}

The update rule for the controller is based on a policy gradient method in reinforcement learning~\cite{Liu2019darts}. The policy gradient method aims to maximize the return by optimizing a parameterized policy, and here the controller $\pmb{\chi}_{\Phi}$ plays the role of that policy. To update the parameter $\Phi$, we apply gradient ascent method to the target function~\eqref{score_FEX}. To this end, we first compute the derivative of~\eqref{score_FEX} with respect to $\Phi$ and obtain
\begin{equation}
\label{score_derivative_FEX}
\nabla_\Phi J(\Phi)
=
\mathbb{E}_{\pmb{e}\sim\pmb{\chi}_\Phi}
\left[
S(\pmb{e})\sum_{i=1}^s \nabla_\Phi \log\bigl(p_\Phi^i(e_i)\bigr)
\right],
\end{equation}
where $s$ is the total number of nodes in the tree. In practice, the expectation is approximated by an average over a batch of sampled operator sequences
$ \pmb{e}^{(1)}, \pmb{e}^{(2)}, \dots, \pmb{e}^{(N)}$ where $N$ is the batch size. Then, an approximation of~\eqref{score_derivative_FEX} is given by
\begin{equation}
\label{appro_score_derivative_FEX}
\nabla_\Phi J(\Phi)
\approx
\frac{1}{N}\sum_{k=1}^N
\left[
S\bigl(\pmb{e}^{(k)}\bigr)
\sum_{i=1}^s
\nabla_\Phi \log\bigl(p_\Phi^i(e_i^{(k)})\bigr)
\right].
\end{equation}
The controller is then updated by gradient ascent with learning rate $\eta$ using the approximate score derivative~\eqref{appro_score_derivative_FEX}:
\begin{equation}
\label{GD_update}
\Phi \leftarrow \Phi + \eta \nabla_\Phi J(\Phi),
\end{equation}
where $\eta>0$ denotes the learning rate.

To increase the probability of obtaining the best score for the sequence $\pmb{e}$, the objective function proposed in~\cite{Petersen2021deep} is employed to seek the optimal solution via:
\begin{equation}
\label{quantile_score_FEX}
J(\Phi)=\mathbb{E}_{\pmb{e}\sim\pmb{\chi}_\Phi}\{S(\pmb{e})\mid S(\pmb{e})\ge S_{\nu,\Phi}\},
\end{equation}
where $S_{\nu,\Phi}$ denotes the $(1-\nu)\times 100\%$-quantile of the score distribution under $\pmb{\chi}_{\Phi}$. As a result, instead of rewarding all samples equally, this objective emphasizes samples whose scores lie above a prescribed quantile threshold. Using a finite batch of samples, the corresponding gradient approximation of~\eqref{quantile_score_FEX} becomes
\begin{equation}
\fontsize{7.8pt}{7.8pt}\selectfont
\label{discrete_quantile_score_FEX}
\nabla_\Phi J(\Phi)\approx
\frac{1}{N}\sum_{k=1}^N
\left[
\bigl(S(\pmb{e}^{(k)})-\hat S_{\nu,\Phi}\bigr)
\mathbf{1}_{\{S(\pmb{e}^{(k)})\ge \hat S_{\nu,\Phi}\}}
\sum_{i=1}^s \nabla_\Phi \log\bigl(p_\Phi^i(e_i^{(k)})\bigr)
\right].
\end{equation}
\subsubsection*{Candidate Optimization}
As discussed above, score computation involves a nonconvex optimization problem over $\theta$, and this optimization starts from a random initialization. Consequently, a genuinely good operator sequence may sometimes receive only a mediocre score. To reduce the risk of missing a good symbolic structure, we maintain a pool of $K$ strong candidates rather than greedily keeping only the current best one. During the search stage, if the candidate pool $\mathbb{P}$ contains fewer than $K$ elements, then a newly generated sequence $\pmb{e}$ is added to the pool. If the pool is already full, then the new sequence $\pmb{e}$ is added only when its score is larger than the smallest score currently in the pool, in which case the worst candidate is removed. In this way, the pool stores the top-$K$ high-score operator sequences encountered during the search. After the search stage is completed, we do not immediately select the sequence with the highest current score as the final solution. Instead, for every candidate $\pmb{e} \in \mathbb{P}$, we re-optimize $\mathcal{L}(u(\cdot;\pmb{\mathcal{T}}, \pmb{e},\pmb{\theta}))$ over the continuous parameter $\pmb{\theta}$ using a first-order optimization method with a small learning rate for $T_3$ additional iterations. Hence, once a symbolic structure appears promising, additional effort is devoted to refining its continuous coefficients more carefully. We summarize the training procedure for FEX in Algorithm~\ref{alg:fex_fixed_tree} below.
\begin{algorithm}[h!]
\caption{FEX with a fixed tree}
\label{alg:fex_fixed_tree}
\begin{algorithmic}[1]
\Require PDE and associated functional $\mathcal{L}$; fixed tree $\pmb{\mathcal{T}}$; search-loop iterations $T$; coarse-tuning iterations $T_1$ (Adam) and $T_2$ (BFGS); fine-tuning iterations $T_3$ (Adam); pool size $K$; batch size $N$.
\Ensure Approximate solution $u(x;\pmb{\mathcal{T}}, \hat{\pmb{e}},\hat{\pmb{\theta}})$
\State Initialize the controller $\chi$ for the tree $\pmb{\mathcal{T}}$
\State Initialize the candidate pool $\mathbb{P} \gets \{\}$
\For{$t=1,\dots,T$}
\State Sample $N$ sequences $\{\pmb{e}^{(1)},\pmb{e}^{(2)},\dots,\pmb{e}^{(N)}\}$ from $\chi$
\For{$n=1,\dots,N$}
\State Coarse-tune $\mathcal{L}(u(x;\pmb{\mathcal{T}},\pmb{e}^{(n)},\pmb{\theta}))$ for $T_1+T_2$ iterations to get score $S(\pmb{e}^{(n)})$
\If{$\pmb{e}^{(n)}$ belongs to the current top-$K$ candidates in $\mathbb{P}$}
\State Add $\pmb{e}^{(n)}$ to $\mathbb{P}$
\If{$|\mathbb{P}|>K$}
\State Remove from $\mathbb{P}$ the candidate with the smallest score
\EndIf
\EndIf
\EndFor
\State Update the controller $\chi$ using~\eqref{discrete_quantile_score_FEX}
\EndFor
\ForAll{$\pmb{e}\in \mathbb{P}$}
\State Fine-tune $L(u(x;\pmb{\mathcal{T}},\pmb{e},\pmb{\theta}))$ for $T_3$ iterations
\EndFor
\State \Return the candidate expression with the smallest fine-tuning error
\end{algorithmic}
\end{algorithm}

% \begin{algorithm}[H]
% \caption{FEX with progressively expanding trees}
% \label{alg:fex_expanding_tree}
% \begin{algorithmic}[1]
% \Require Tree family $\{\mathcal{T}_1,\mathcal{T}_2,\dots\}$; error tolerance $\varepsilon$
% \Ensure Approximate solution $u(x;\hat{\mathcal{T}},\tilde{\pmb{e}},\tilde{\pmb{\theta}})$
% \ForAll{$\mathcal{T}\in\{\mathcal{T}_1,\mathcal{T}_2,\dots\}$}
%     \State Initialize the controller $\chi$ for the tree $\mathcal{T}$
%     \State Obtain $u(x;\mathcal{T},\hat{\pmb{e}},\hat{\pmb{\theta}})$ by running Algorithm~\ref{alg:fex_fixed_tree}
%     \If{$L(u(\cdot;\mathcal{T},\hat{\pmb{e}},\hat{\pmb{\theta}})) \le \varepsilon$}
%         \State \textbf{break}
%     \EndIf
% \EndFor
% \State \Return the expression with the smallest functional value
% \end{algorithmic}
% \end{algorithm}

We construct the candidate pool for FEX using functions generated by the following single-hidden-layer fully connected neural network:
\begin{equation}
\label{func_transnet}
u_{\mathrm{NN}}(\pmb{y}) := \sum\limits^M_{m=1} \alpha_m\psi_m(\pmb{y})+\alpha_0 = \sum_{m=1}^{M} \alpha_m \sigma\left(\pmb{w}_m^{\top}\pmb{y} + b_m\right) + \alpha_0,
\end{equation}
where $M$ is the number of hidden neurons, $\psi_m(\pmb{y}) := \sigma\left(\pmb{w}_m^{\top}\pmb{y} + b_m\right)$, the vector $\pmb{w}_m = (w_{m,1}, \hdots, w_{m,d})^{\top}$ and the scalar $b_m$ are the weights and bias of the $m$th hidden neuron, the vector $\pmb{\alpha} = (\alpha_0, \hdots, \alpha_M)^{\top}$ contains the weights and bias of the output layer, and $\sigma(\cdot)$ is the activation function. To ensure that each approximation function is meaningful, we train the network~\eqref{func_transnet} so that it approximates a pre-selected target function. As the training mechanism, we employ the Transferable Neural Network (TransNet) developed in~\cite{Zhang2024TransNet}, which accelerates the training process while still providing an accurate approximation.

TransNet is a transfer-learning neural network that, once trained, can be reused to solve PDEs with different coefficients, on different domains, or even across different PDE classes. Numerical results in~\cite{Zhang2024TransNet, Zhang2024} show that TransNet can be transferred effectively to a variety of PDEs and domains with low computational cost. In addition, it makes more effective use of the expressive power of a given neural network and outperforms baseline methods such as PINNs. However, extending this approach to very high-dimensional problems is not straightforward. In this work, we overcome this difficulty by integrating TransNet into the FEX framework. In the next section, we briefly recall the key ingredients of TransNet.
\section{Transferrable Neural Networks}
\label{Sec3_TransNet}
If we treat the set of hidden neurons $\left\{\psi_m(\pmb{y})\right\}^M_{m=1}$ as a globally supported basis in $\mathbb{R}^d$, then the function generated by~\eqref{func_transnet} lives in a linear space extended by the basis $\left\{\psi_m(\pmb{y})\right\}^M_{m=1}$, namely the \textit{neural feature space},
\begin{equation}
\label{feature_space}
\mathcal{P}_{\mathrm{NN}} = \text{span} \left\{1, \psi_1(\pmb{y}), \hdots, \psi_M(\pmb{y})\right\}.
\end{equation}
In the standard PINN method, the coefficient $\pmb{\alpha}$ and the space $\mathcal{P}_{\mathrm{NN}}$ are simultaneously trained using SGD methods, which often leads to a non-convex and ill-conditioned optimization problem. It has been shown that the non-convexity and ill-conditioning in the neural network training are major reasons for the unsatisfactory accuracy of the trained neural network~\cite{Basri2019convergence, Cao2021spectral, Xu2020}. The main idea of TransNet is to separate the training of $\mathcal{P}_{\mathrm{NN}}$ from that of $\pmb{\alpha}$ and to construct a single neural feature space $\mathcal{P}_{\mathrm{NN}}$ that can be used for various PDEs. As a result, different target functions can be approximated without restarting the entire training process. The main procedure of TransNet can be summarized in four steps.
\subsection{Reparameterization of $\mathcal{P}_{\mathrm{NN}}$}

The first step is to reparameterize the weights and bias of each hidden neuron $\sigma\left(\pmb{w}_m^{\top}\pmb{y}+b_m\right)$ into distinct parameters that determine the \textit{location} and \textit{shape} of the neuron. More precisely, we rewrite each hidden neuron as
\begin{equation}
\sigma\left(\pmb{w}_m^{\top}\pmb{y}+b_m\right) = \sigma\left(\gamma_m \left(\pmb{a}_m^{\top}\pmb{y}+r_m\right)\right),
\end{equation}
where $(\pmb{a}_m$, $r_m)$ is location parameters and $\gamma_m$ is the shape parameter and are given by
\begin{equation}
\label{reparametrize}
\begin{array}{l}
\pmb{a}_m = \dfrac{\pmb{w}_m}{\norm{\pmb{w}_m}_2}, \vspace{0.1cm} \\
r_m = \dfrac{b_m}{\norm{\pmb{w}_m}_2}, \vspace{0.1cm} \\
\gamma_m = \norm{\pmb{w}_m}_2.
\end{array}
\end{equation}
After reparameterization, each of these parameters has a distinct geometric interpretation for the basis function $\psi_m(\pmb{y}) := \sigma\left(\pmb{w}_m^{\top}\pmb{y}+b_m\right)$. More precisely, when $\sigma$ is the ReLU activation function, each hidden neuron is associated with a \textit{partition hyperplane}, which separates the activated and unactivated regions and is defined by $\pmb{w}_m^{\top}\pmb{y}+b_m=0$. This partition hyperplane can be equivalently represented by the location parameter through the equation $\pmb{a}_m^{\top}\pmb{y}+r_m=0$. When $\sigma$ is the $\tanh(\cdot)$ activation function, the partition hyperplane $\pmb{a}_m^{\top}\pmb{y}+r_m=0$ represents the transition rate from $-1$ to $1$ on the basis function $\psi_m(\pmb{y})$. For further discussion of the geometric roles of the location and shape parameters with respect to this activation function, we refer the interested reader to~\cite{Zhang2024TransNet}.

The purpose of the re-parametrization step is to effectively segregate location and shape parameters, each possessing almost independent geometric significance, making it natural to handle them differently.
\subsection{Generating uniformly distributed neurons for $\mathcal{P}_{\mathrm{NN}}$}
The second step in constructing $\mathcal{P}_{\mathrm{NN}}$ is to determine the location parameters $\left\{(\pmb{a}_m, r_m)\right\}_{m=1}^M$ in~\eqref{reparametrize} so that the resulting partition hyperplanes are uniformly distributed over the computational domain $\Omega$. 

The idea to handle this task is motivated by ReLU networks, whose expressive power is related to the number of linear regions being created in the input domain~\cite{Arora2016, Daubechies2022}. In addition to the number of linear pieces, their distribution also plays a crucial role, as a denser distribution in a specific region of the domain may enhance approximation quality for that region~\cite{Zhang2024TransNet}. Therefore, to achieve a strong and transferable approximation across all parts of the domain, the goal is to establish \textit{uniformly distributed} linear pieces throughout the domain. Based on the study in~\cite{Hanin2019}, the authors in~\cite{Zhang2024TransNet} proposed to quantify the expressive power of the ReLU networks through the density of the partition hyperplanes. 

For a detailed discussion of how to generate the parameters $\left\{(\pmb{a}_m, r_m)\right\}_{m=1}^M$ in this step, we refer the interested reader to~\cite{Zhang2024TransNet}. The complete procedure is to sample $\pmb{a}_m$ and $r_m$ as follows:
\begin{equation}
\label{step3_sample}
\pmb{a}_m = \dfrac{X_m}{\norm{X_m}_2}, \qquad r_m = U_m, \qquad m=1, \hdots, M,
\end{equation}
where $X_m$ are i.i.d. standard Gaussian random vectors and $U_m$ are i.i.d. uniform random variables on $[0,1]$.
\subsection{Turning the shape of the neurons in $\mathcal{P}_{\mathrm{NN}}$ using auxiliary functions}
The third step is to tune the shape parameters $\{\gamma_m\}_{m=1}^M$ in~\eqref{reparametrize}, which control the shapes of the neural basis functions $\psi_m(\pmb{y})$. For simplicity, we assume that all neurons share the same shape parameter, that is, $\gamma_m = \gamma$ for $m=1, \hdots, M$. The general idea is to construct a collection of auxiliary functions and then determine the value of $\gamma$ that minimizes the fitting error over this collection.

Since the goal is to construct a feature space $\mathcal{P}_{\mathrm{NN}}$ that can approximate a broad class of functions, the authors in~\cite{Zhang2024TransNet} propose using realizations of Gaussian random fields (GRFs) as auxiliary functions for tuning $\gamma$. To this end, let $G(\pmb{y} \mid \omega, \eta)$ denote the GRF used to generate the auxiliary functions, where $\omega$ represents the underlying randomness and $\eta$ is a fixed correlation length. We first generate $K$ realizations of the GRF, denoted by $\left\{G(\pmb{y} \mid \omega_k, \eta)\right\}_{k=1}^K$. For each realization, we evaluate its values
\[
g_j^k = G(\pmb{y}_j^k \mid \omega_k, \eta), \qquad j=1,\hdots,J,
\]
at $J$ sample points $\left\{\pmb{y}_j^k\right\}_{j=1}^J$ uniformly distributed in the unit ball $B_1(\pmb{0})$. The corresponding fitting error for the $k$th realization is defined by the least-squares mean squared error
\begin{equation}
\label{MSE_fitting_err}
\mathrm{MSE}_{k}(\gamma) = \min\limits_{\alpha} \sum\limits_{j=1}^J \left[g_j^k - \sum\limits_{m=0}^M \alpha_m \psi_m(\pmb{y}_j^k)\right]^2,
\end{equation}
where $\pmb{a}_m$ and $r_m$ are generated according to~\eqref{step3_sample} and $\psi_0(\cdot) \equiv 1$. The optimal shape parameter $\gamma^{\mathrm{opt}}$ is then defined as the minimizer of the average fitting error, namely,
\begin{equation}
\label{opt_fitting}
\gamma^{\mathrm{opt}} = \arg\min\limits_{\gamma} \frac{1}{K}\sum\limits_{k=1}^K \mathrm{MSE}_k(\gamma).
\end{equation}

The minimization problem~\eqref{opt_fitting} can be solved by grid search method since it is one dimensional. The main steps for solving $\gamma^{\mathrm{opt}}$ is summarized in Algorithm~\ref{alg:tuning_gamma}.
\begin{algorithm}[h]
\caption{Tuning shape parameter $\gamma$}
\label{alg:tuning_gamma}
\textbf{Input} Number of basis functions $M$; correlation length $\eta$;  number of GRF realizations $K$; grid search range $(\gamma^{\min}, \gamma^{\max})$ and mesh size $S$.\\
\textbf{Output} Optimal shape parameter $\gamma^{\mathrm{opt}}$.
\begin{algorithmic}[1]
\For{k=1:K}
\State Sample $\left\{y^{k}_j\right\}^J_{j=1}$ uniformly in the unit ball.
\State Evaluate a realization of GRF on $\left\{y^{k}_j\right\}^J_{j=1}$: $\left\{g^{k}_j = G(\pmb{y}_j \vert \omega_k, \eta)\right\}^J_{j=1}.$
\EndFor
\State Generate uniform mesh for $\gamma$: $\gamma^{\min} = \gamma_1 < \hdots < \gamma_S = \gamma^{\max}$.
\For{$s = 1:S$}
\For{$k = 1:K$}
\State Set $\psi_0(\pmb{y}) = 1$.
\For{$m=1:M$}
\State Generate location parameters $(\pmb{a}_m, r_m)$ according to Eq.~\eqref{reparametrize}.
\State Calculate weight and bias: \[\pmb{w}_m = \gamma_s \pmb{a}_m, \; b_m = \gamma_s r_m.\] \vspace{-0.5cm}
\State Construct the neural basis function: \[\psi_m(\pmb{y}) = \sigma\left(\pmb{w}^{\top}_m\pmb{y}+b_m\right).\]\vspace{-0.5cm}
\EndFor
\State Compute the least square MSE:
\begin{equation*}
\mathrm{MSE}_{s, k} = \min\limits_{\alpha} \sum\limits_{j=1}^J \left[g_j^k - \sum\limits_{m=0}^M \alpha_m \psi_m(\pmb{y}_j^k)\right]^2.
\end{equation*}
\EndFor
\State Compute the average MSE: \[\mathrm{avgMSE}_s = \dfrac{1}{K}\sum^K_{k=1}\mathrm{MSE}_{s, k}.\] \vspace{-0.5cm}
\EndFor
\State Find the smallest $\mathrm{avgMSE}_s$: 
\[s^{*} = \arg \min^S_{s=1} \mathrm{avgMSE}_s.\]
\Return $\gamma^{\mathrm{opt}} = \gamma_{s^{*}}.$
\end{algorithmic}
\end{algorithm}

\subsection{Applying TransNet to approximate a function}
Assume that the feature space $\mathcal{P}_{\mathrm{NN}}$ is well-tuned. Let $f$ be the target function that we aim to approximate and let $\{\pmb{y}_j\}^J_{j=1}$ denote the samples from $\Omega \cup \partial\Omega$, the approximation of $f$ from the space $\mathcal{P}_{\mathrm{NN}}$ is obtained by solving the following optimization problem
\begin{equation}
\label{step4_optimize}
\min\limits_{\alpha}\left\{ \dfrac{1}{J}\sum\limits^J_{j=1} \left[f\left(\pmb{y}_j\right) - \sum\limits^M_{m=0}\alpha_m \psi_m\left(\pmb{y}_j\right)\right]\right\}.
\end{equation}
We denote by $\mathrm{TN}_f$ the approximate function to $f$ generated by solving the optimization~\eqref{step4_optimize}. We shall provide the details of the candidate pool for each test case in the numerical section.
\section{Numerical Experiments}
\label{Sec4_Results}
In this section, we present three high-dimensional time-independent PDEs to demonstrate the performance of the proposed algorithm in recovering the corresponding true solutions. In particular, we consider the Poisson equation and the reaction-diffusion equation in 60 dimensions, followed by a semi-linear elliptic equation in 55 dimensions. All equations are equipped with Dirichlet boundary conditions. For the last test case, we slightly reduce the dimension in order to increase the number of training points used in the searching loop.

We use the following loss functional to compute the score in \eqref{practice_score} for FEX:
\begin{equation}
\label{MSE_functional}
\mathcal{L}(u):= \norm{\mathcal{D}u(x) - f(u(x), x)}^2_{L^2(\Omega)} + \lambda \norm{\mathcal{B}u(x) - g(x)}^2_{L^2(\partial{\Omega})},
\end{equation}
where $\lambda$ is a positive coefficient used to enforce the boundary constraint. We choose $\lambda = 100$ for all three test cases and use 1000 boundary points to approximate the second integral in \eqref{MSE_functional}. For the first integral in \eqref{MSE_functional}, we use 500 interior points in the first two cases and 600 interior points in the last case.
\subsection{Test case 1: 60D Poisson equation}
We first consider the Poisson equation
\begin{equation}
\label{60D_Poisson}
\begin{array}{l}
-\Delta u(\pmb{x}) = f(\pmb{x}), \quad \pmb{x} \in \Omega,
\end{array}
\end{equation}
where $\pmb{x} = (x_1, \hdots, x_d) \in \Omega = (-1, 1)^d$ with $d = 60$. The true solution to~\eqref{60D_Poisson} is chosen as
\[
u(\pmb{x}) = 0.5\sum_{i=1}^d x_i^2.
\]
We study two scenarios: one in which the candidate pool contains the candidate associated with the key nonlinear operator $x^2$ appearing in the true solution, and one in which it does not. More precisely, we consider the following two candidate pools for the FEX search loop:
\begin{equation}
\label{Poisson_pool1}
\mathbb{P}_1 := \{0, 1, \mathrm{Id}, \mathrm{TN}_{x^2}, \mathrm{TN}_{x^3}, \mathrm{TN}_{x^4}, \mathrm{TN}_{\exp}, \mathrm{TN}_{\sin(x)}, \mathrm{TN}_{\cos(x)} \},
\end{equation}
and
\begin{equation}
\label{Poisson_pool2}
\mathbb{P}_2 := \{0, 1, \mathrm{Id}, \mathrm{TN}_{\sin(x^2)}, \mathrm{TN}_{x^3}, \mathrm{TN}_{x^4}, \mathrm{TN}_{\exp}, \mathrm{TN}_{\sin(x)}, \mathrm{TN}_{\cos(x)} \}.
\end{equation}
In the second candidate pool, we exclude the candidate $\mathrm{TN}_{x^2}$. To keep the two candidate pools the same size, we replace it with $\mathrm{TN}_{\sin(x^2)}$, which approximates $\sin(x^2)$. Since $\sin(x^2)$ behaves similarly to $x^2$ near $x=0$, this replacement preserves similar local nonlinear behavior while removing the candidate associated with the true quadratic structure.

We run the search loop (that is, the controller update) for $700$ iterations for the pool $\mathbb{P}_1$ and for $1200$ iterations for the pool $\mathbb{P}_2$. Due to the absence of the true operator, we allow more iterations for the second pool so that FEX can identify a suitable replacement. For both pools, the relevant score-related parameters are updated first using Adam with learning rate $0.001$ for $T_1 = 2$ iterations and then using BFGS with learning rate $1$ for at most $T_2 = 20$ iterations. Finally, the parameter vector $\theta$ is optimized using Adam with initial learning rate $0.01$ for $T_3 = 15000$ iterations. We fix the depth of the binary tree to be $2$.

\begin{figure}[h!]
\begin{minipage}{0.333\textwidth} 
\includegraphics[scale= 0.185]{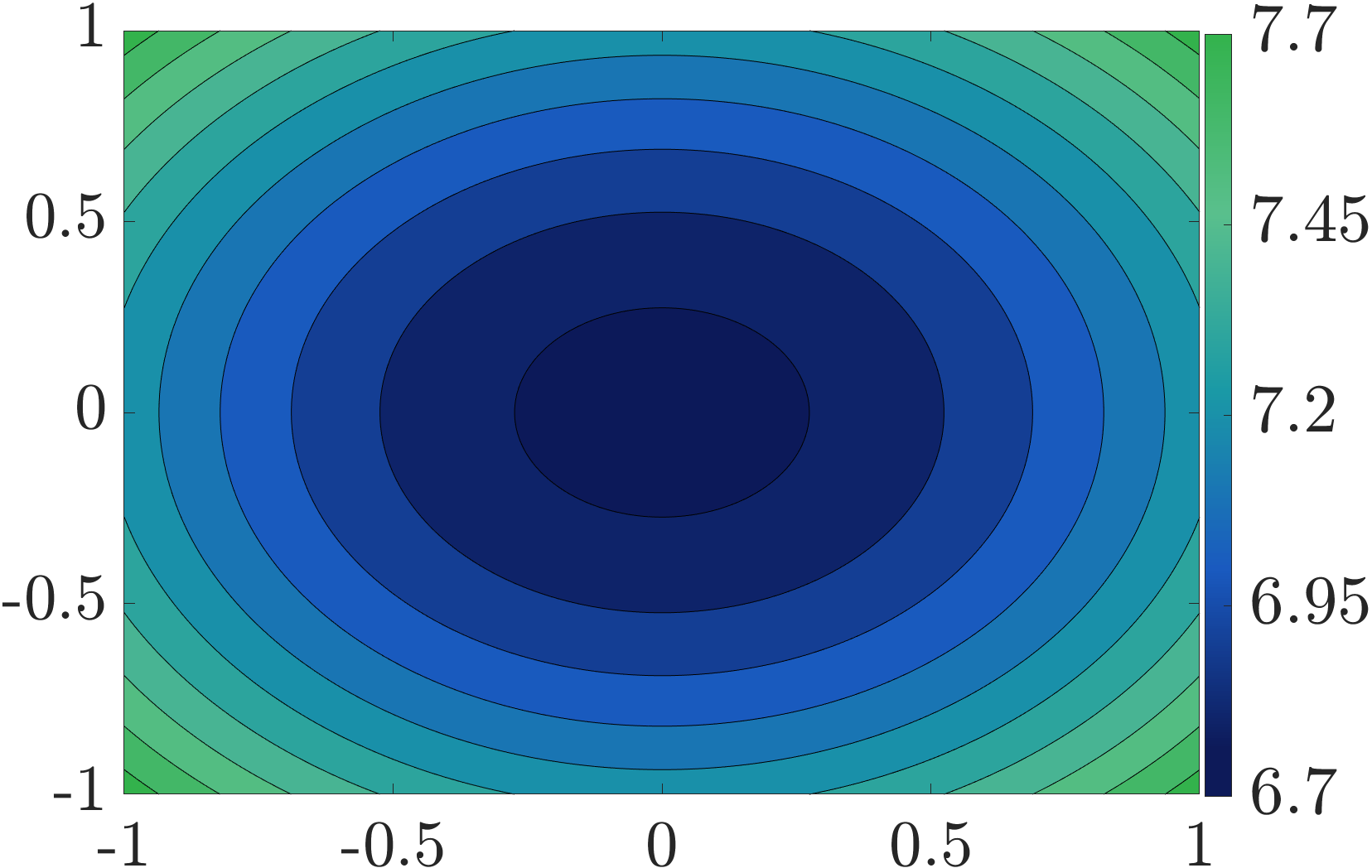}
\end{minipage}%
\begin{minipage}{0.333\textwidth}
\includegraphics[scale= 0.185]{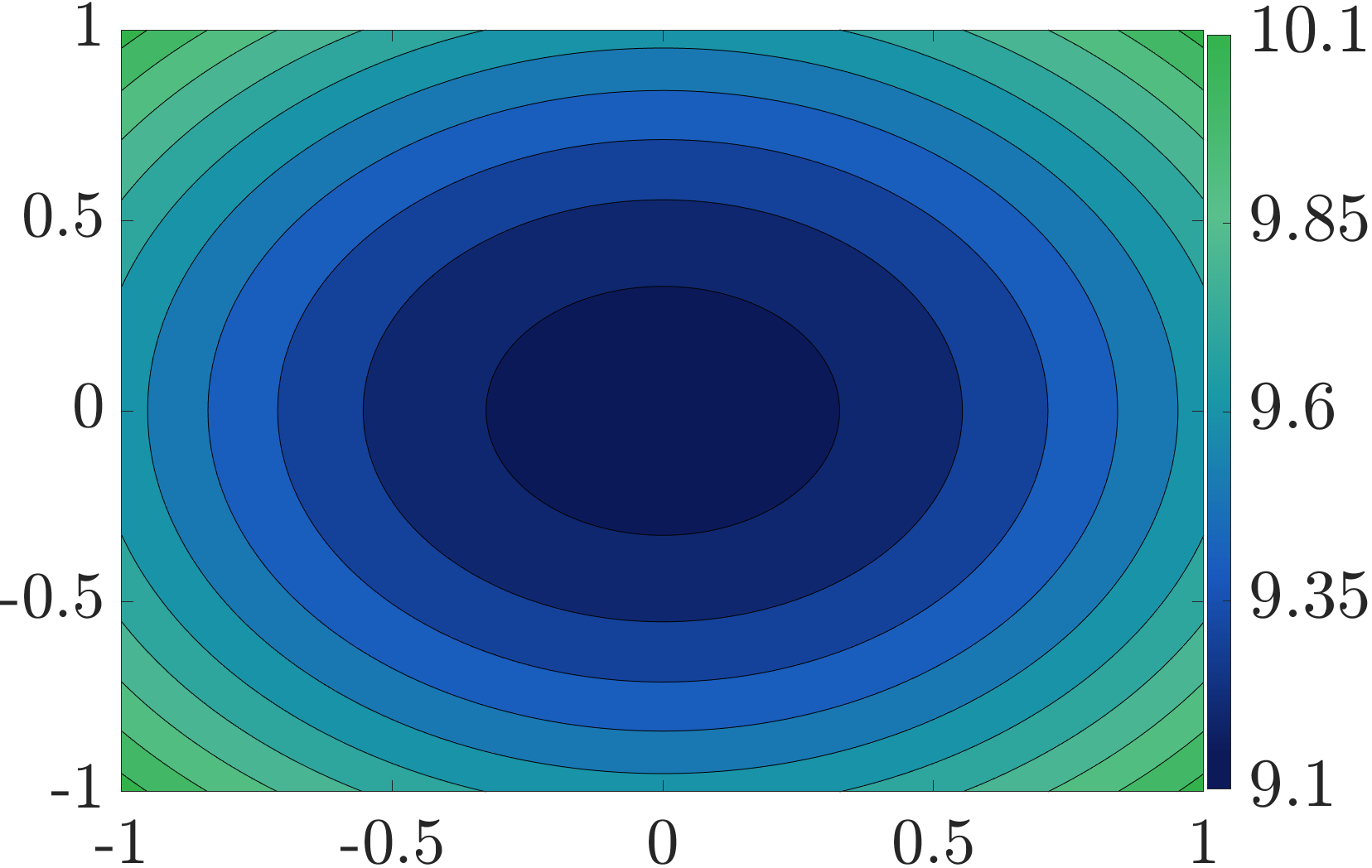}
\end{minipage}%
\begin{minipage}{0.333\textwidth} 
\includegraphics[scale= 0.185]{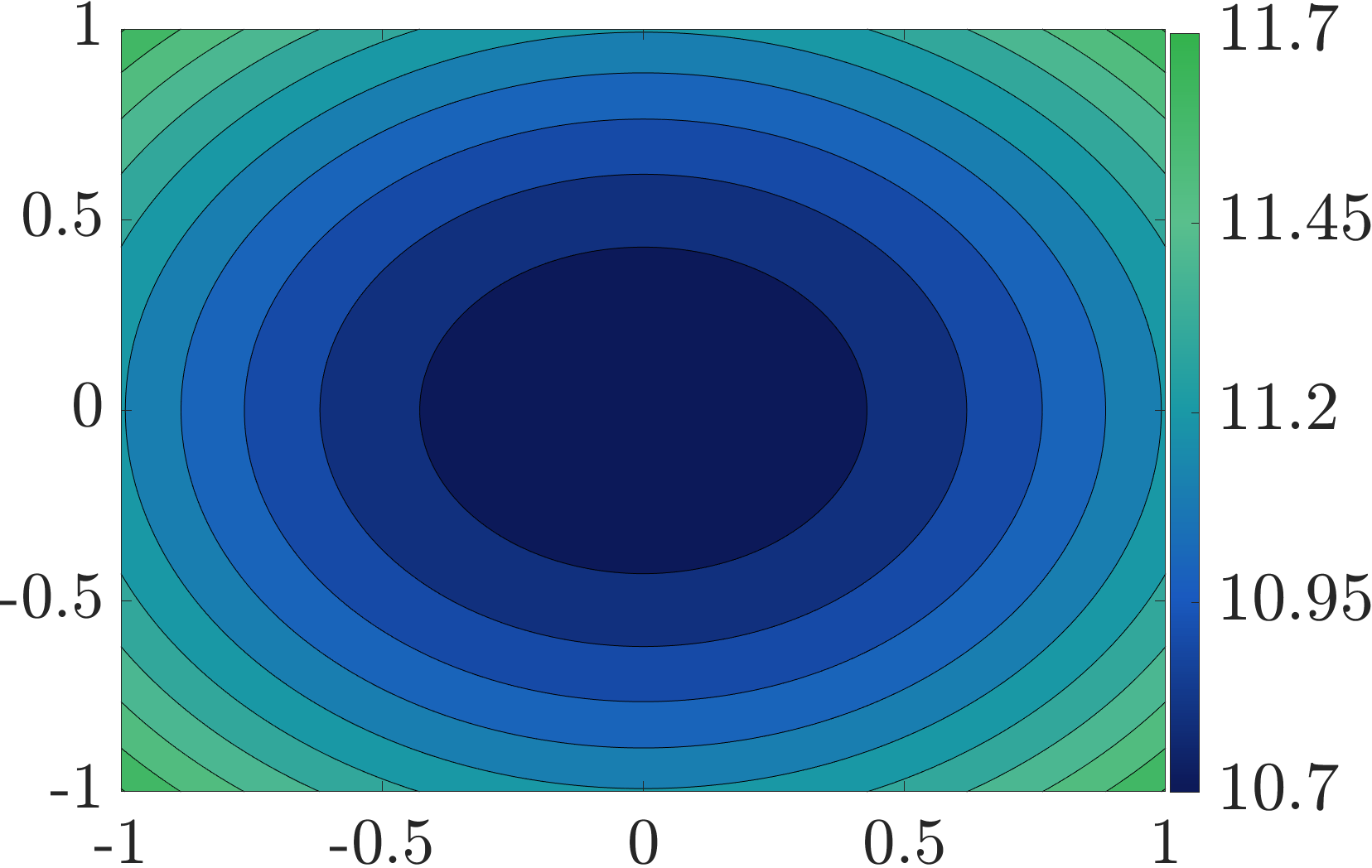}
\end{minipage}
\caption{\footnotesize [60D Poisson] Heatmaps of the reference solution on two-dimensional slices, with the remaining $58$ dimensions fixed at predefined values. (First) Dimensions $(22, 37)$. (Second) Dimensions $(30, 35)$. (Third) Dimensions $(41, 18)$.}
\label{Poisson_Ref}
\vspace{-0.4cm}
\end{figure}
\begin{figure}[h!]
\begin{minipage}{0.333\textwidth} 
\includegraphics[scale= 0.185]{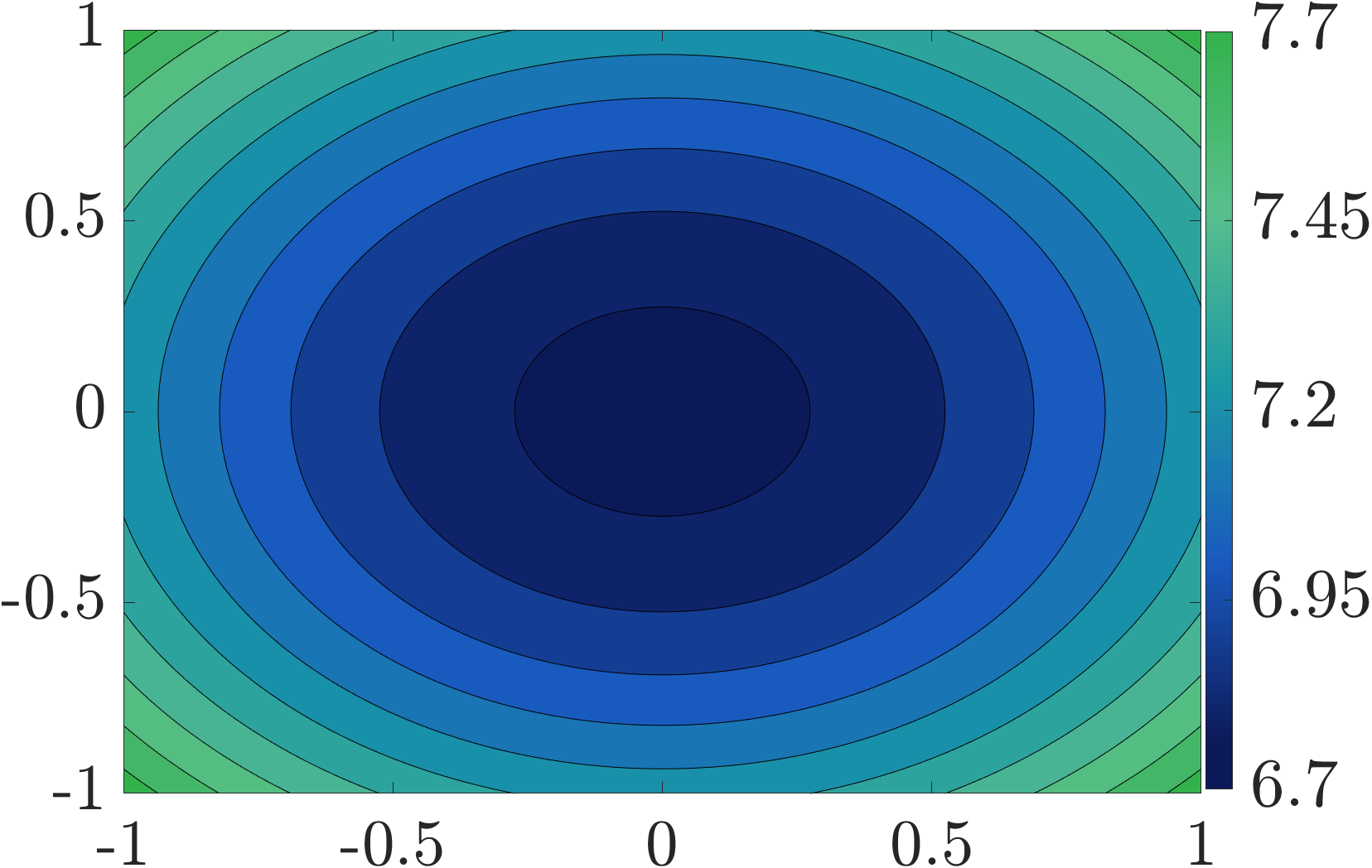}
\end{minipage}%
\begin{minipage}{0.333\textwidth}
\includegraphics[scale= 0.185]{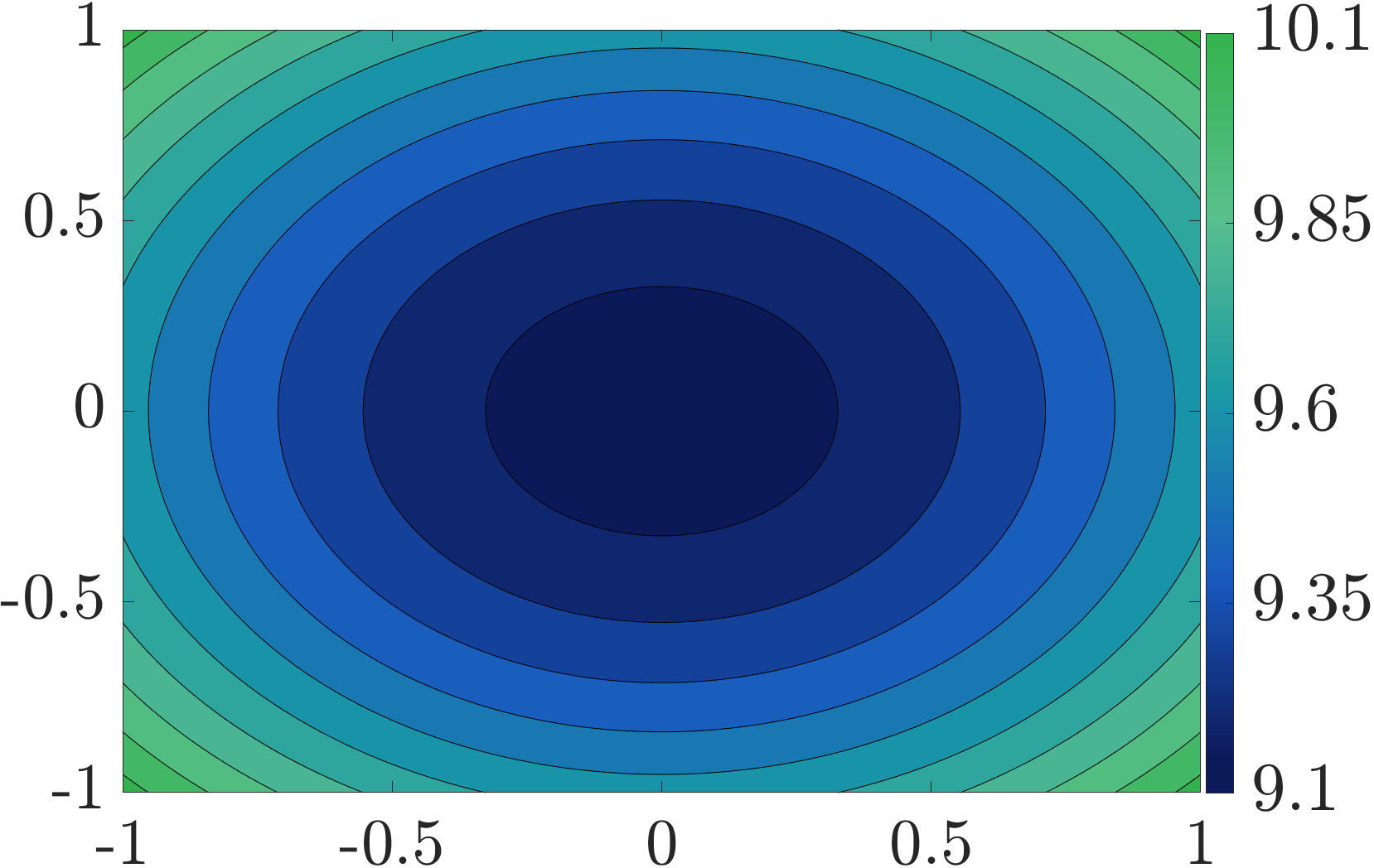}
\end{minipage}%
\begin{minipage}{0.333\textwidth} 
\includegraphics[scale= 0.185]{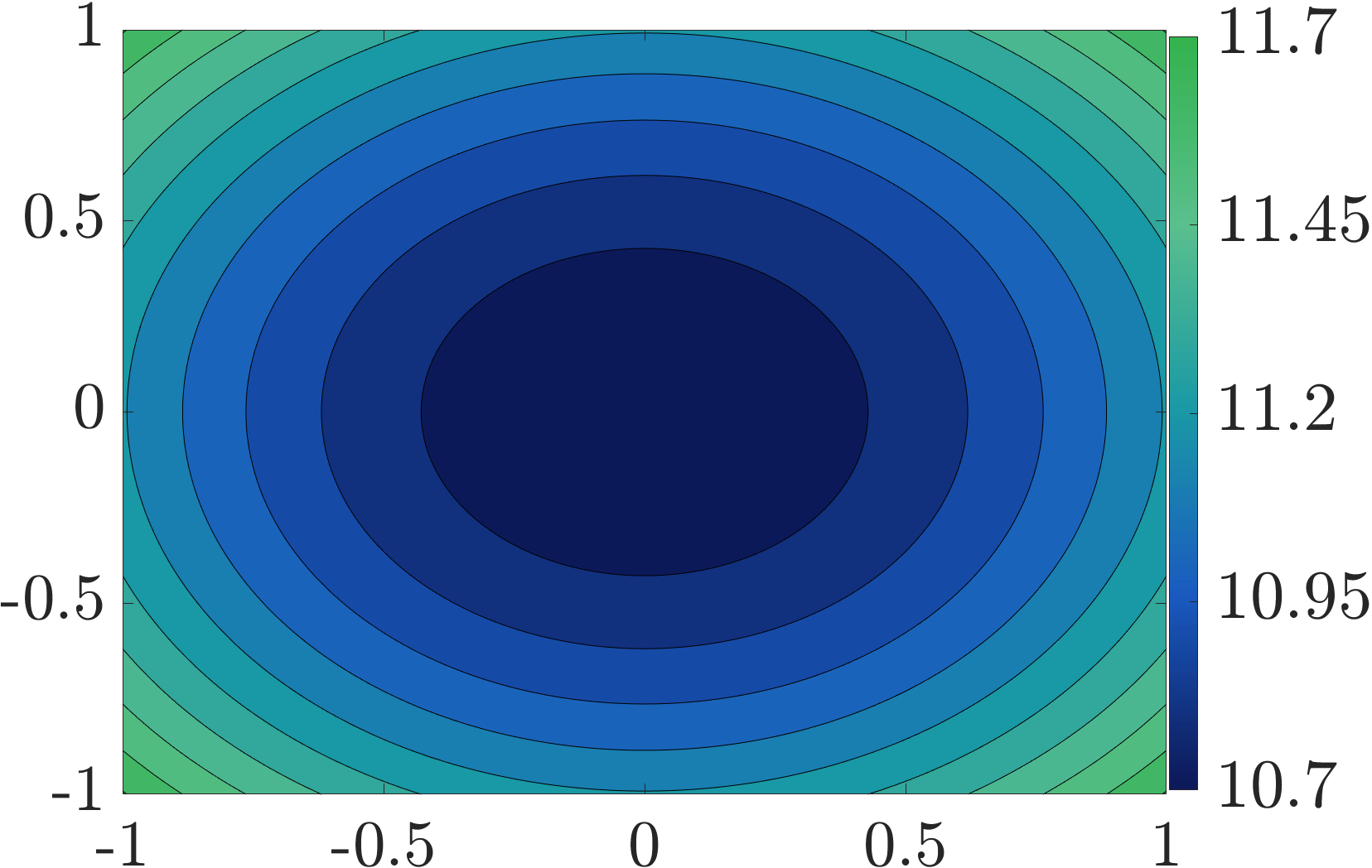}
\end{minipage}

\begin{minipage}{0.333\textwidth} 
\includegraphics[scale= 0.185]{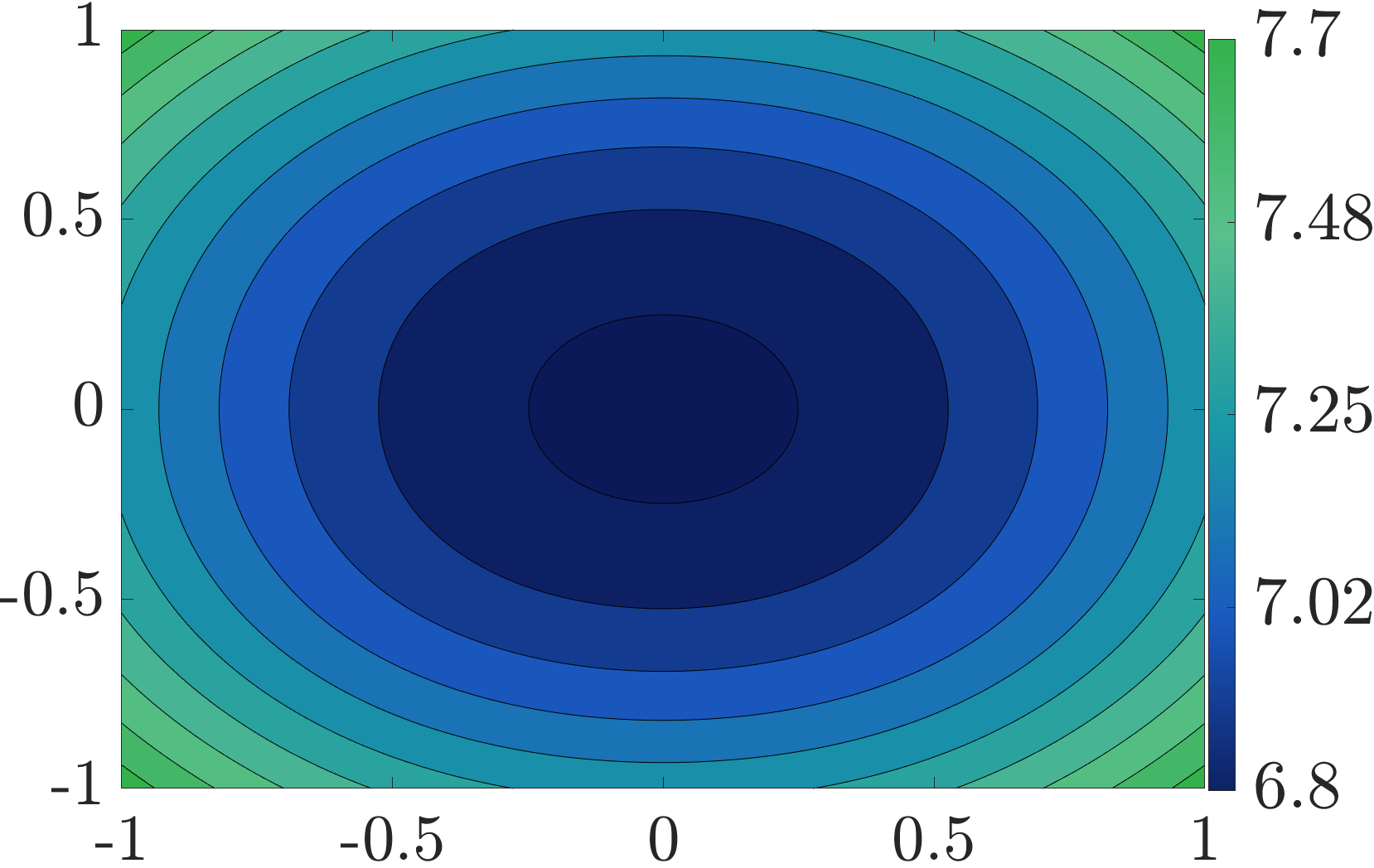}
\end{minipage}%
\begin{minipage}{0.333\textwidth}
\includegraphics[scale= 0.185]{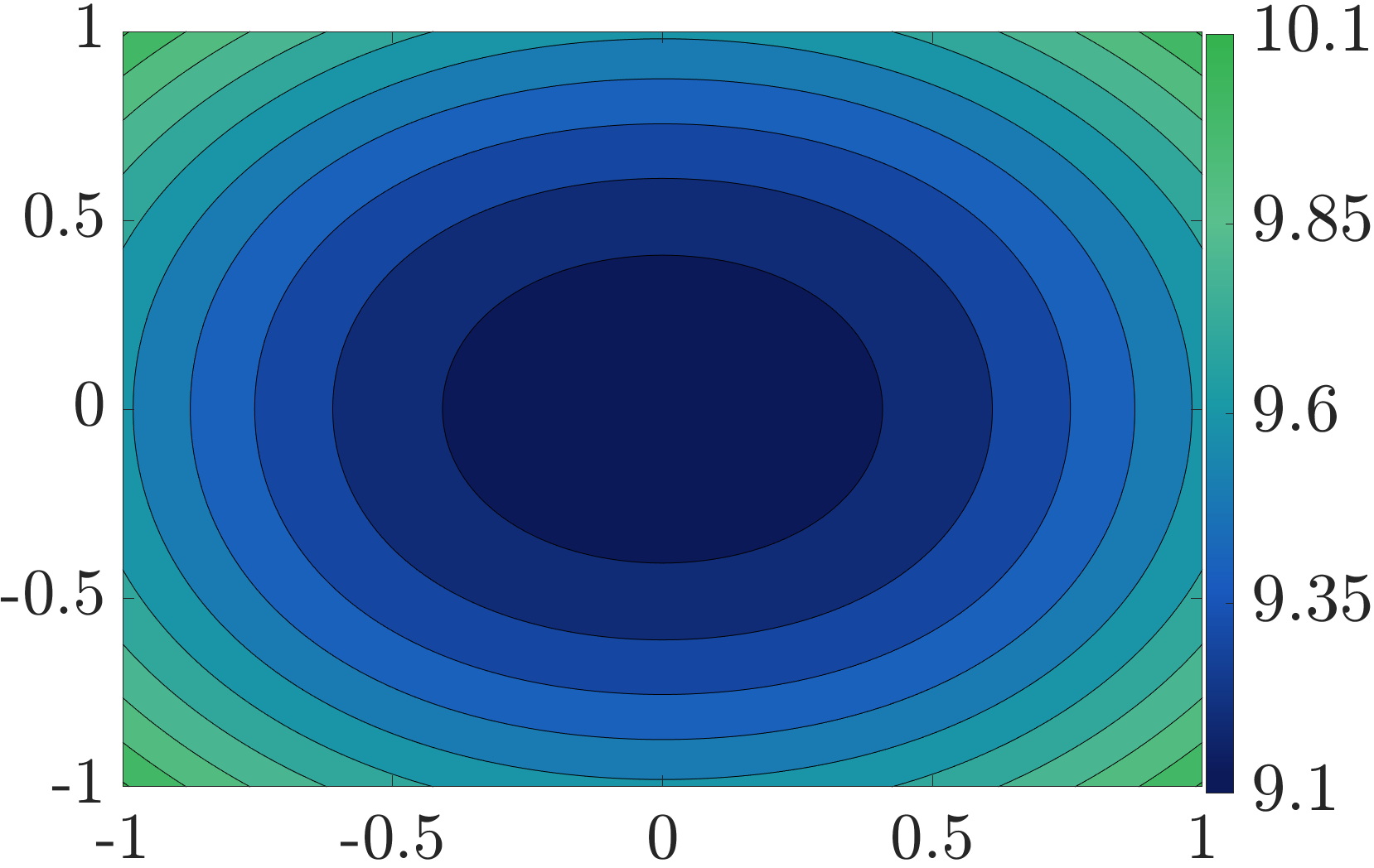}
\end{minipage}%
\begin{minipage}{0.333\textwidth} 
\includegraphics[scale= 0.185]{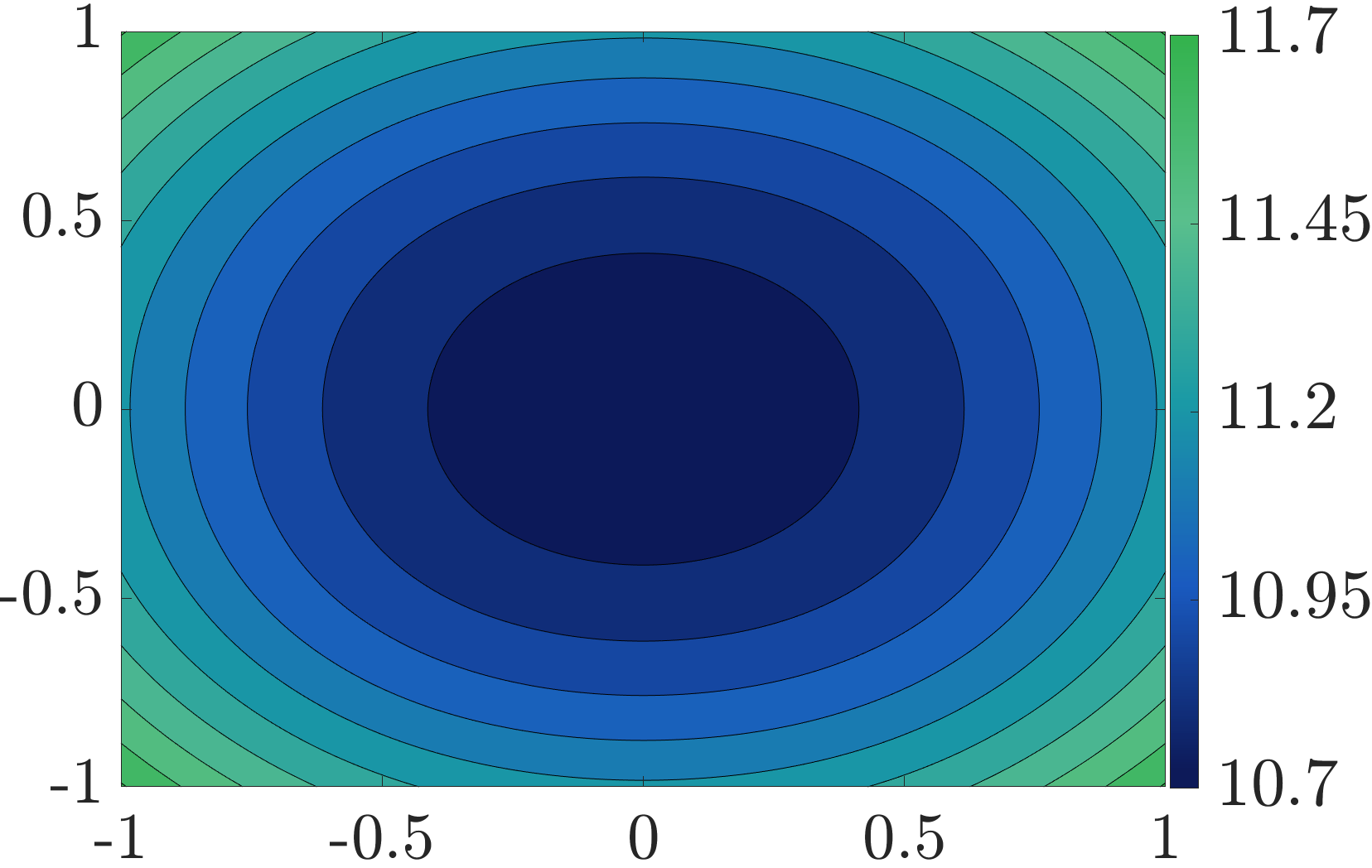}
\end{minipage}
\caption{\footnotesize [60D Poisson] Heatmaps of the predicted solution with pool $\mathbb{P}_1$ (first row) and pool $\mathbb{P}_2$ (second row). (First column) Dimensions $(22, 37)$. (Second column) Dimensions $(30, 35)$. (Third column) Dimensions $(41, 18)$.}
\label{Poisson_Pred}
\end{figure}

\begin{figure}[h!]
\begin{minipage}{0.333\textwidth} 
\includegraphics[scale= 0.185]{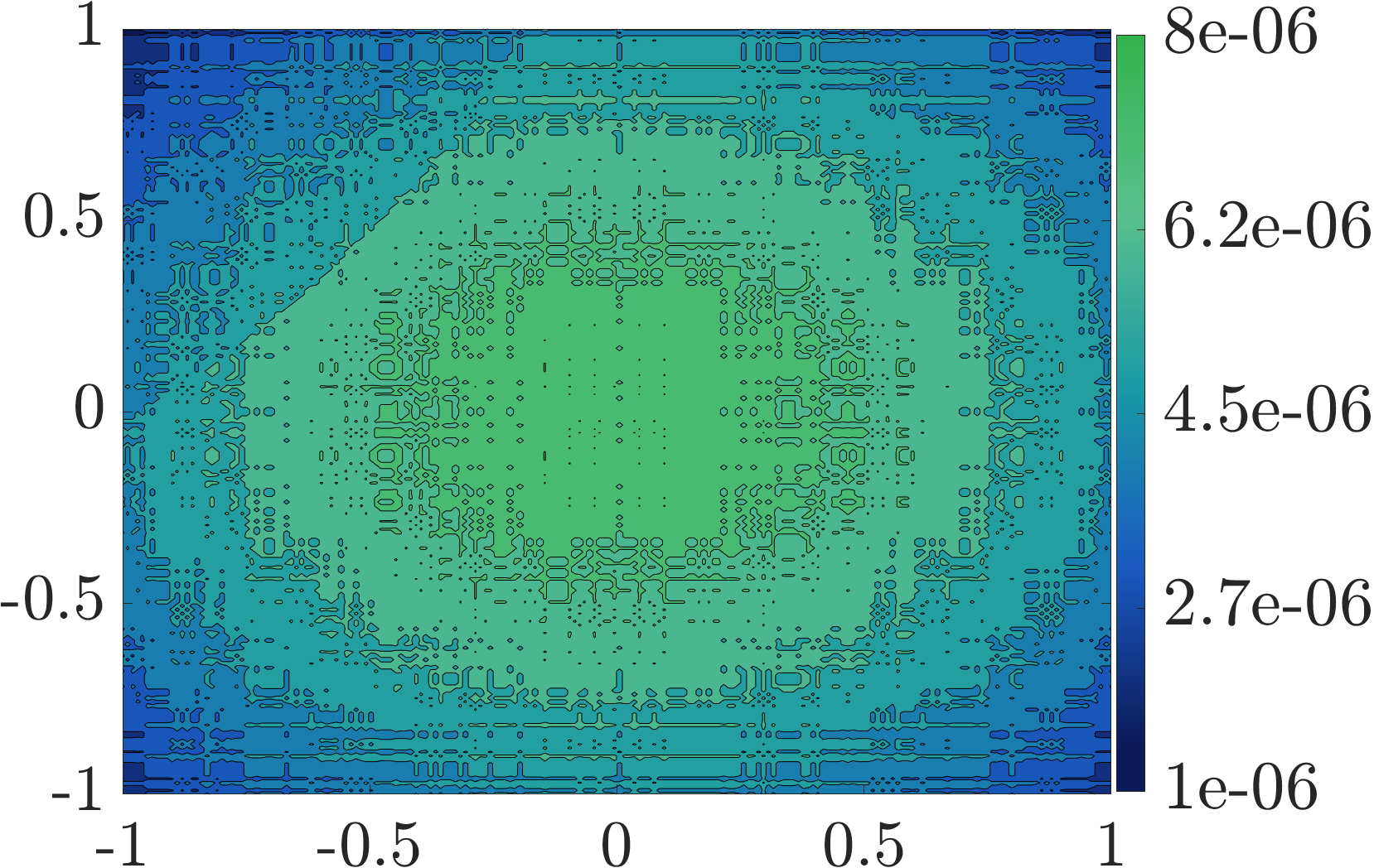}
\end{minipage}%
\begin{minipage}{0.333\textwidth}
\includegraphics[scale= 0.185]{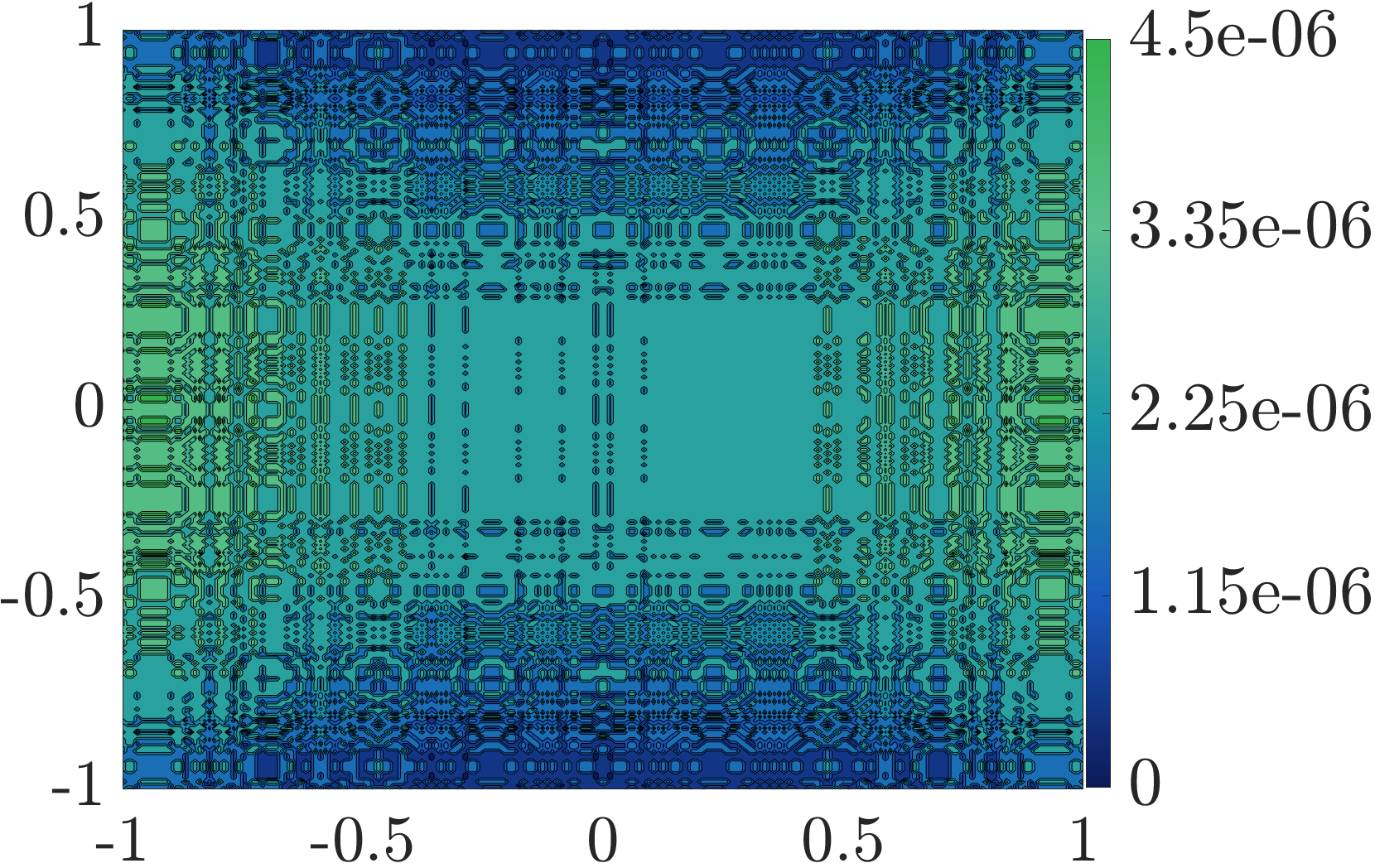}
\end{minipage}%
\begin{minipage}{0.333\textwidth} 
\includegraphics[scale= 0.185]{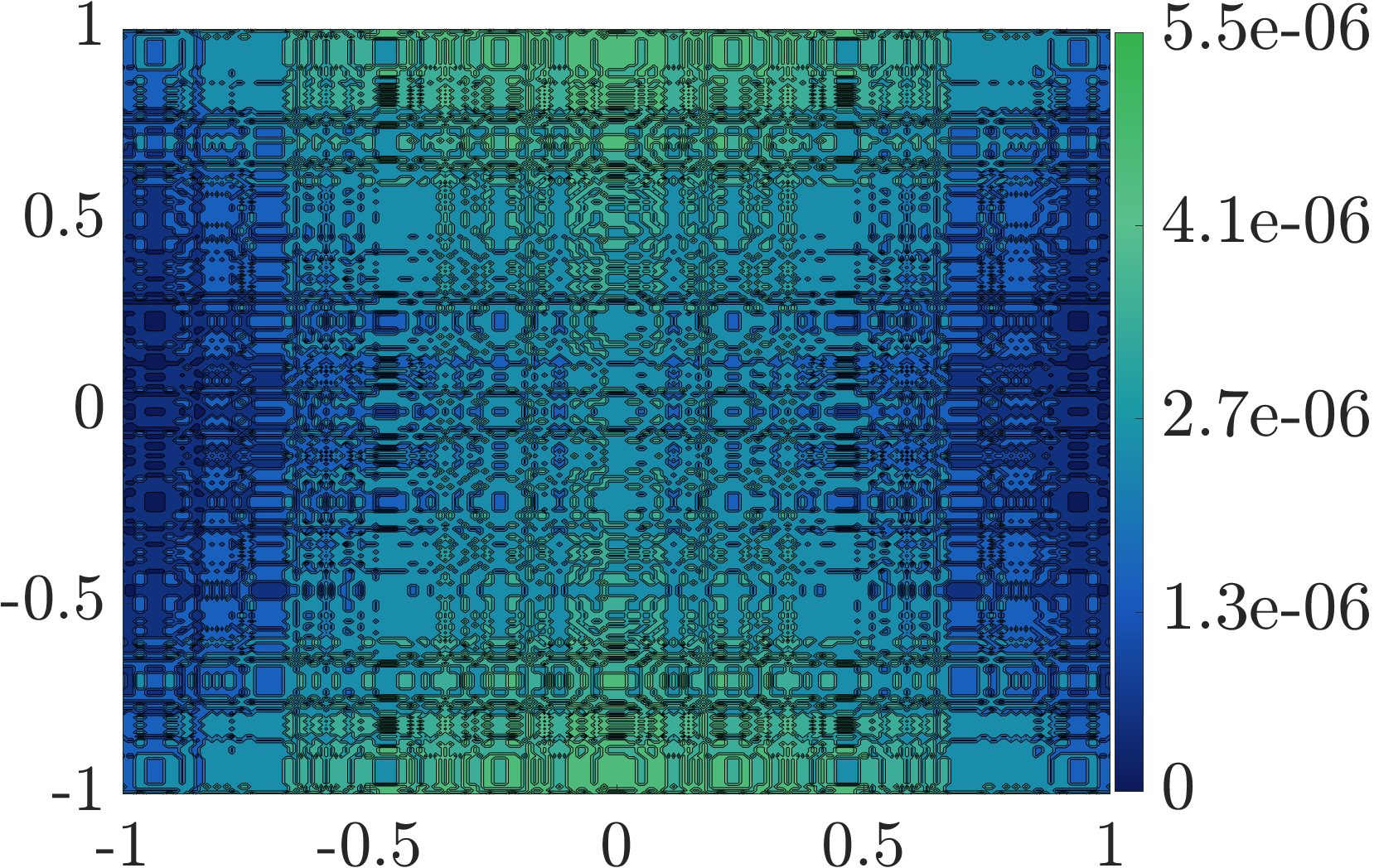}
\end{minipage}

\begin{minipage}{0.333\textwidth} 
\includegraphics[scale= 0.185]{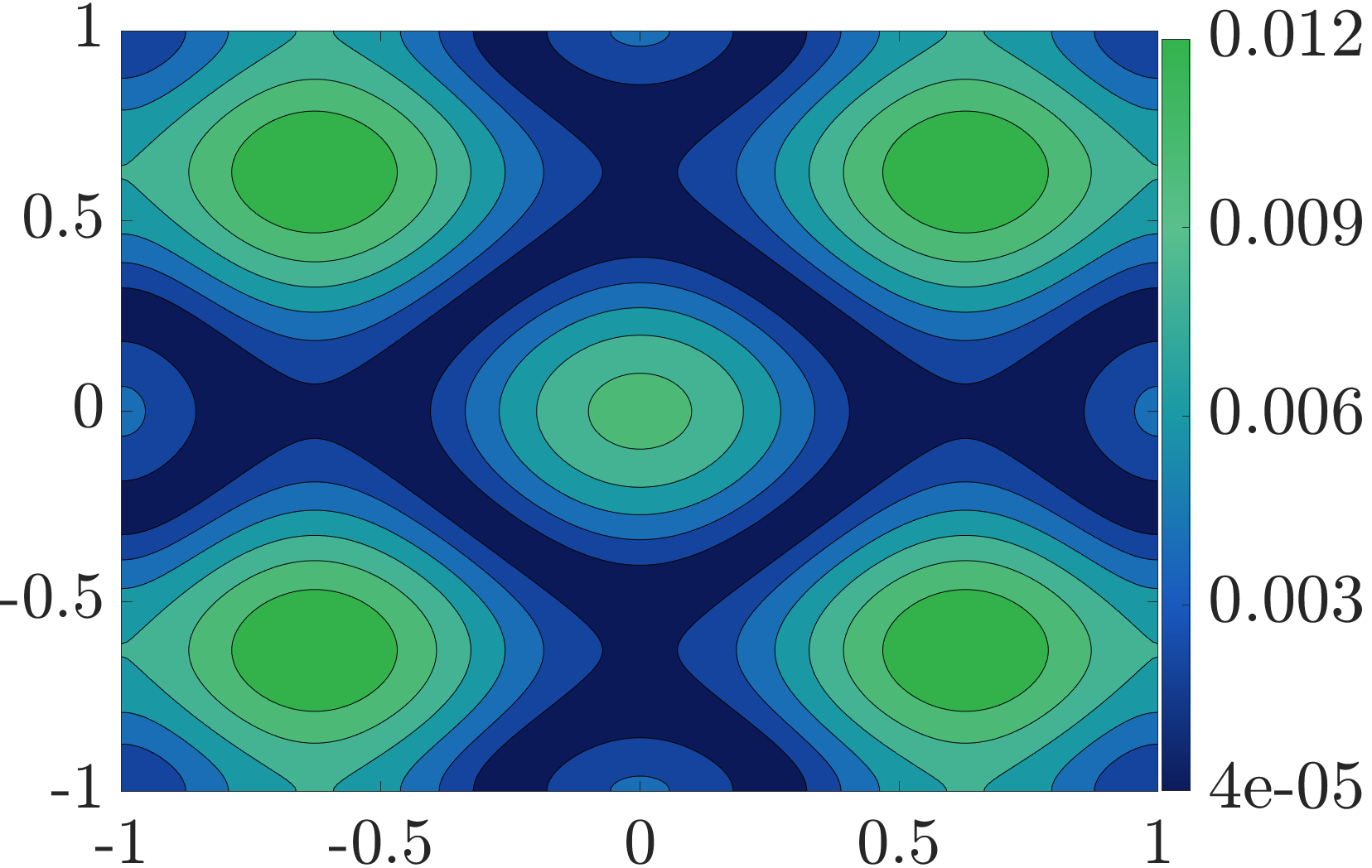}
\end{minipage}%
\begin{minipage}{0.333\textwidth}
\includegraphics[scale= 0.185]{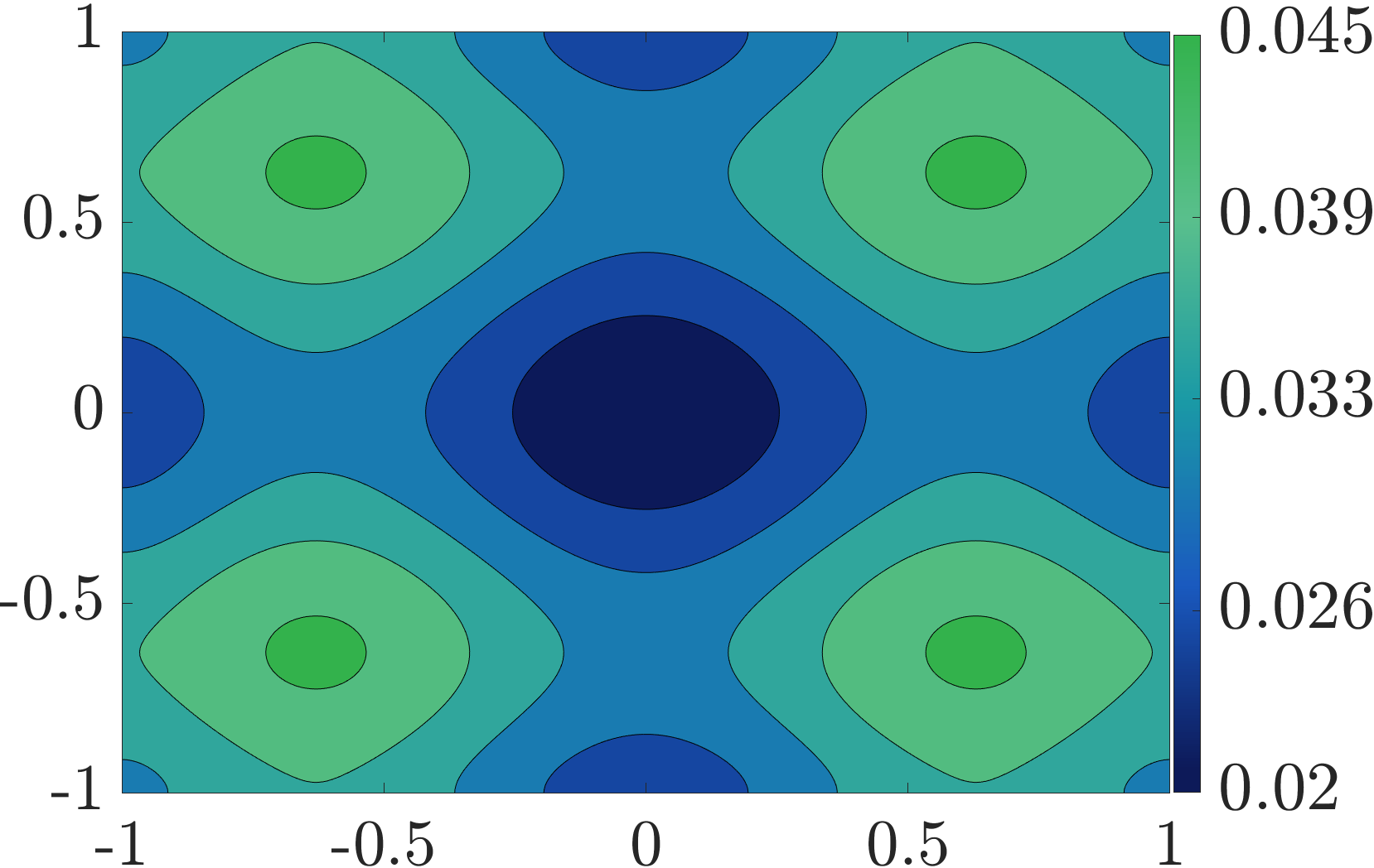}
\end{minipage}%
\begin{minipage}{0.333\textwidth} 
\includegraphics[scale= 0.185]{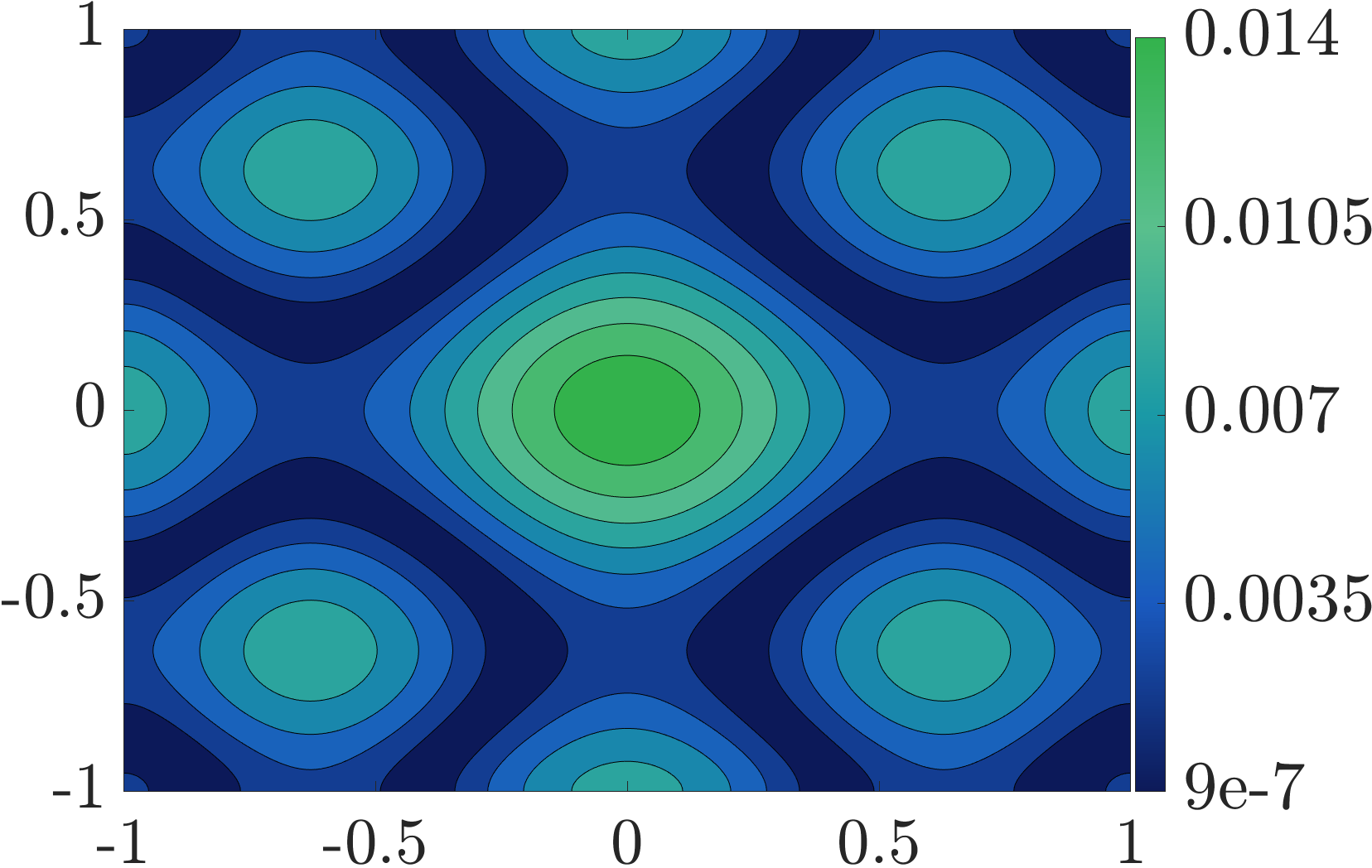}
\end{minipage}
\caption{\footnotesize [60D Poisson] Heatmaps of absolute pointwise error with pool $\mathbb{P}_1$ (first row) and pool $\mathbb{P}_2$ (second row). (First column) Dimensions $(22, 37)$. (Second column) Dimensions $(30, 35)$. (Third column) Dimensions $(41, 18)$.}
\label{Poisson_Error}
\vspace{-0.4cm}
\end{figure}
The true solution on three selected pairs of dimensions, with resolution $200 \times 200$, is displayed in Figure~\ref{Poisson_Ref}. The corresponding predicted solutions generated from the pools $\mathbb{P}_1$ and $\mathbb{P}_2$ are shown in Figure~\ref{Poisson_Pred}. After the search loop terminates, FEX constructs an expression involving the candidates associated with $(\cdot)^2$ and $\sin(\cdot)$ when trained with the pool $\mathbb{P}_1$. When the candidate associated with $(\cdot)^2$ is excluded, FEX instead constructs an expression using a linear combination of the candidates associated with $(\cdot)^4$ and $\sin((\cdot)^2)$. In both cases, the coefficients in the resulting expressions are further optimized so that the final approximation better matches the true solution. As seen from the figures, FEX accurately recovers the true solution when the candidate associated with $(\cdot)^2$ is included in the pool, since this allows the model to represent the solution more directly and makes the optimization of $\pmb{\theta}$ easier. When this candidate is excluded, FEX must instead rely on the candidates associated with $(\cdot)^4$ and $\sin((\cdot)^2)$ to mimic the role of $(\cdot)^2$, which makes it more difficult to tune the corresponding coefficients so that their combination matches the true solution. As a result, FEX still captures the overall shape of the solution on these dimension pairs, although a noticeable discrepancy remains.

% The true solution on three selected pairs of dimensions, with a resolution of $200 \times 200$, is displayed in Figure~\ref{Poisson_Ref}. The corresponding predicted solutions generated from the pools $\mathbb{P}_1$ and $\mathbb{P}_2$ are shown in Figure~\ref{Poisson_Pred}. After the search loop terminates, FEX constructs an expression involving $(\cdot)^2$ and $\sin(\cdot)$ when trained with the pool $\mathbb{P}_1$. When the operator $(\cdot)^2$ is excluded, FEX instead builds an expression using a linear combination of $(\cdot)^4$ and $\sin((\cdot)^2)$, since these functions have shapes similar to that of $(\cdot)^2$. In both cases, the resulting candidates are further optimized to adjust their magnitudes so that they better match the true solution. As seen in both figures, FEX accurately recovers the true solution when $(\cdot)^2$ is included in the candidate pool, since this allows the model to represent the solution more directly and makes the optimization of $\pmb{\theta}$ easier. When this operator is excluded, FEX must instead rely on $(\cdot)^4$ and $\sin((\cdot)^2)$ to mimic the role of $(\cdot)^2$, which makes it more difficult to adjust the corresponding coefficients so that their combination matches the true solution. As a result, FEX still captures the overall shape of the solution on these dimension pairs, although a noticeable discrepancy remains.

To further compare the two approximate solutions, we plot the corresponding absolute pointwise errors in Figure~\ref{Poisson_Error}. The results show that FEX reconstructs the true solution with very high accuracy when the candidate pool contains the true operator. The predicted solution generated from pool $\mathbb{P}_2$ is less accurate than that from pool $\mathbb{P}_1$, but the error remains moderate. This behavior is consistent with the numerical results reported in~\cite{Liang2025finite}.

Finally, we compute the relative $L^2$ errors for both approximate solutions. Since computing these errors requires evaluating high-dimensional integrals, we approximate them using Monte Carlo integration. Because the samples used in Monte Carlo integration are randomly generated, we repeat the computation $50$ times and report the sample means and standard deviations in Table~\ref{L2_Errors}. These results further confirm the performance of FEX for both candidate pools.
\subsection{Test case 2: 60D Stationary reaction-diffusion equation}
We next consider the stationary reaction--diffusion equation
\begin{equation}
\label{60D_ReactDiff}
-\nu \Delta u + \mu u = f(\pmb{x}), \quad \pmb{x} \in \Omega,
\end{equation}
where $\pmb{x} = (x_1, \hdots, x_d) \in \Omega = (0,1)^d$ with $d = 60$. We assume that the true solution to~\eqref{60D_ReactDiff} is given by
\[
u(\pmb{x}) = \sum_{i=1}^d x_i^3.
\]
Similar to the previous test case, we consider two candidate pools, only one of which contains the candidate $\mathrm{TN}_{x^3}$. More precisely, the first candidate pool is defined as
\begin{equation}
\label{ReactDiff_pool1}
\mathbb{P}_1 := \{0, 1, \mathrm{Id}, \mathrm{TN}_{x^2}, \mathrm{TN}_{x^3}, \mathrm{TN}_{x^4}, \mathrm{TN}_{\exp}, \mathrm{TN}_{\sin(x)}, \mathrm{TN}_{\cos(x)} \}.
\end{equation}
For the second candidate pool, we remove the candidate $\mathrm{TN}_{x^3}$. To keep the two candidate pools the same size, we replace it with $\mathrm{TN}_{x\sin(x)}$, which approximates $x\sin(x)$. Unlike the replacement used in the previous test case, $x\sin(x)$ does not exhibit local behavior similar to that of $x^3$. Thus, the second candidate pool is defined as
\begin{equation}
\label{ReactDiff_pool2}
\mathbb{P}_2 := \{0, 1, \mathrm{Id}, \mathrm{TN}_{x^2}, \mathrm{TN}_{x\sin(x)}, \mathrm{TN}_{x^4}, \mathrm{TN}_{\exp}, \mathrm{TN}_{\sin(x)}, \mathrm{TN}_{\cos(x)} \}.
\end{equation}
We use the same numbers of search iterations for both candidate pools as in the previous test case.

\begin{figure}[h!]
\begin{minipage}{0.333\textwidth} 
\includegraphics[scale= 0.185]{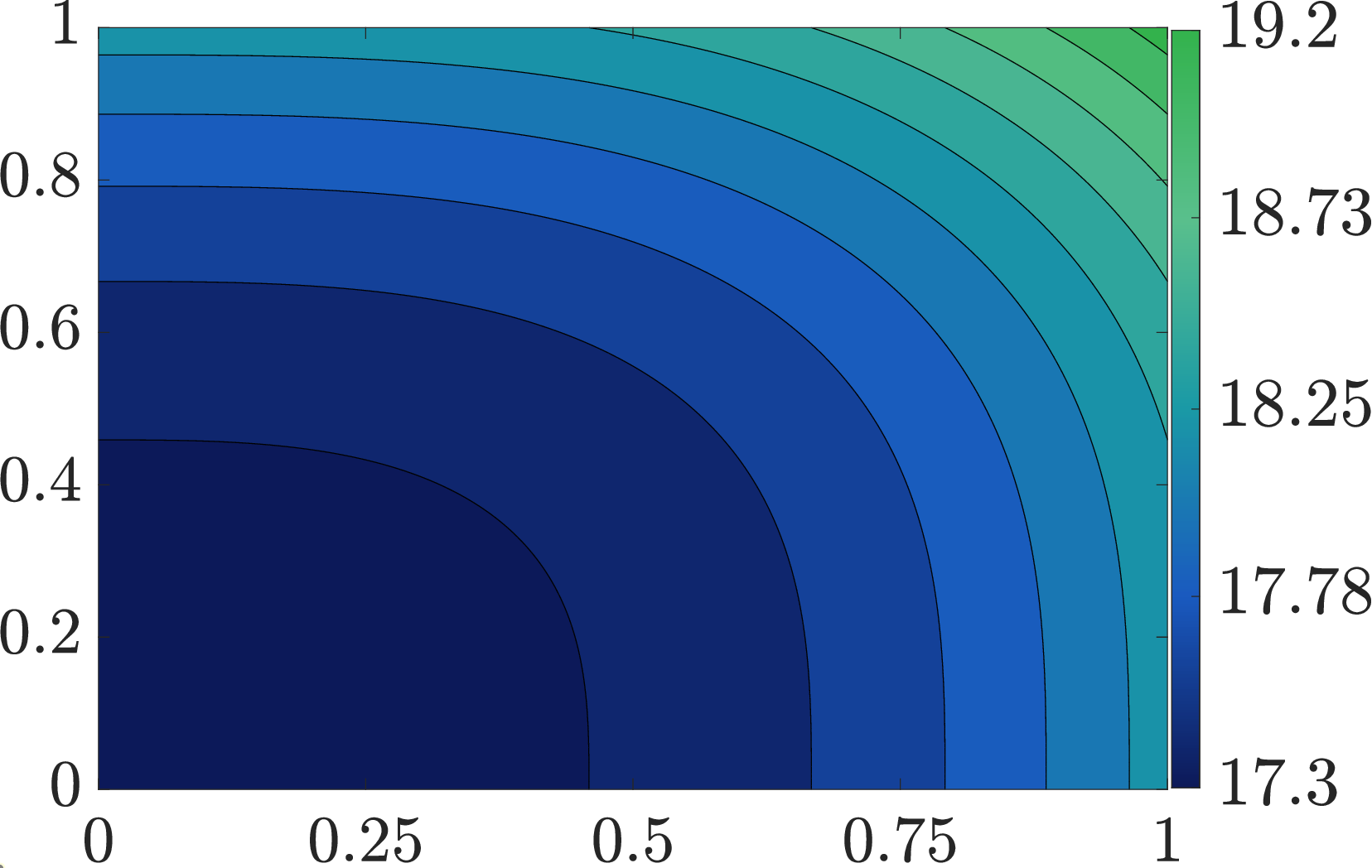}
\end{minipage}%
\begin{minipage}{0.333\textwidth}
\includegraphics[scale= 0.185]{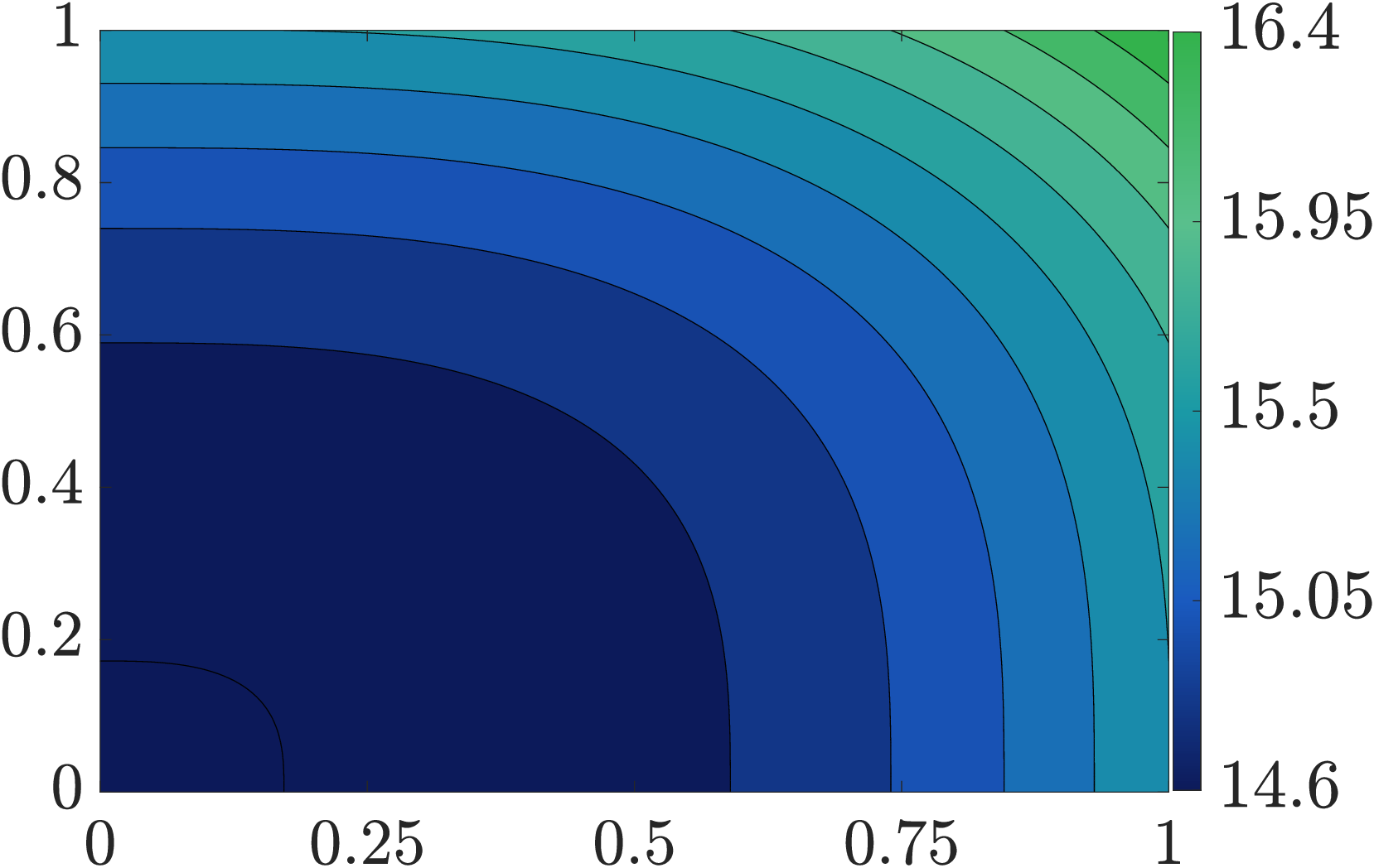}
\end{minipage}%
\begin{minipage}{0.333\textwidth} 
\includegraphics[scale= 0.185]{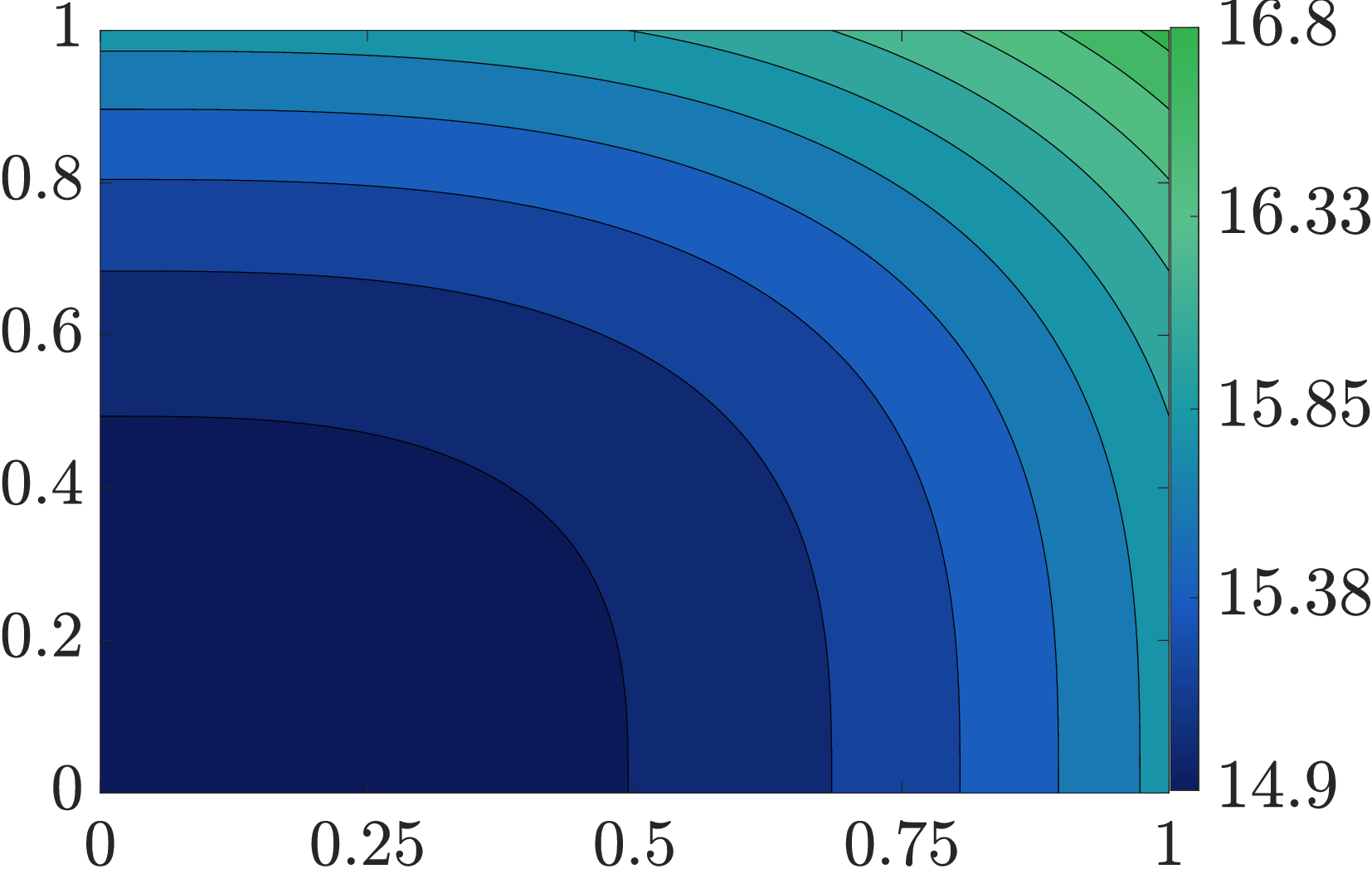}
\end{minipage}
\caption{\footnotesize [60D Reaction-diffusion] Heatmaps of the reference solution on two-dimensional slices, with the remaining $58$ dimensions fixed at predefined values. (First) Dimensions $(17, 33)$. (Second) Dimensions $(21, 56)$. (Third) Dimensions $(52, 19)$.}
\label{ReactDiff_Ref}
\end{figure}

\begin{figure}[h!]
\begin{minipage}{0.333\textwidth} 
\includegraphics[scale= 0.185]{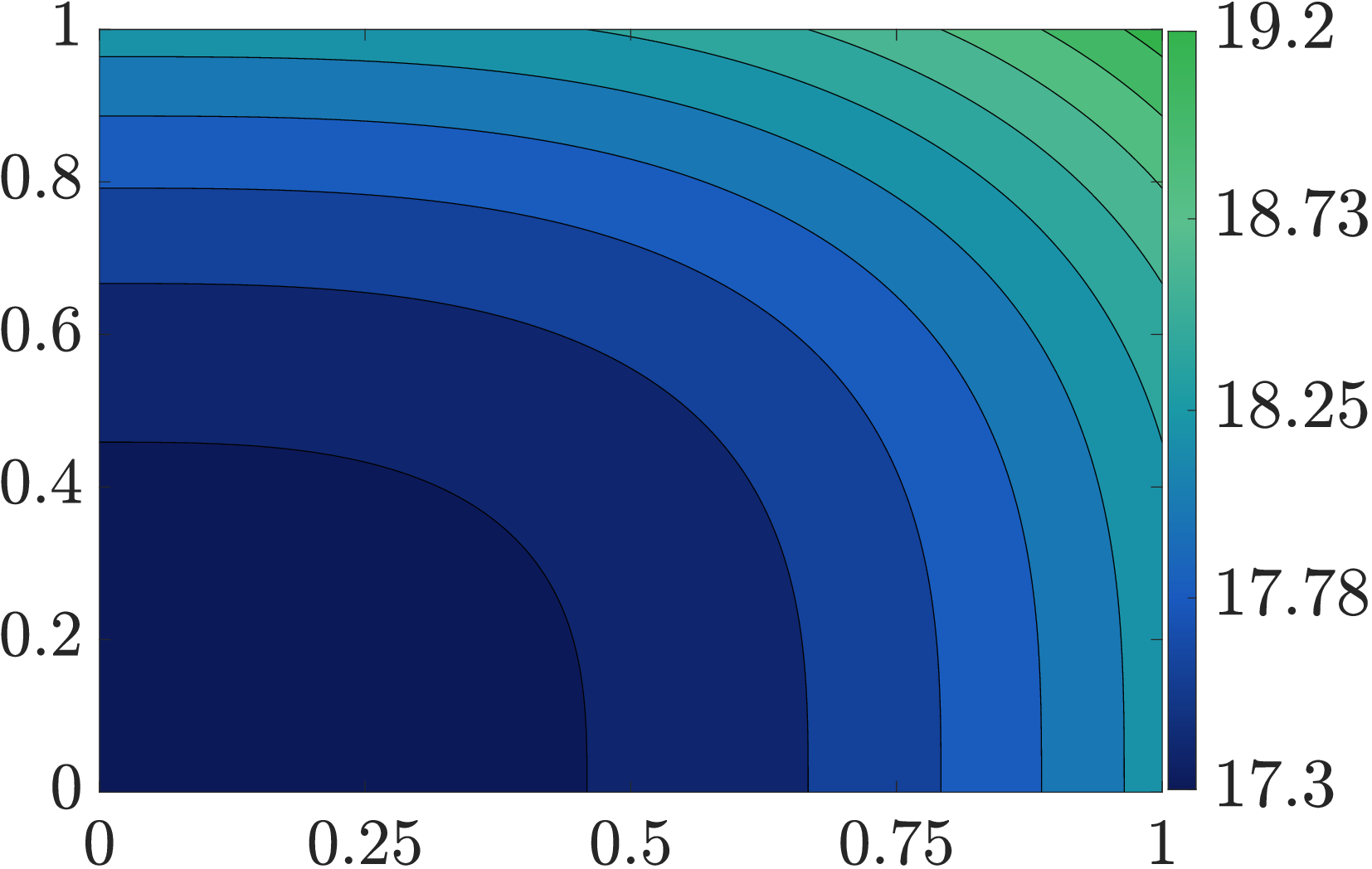}
\end{minipage}%
\begin{minipage}{0.333\textwidth}
\includegraphics[scale= 0.185]{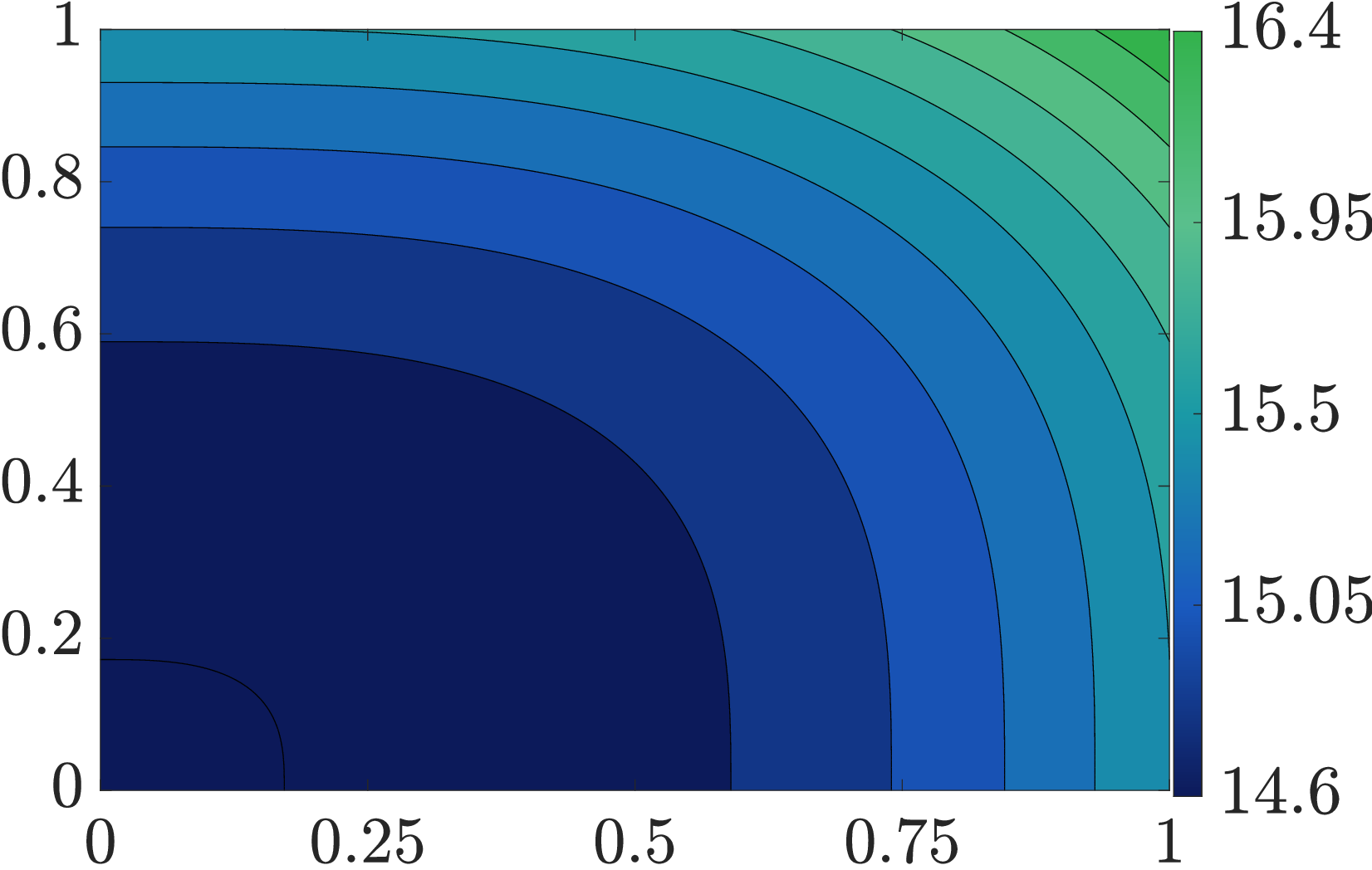}
\end{minipage}%
\begin{minipage}{0.333\textwidth} 
\includegraphics[scale= 0.185]{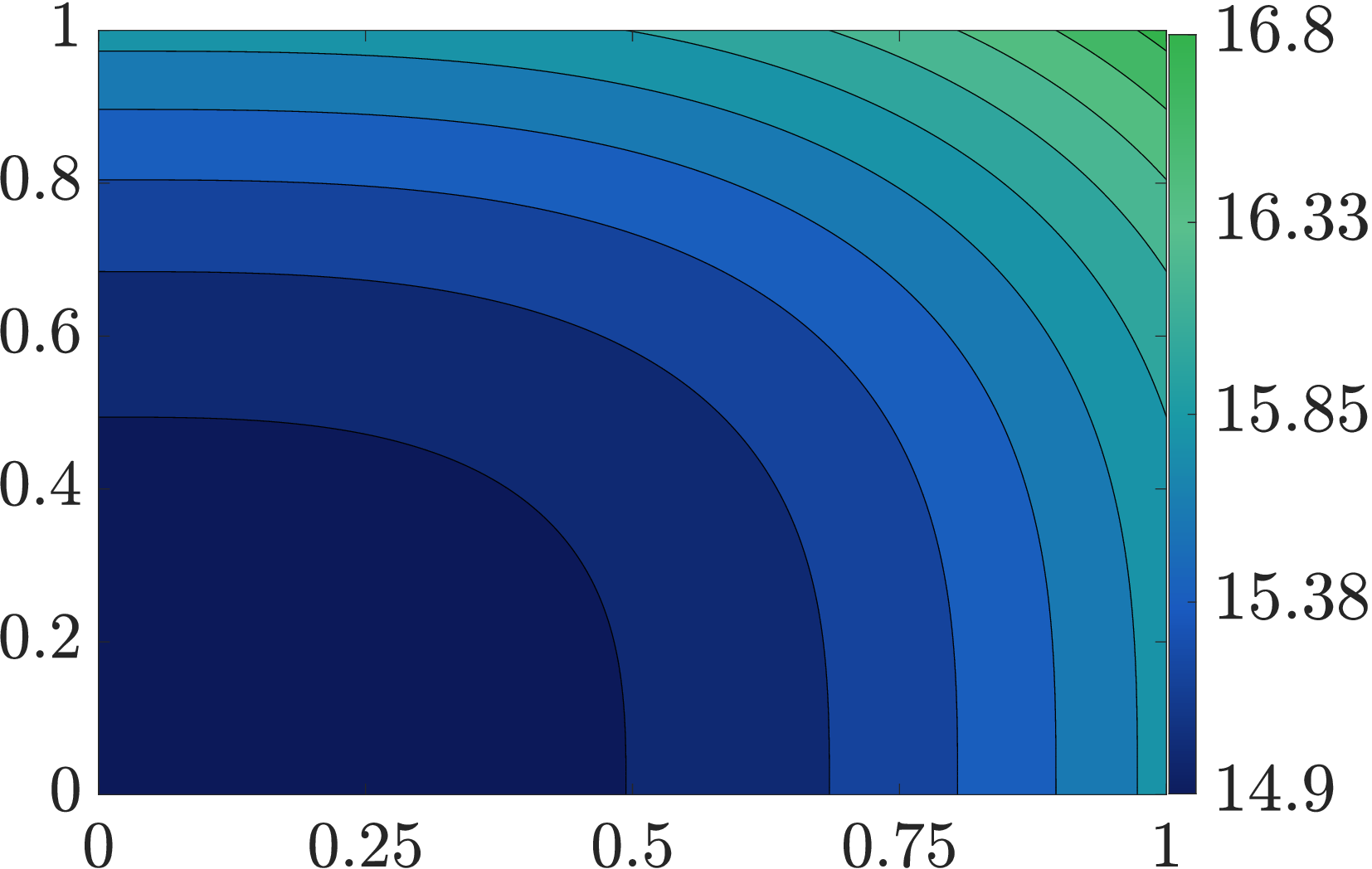}
\end{minipage}

\begin{minipage}{0.333\textwidth} 
\includegraphics[scale= 0.185]{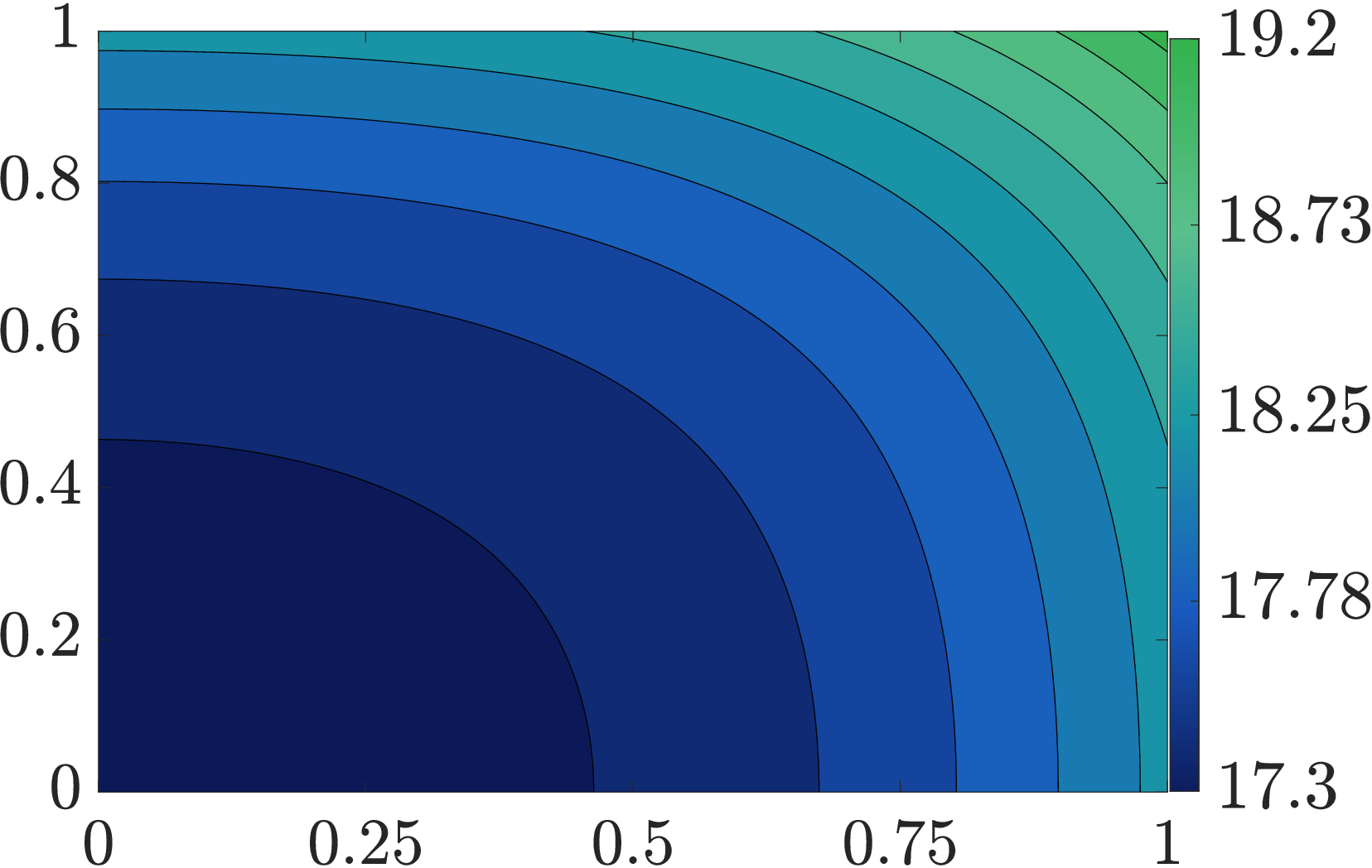}
\end{minipage}%
\begin{minipage}{0.333\textwidth}
\includegraphics[scale= 0.185]{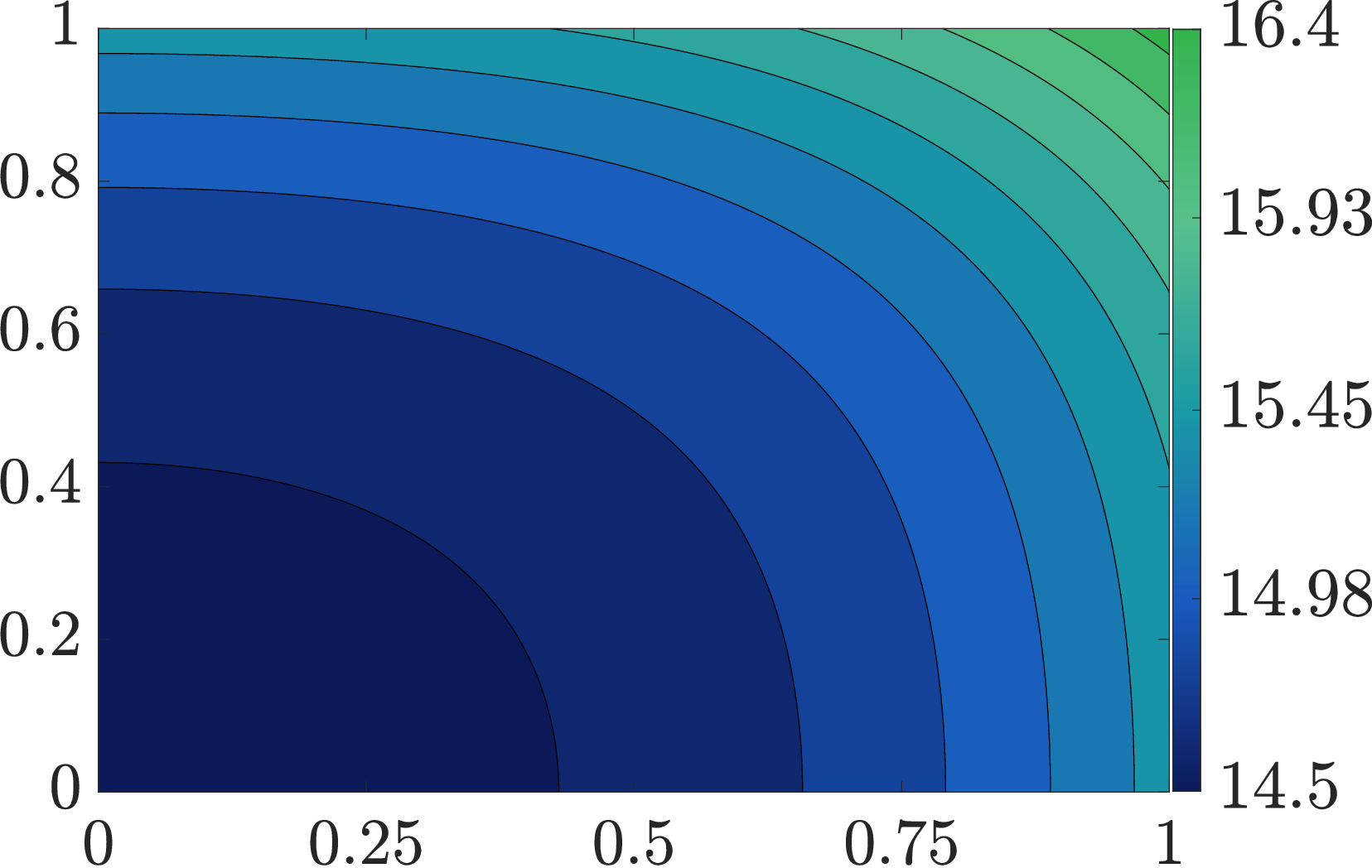}
\end{minipage}%
\begin{minipage}{0.333\textwidth} 
\includegraphics[scale= 0.185]{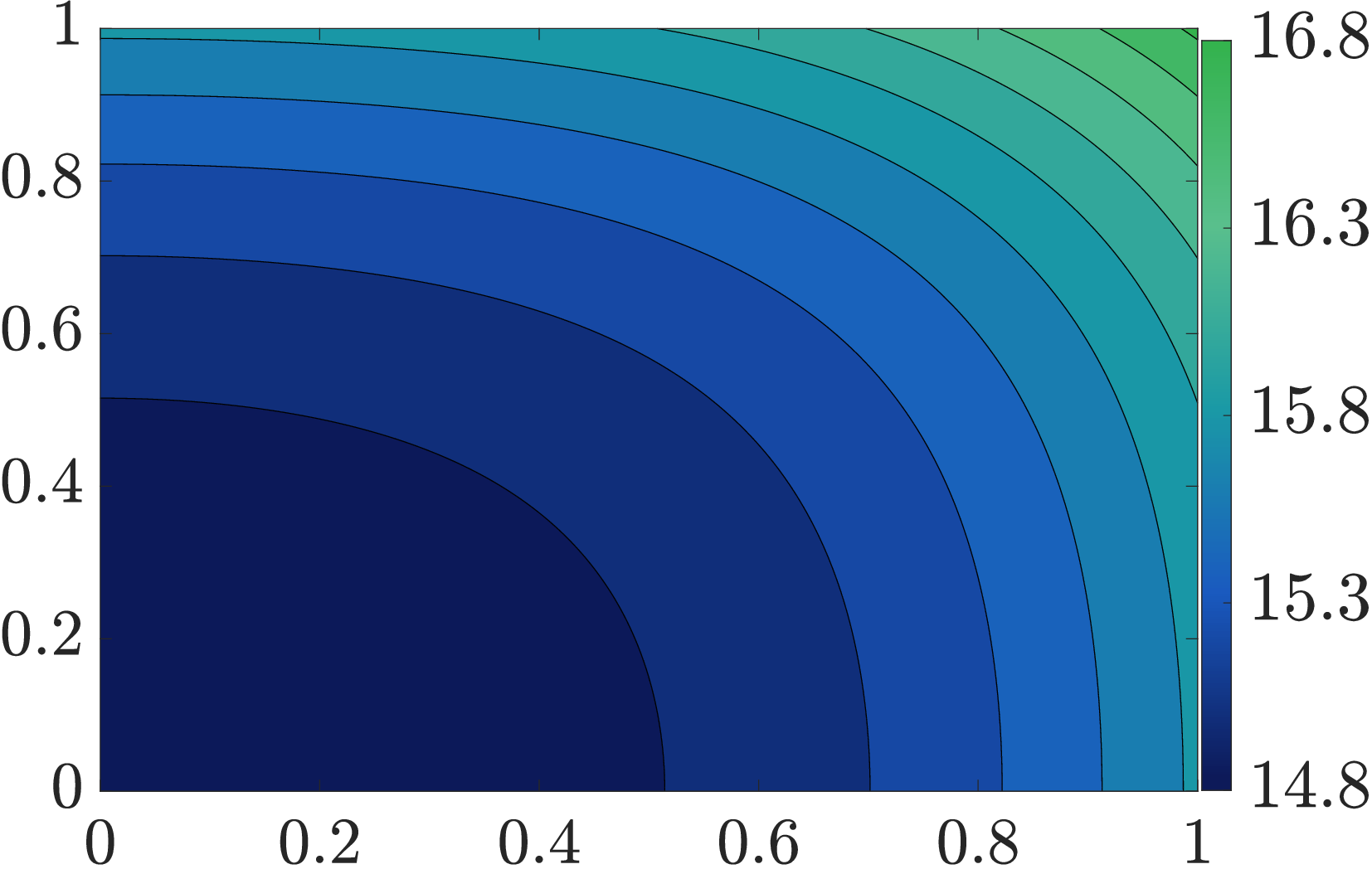}
\end{minipage}
\caption{\footnotesize [60D Reaction-diffusion] Heatmaps of the predicted solution with pool $\mathbb{P}_1$ in~\eqref{ReactDiff_pool1} andwith pool $\mathbb{P}_2$ in~\eqref{ReactDiff_pool2}. (First column) Dimensions $(17, 33)$. (Second column) Dimensions $(21, 56)$. (Third column) Dimensions $(52, 19)$.}
\label{ReactDiff_Pred}
\end{figure}

The true solution to~\eqref{60D_ReactDiff} along three pre-selected pairs of dimensions, with resolution $200 \times 200$, is displayed in Figure~\ref{ReactDiff_Ref}. The predicted solutions generated from the pools $\mathbb{P}_1$ and $\mathbb{P}_2$ are shown in Figure~\ref{ReactDiff_Pred}. Similar to the first test case, FEX accurately recovers the true solution when the candidate associated with $(\cdot)^3$ is included in the pool, whereas it yields only a moderately accurate approximation when this candidate is excluded. More specifically, when the pool $\mathbb{P}_2$ is used, FEX constructs the predicted solution from a combination of the candidates associated with $(\cdot)^4$ and $(\cdot)\sin(\cdot)$. Nevertheless, FEX is still able to fine-tune this combination to obtain a reasonably accurate approximation to the true solution.

\begin{figure}[h!]
\begin{minipage}{0.333\textwidth} 
\includegraphics[scale= 0.185]{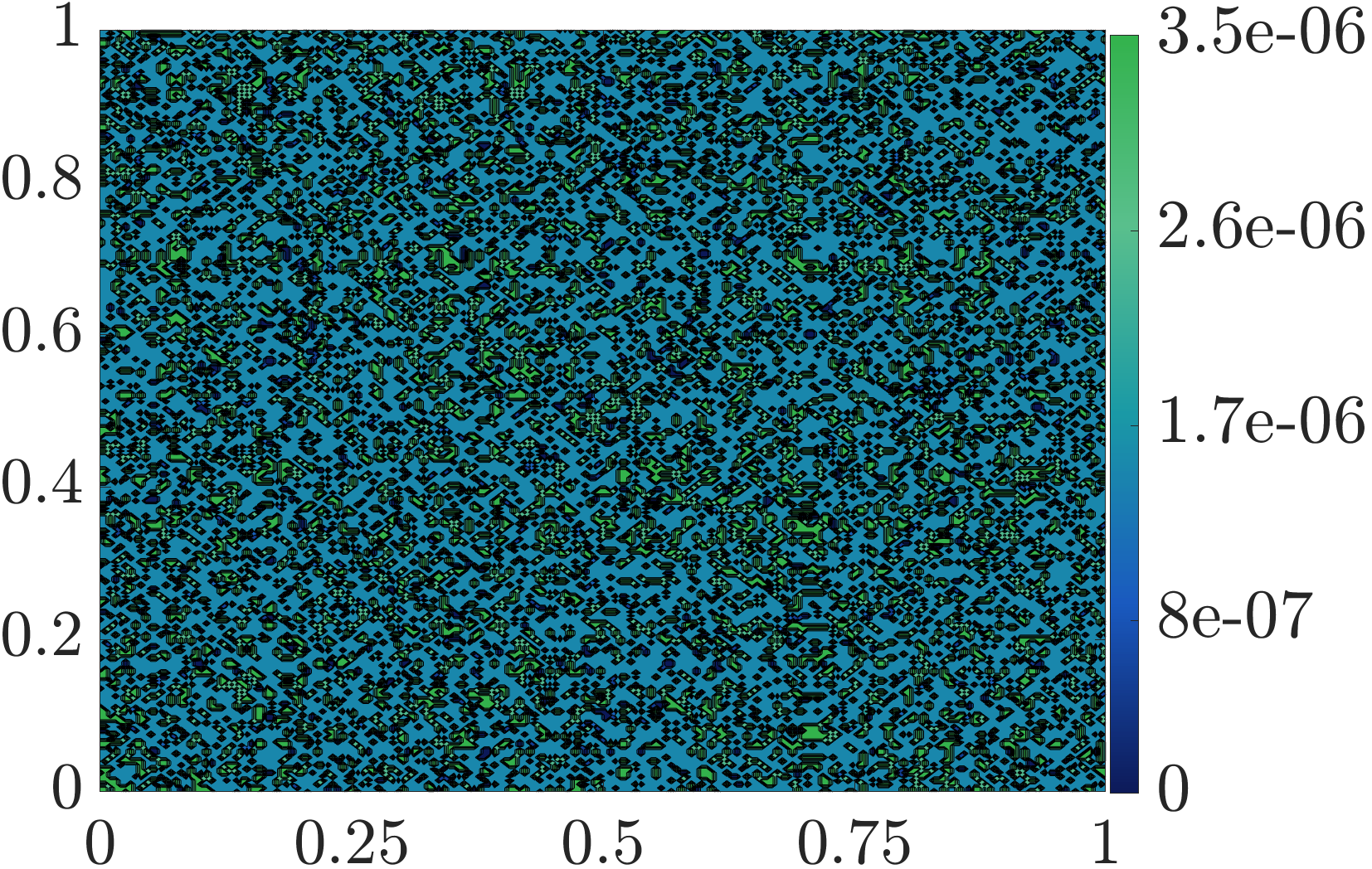}
\end{minipage}%
\begin{minipage}{0.333\textwidth}
\includegraphics[scale= 0.185]{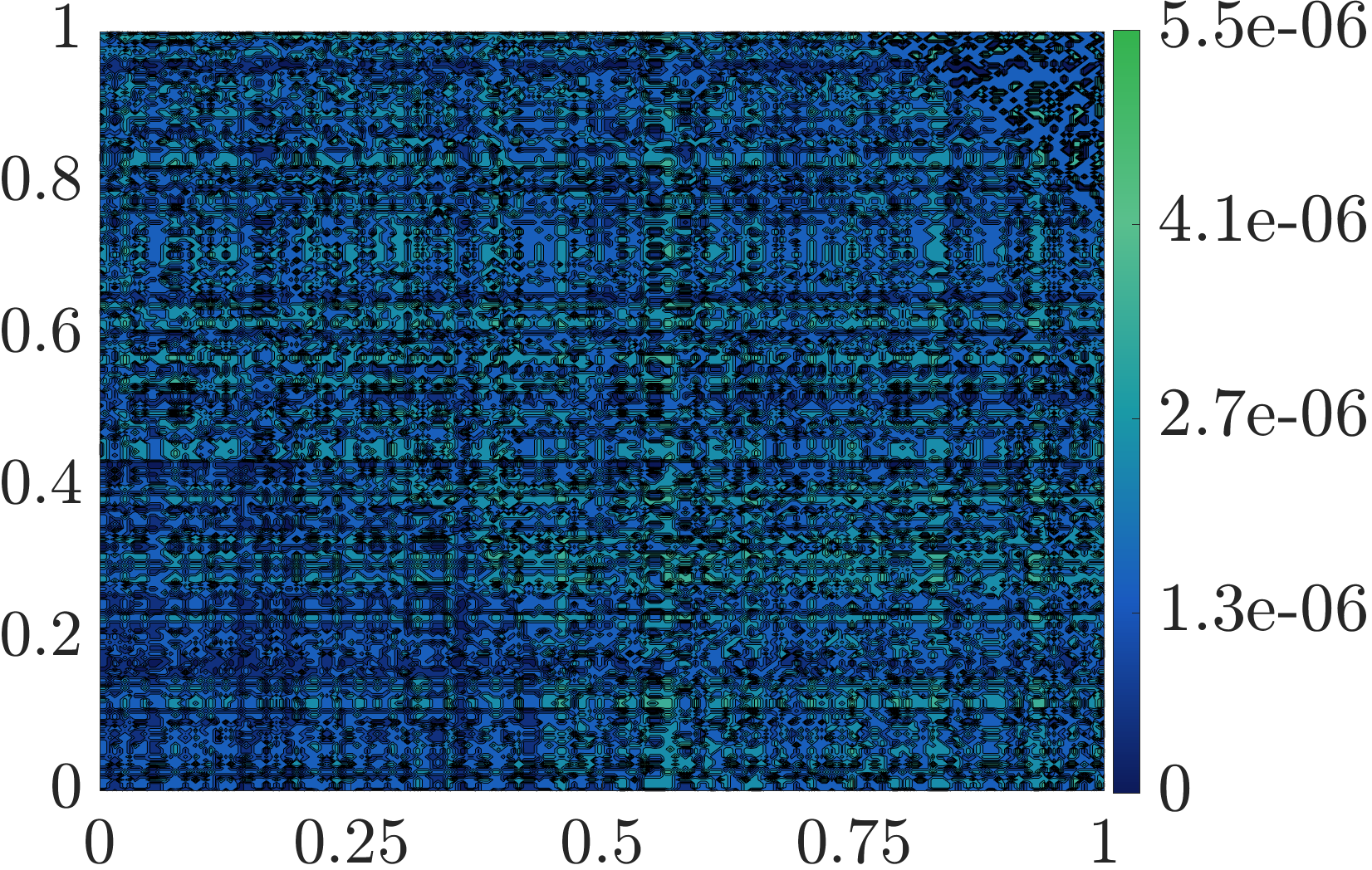}
\end{minipage}%
\begin{minipage}{0.333\textwidth} 
\includegraphics[scale= 0.185]{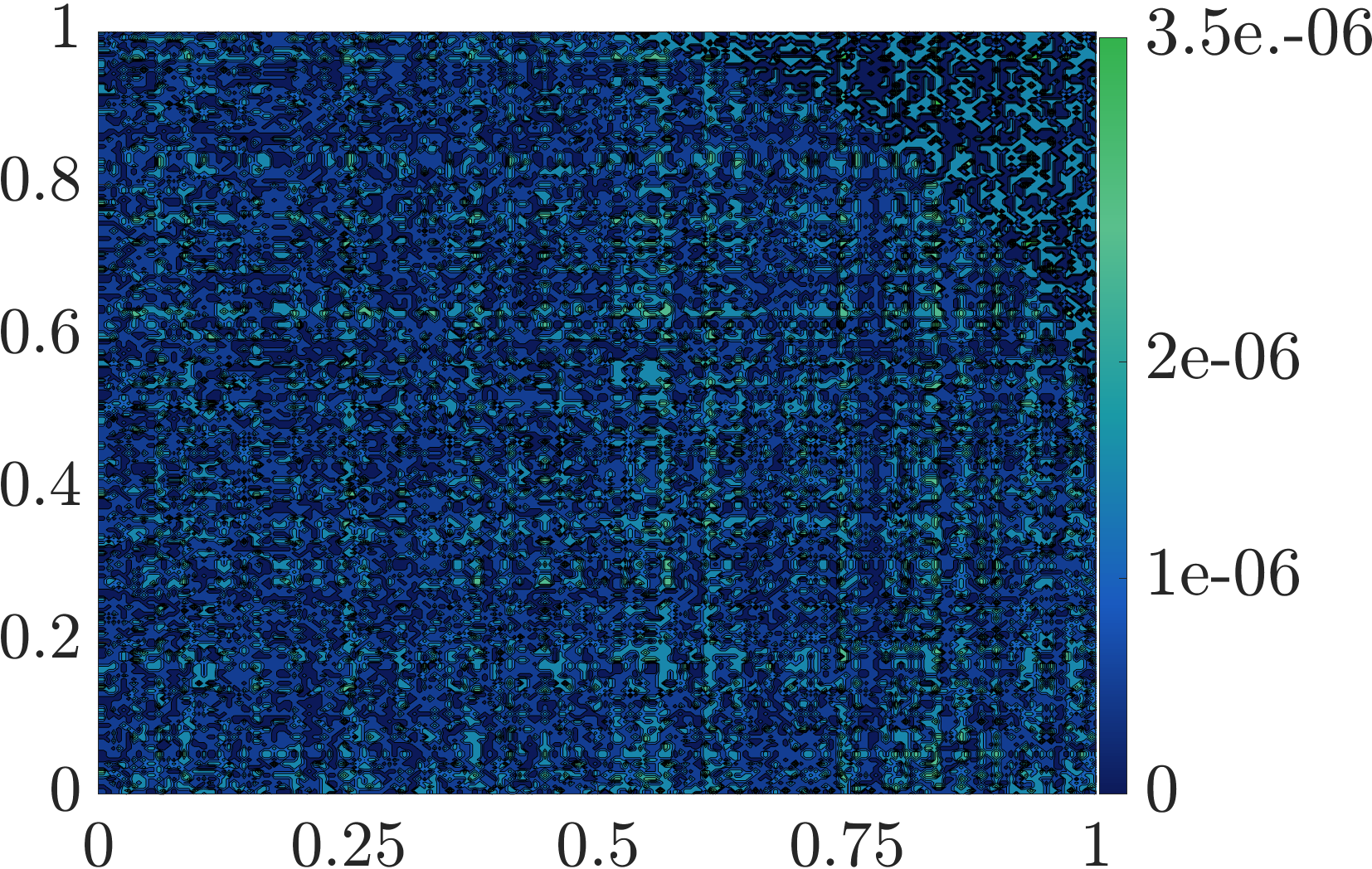}
\end{minipage}

\begin{minipage}{0.333\textwidth} 
\includegraphics[scale= 0.185]{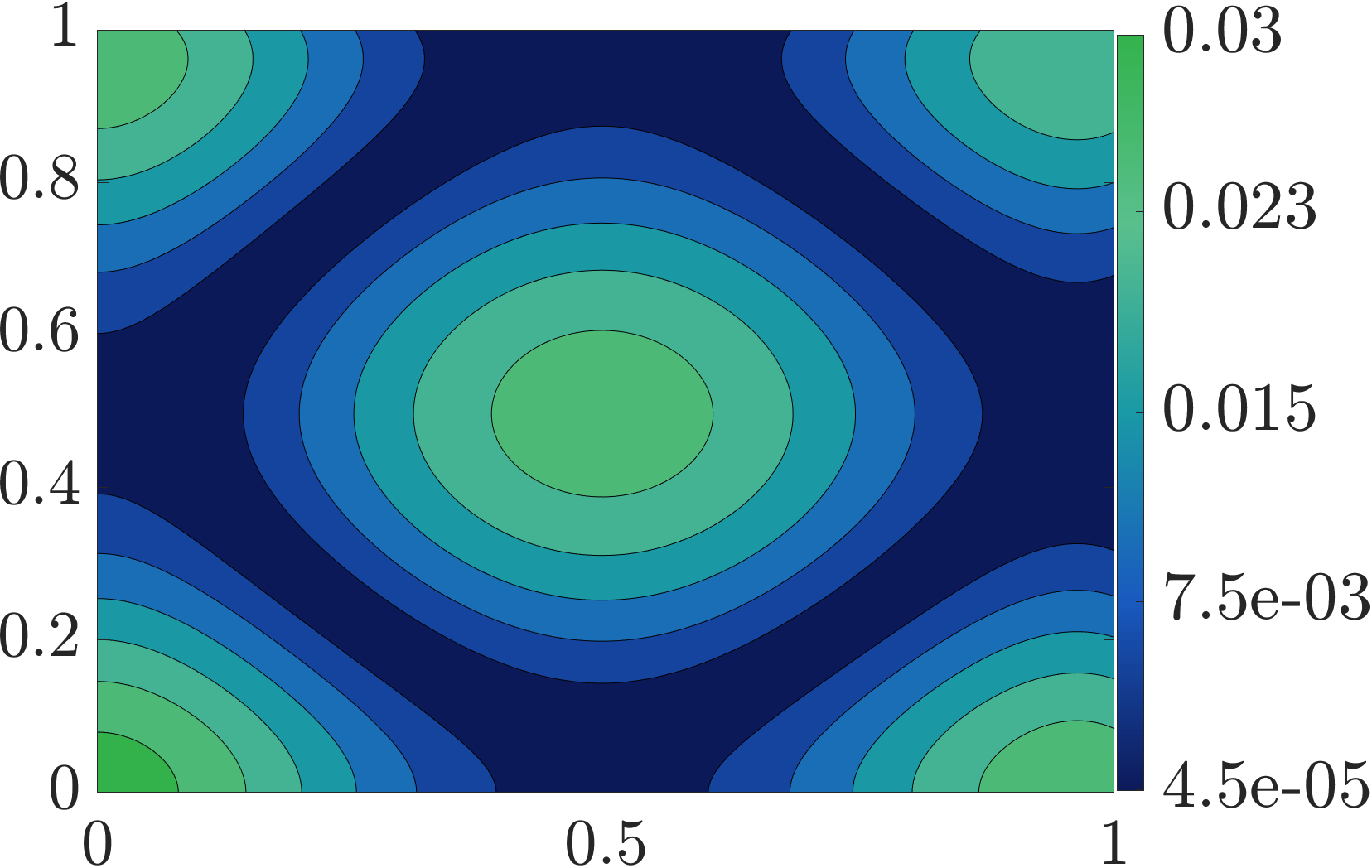}
\end{minipage}%
\begin{minipage}{0.333\textwidth}
\includegraphics[scale= 0.185]{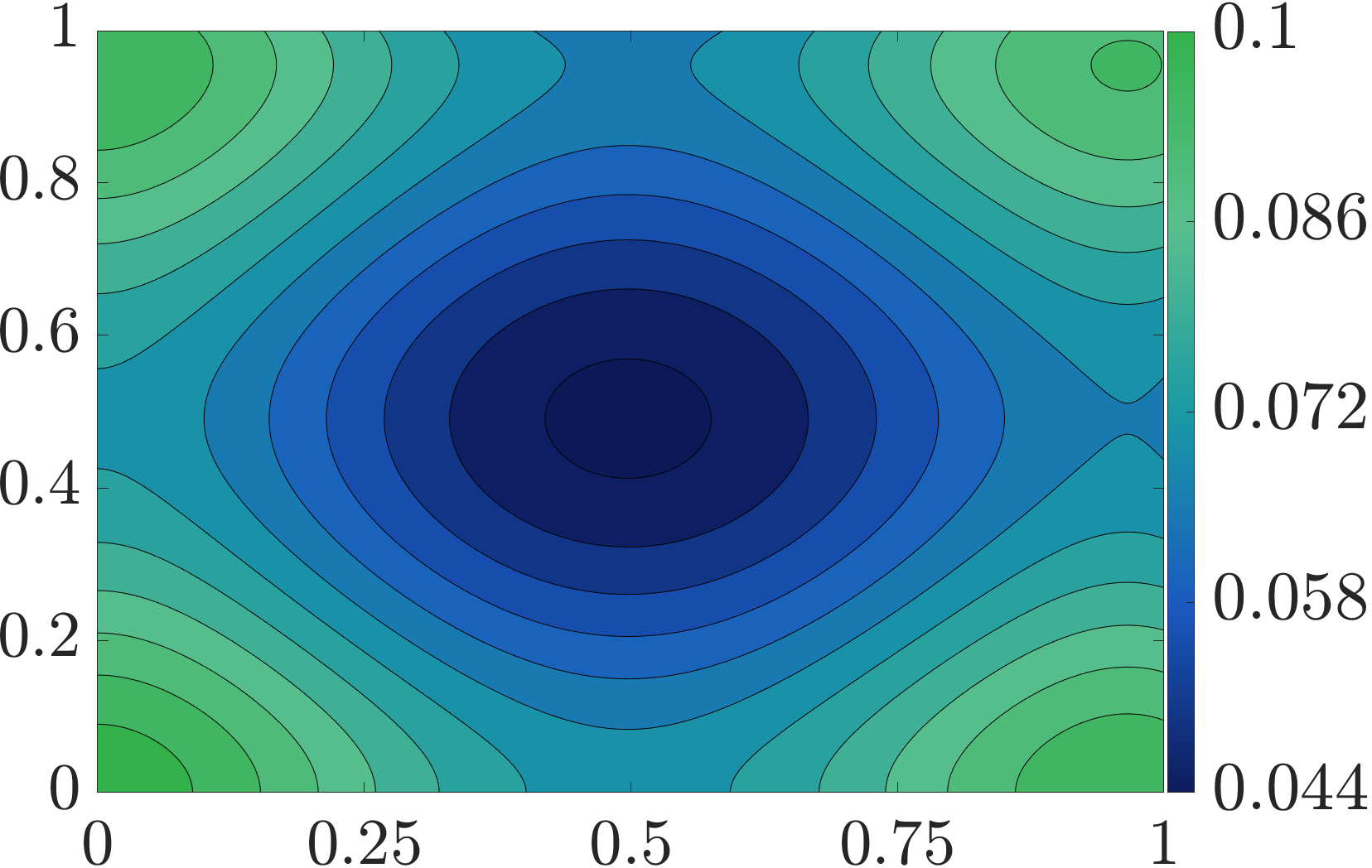}
\end{minipage}%
\begin{minipage}{0.333\textwidth} 
\includegraphics[scale= 0.185]{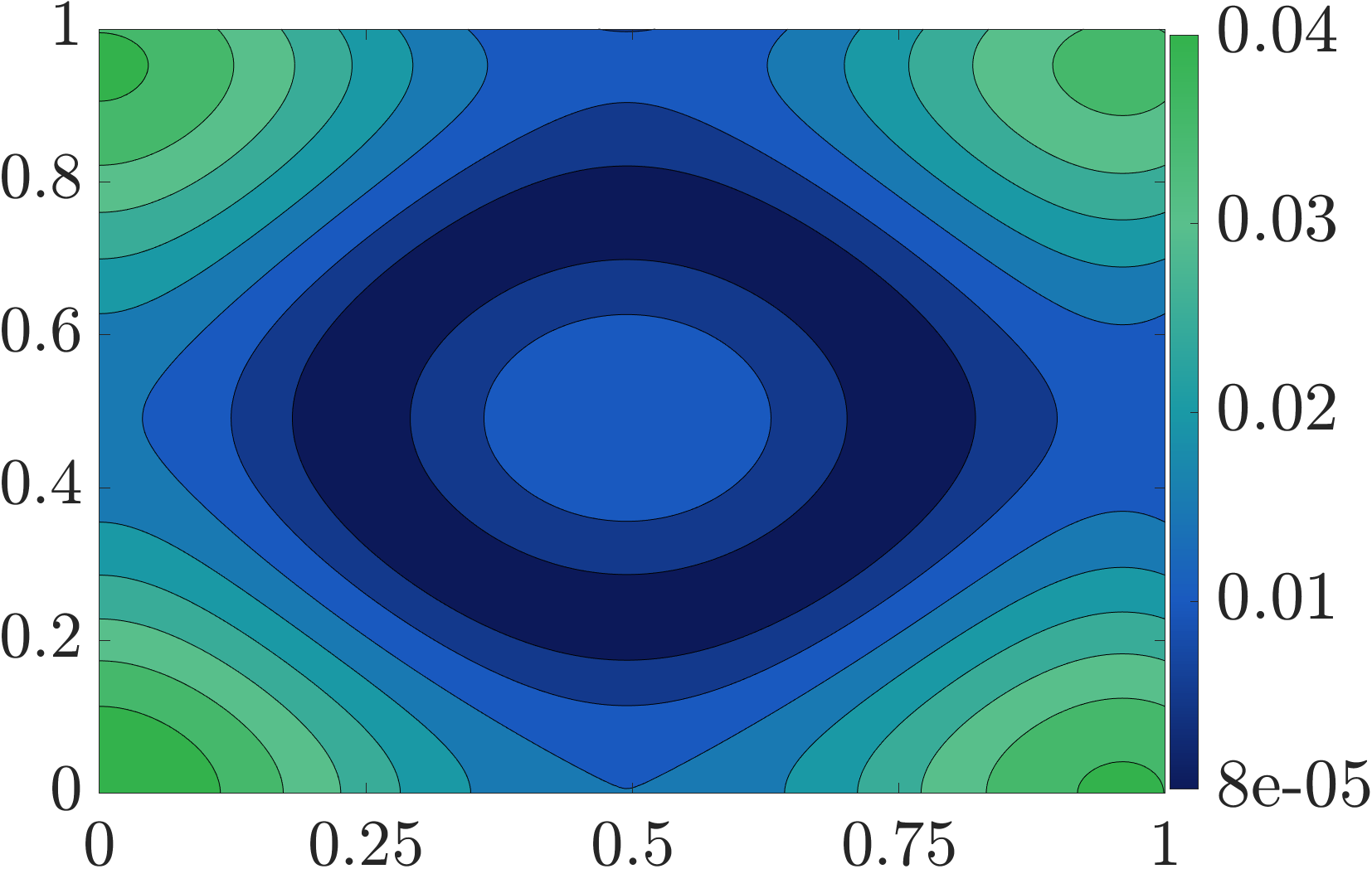}
\end{minipage}
\caption{\footnotesize [60D Reaction-diffusion] Heatmaps of absolute pointwise error with pool $\mathbb{P}_1$ in~\eqref{ReactDiff_pool1} andwith pool $\mathbb{P}_2$ in~\eqref{ReactDiff_pool2}. (First column) Dimensions $(17, 33)$. (Second column) Dimensions $(21, 56)$. (Third column) Dimensions $(52, 19)$.}
\label{ReactDiff_Error}
\end{figure}
We further compute the absolute pointwise errors for both approximate solutions and plot the corresponding heatmaps in Figure~\ref{ReactDiff_Error}. These heatmaps support the behavior observed in Figure~\ref{ReactDiff_Pred}. More precisely, the errors corresponding to the solution generated from pool $\mathbb{P}_1$ are substantially smaller, whereas those corresponding to the solution generated from pool $\mathbb{P}_2$ are larger. Finally, we report the sample means and standard deviations of the Monte Carlo approximations of the relative $L^2$ errors in Table~\ref{L2_Errors}. These values are consistent with the results obtained in the first test case.

\subsection{Test case 3: 55D Semilinear elliptic equation}
In the final test case, we consider a semilinear elliptic equation in $55$ dimensions:
\begin{equation}
\label{55D_Semilinear}
-\nu \Delta u + \mu u + u^2 = f(\pmb{x}), \quad \pmb{x} \in \Omega,
\end{equation}
where $\Omega = (-1,1)^d$ with $d = 55$. As a more challenging example, we assume that the true solution to~\eqref{55D_Semilinear} takes the form
\[
u(\pmb{x}) = \exp\left(\frac{1}{d}\sum_{j=1}^d \cos(x_j)\right).
\]
The candidate pool for generating operator sequences is given by
\begin{equation}
\label{Semilinear_Pool}
\mathbb{P} := \{0, 1, \mathrm{Id}, \mathrm{TN}_{x^2}, \mathrm{TN}_{x^3}, \mathrm{TN}_{x^4}, \mathrm{TN}_{\exp}, \mathrm{TN}_{\sin(x)}, \mathrm{TN}_{\cos(x)} \}.
\end{equation}
We emphasize that, although the candidate pool contains the basic operators needed to represent the true solution, the latter is not a simple additive combination of component-wise nonlinearities, unlike in the previous test cases. Our goal is to investigate whether FEX can still recover this type of composite expression.

\begin{figure}[h!]
\vspace{-0.3cm}
\begin{minipage}{0.333\textwidth} 
\includegraphics[scale= 0.177]{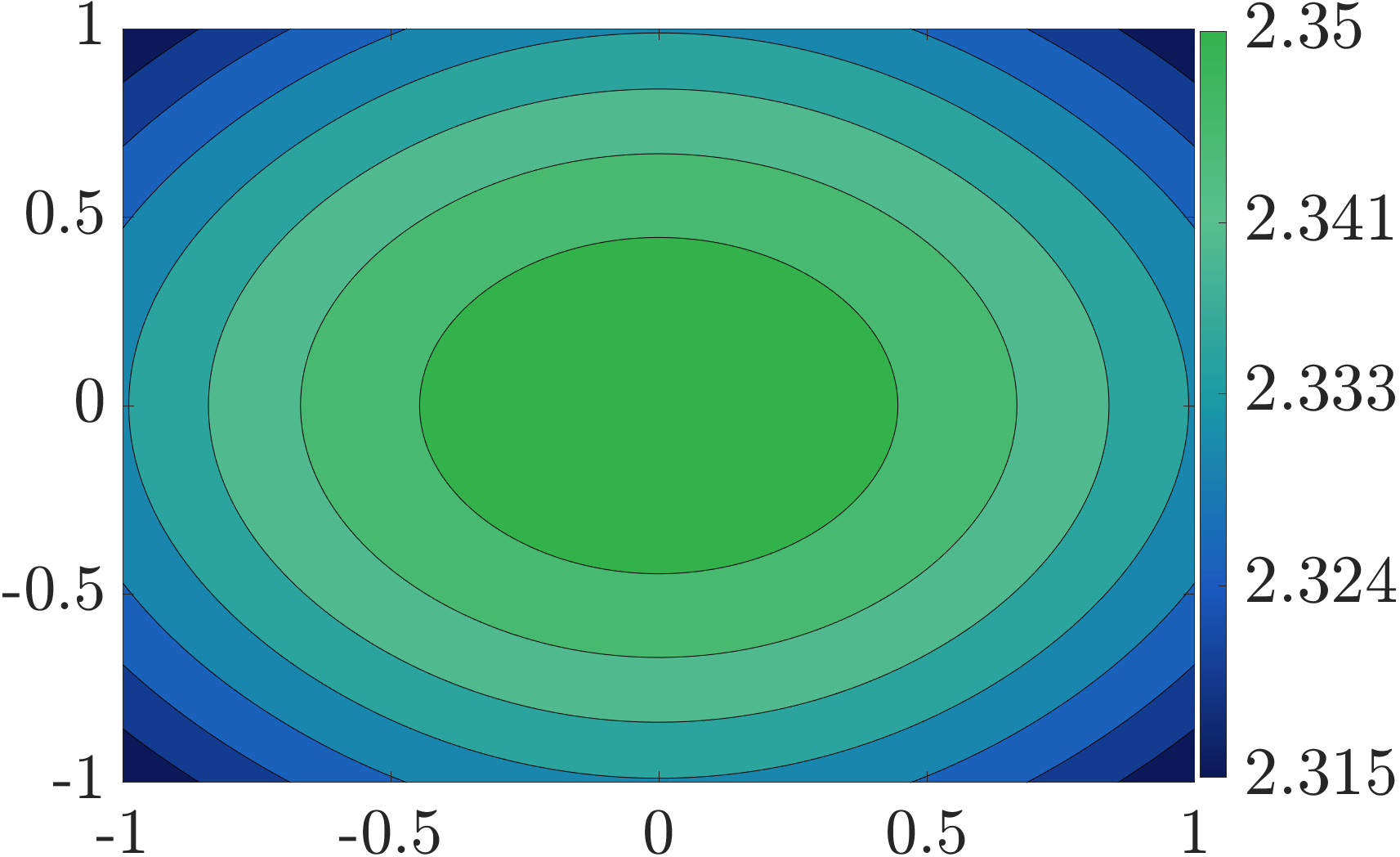}
\end{minipage}%
\begin{minipage}{0.333\textwidth}
\includegraphics[scale= 0.177]{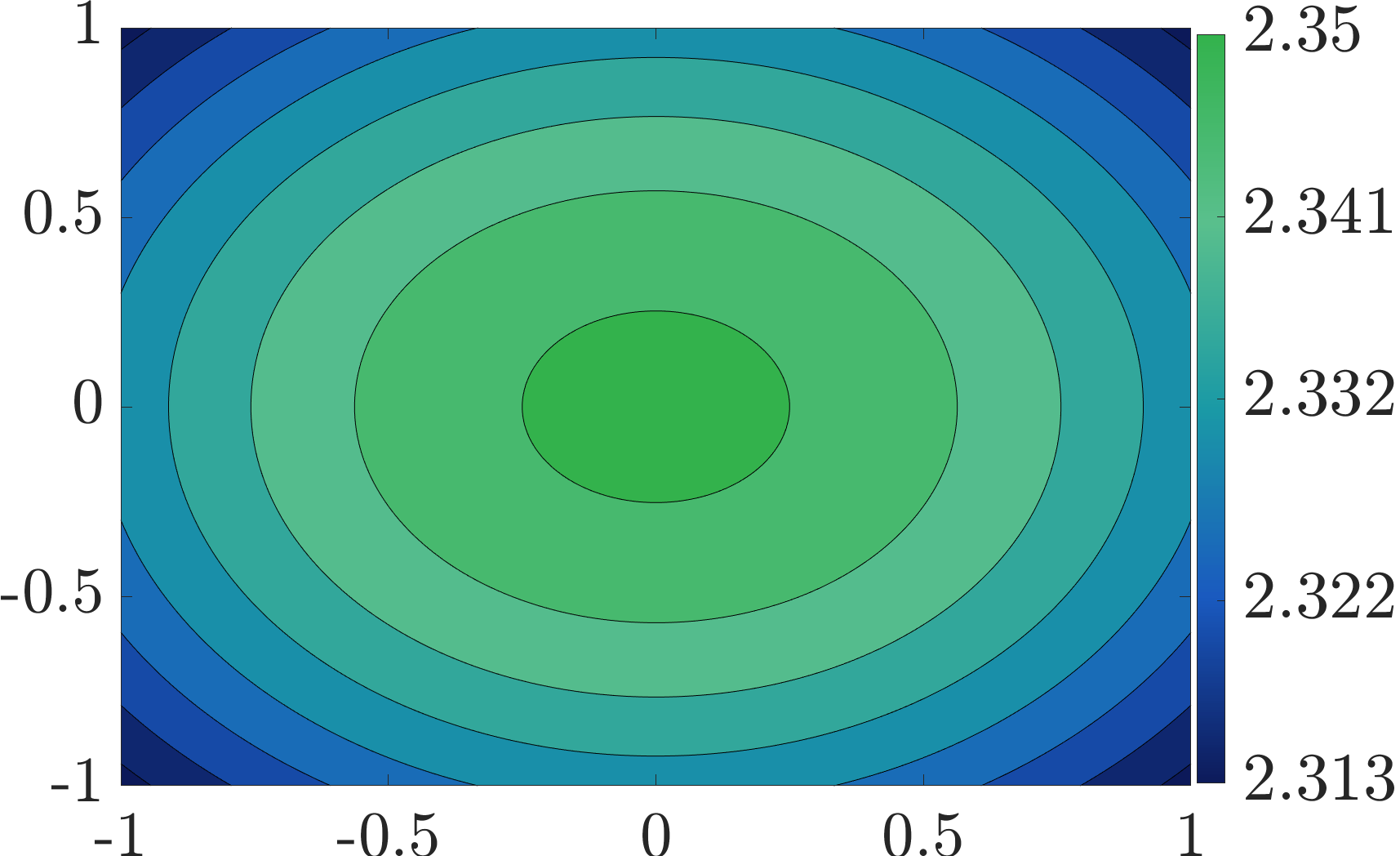}
\end{minipage}%
\begin{minipage}{0.333\textwidth} 
\includegraphics[scale= 0.177]{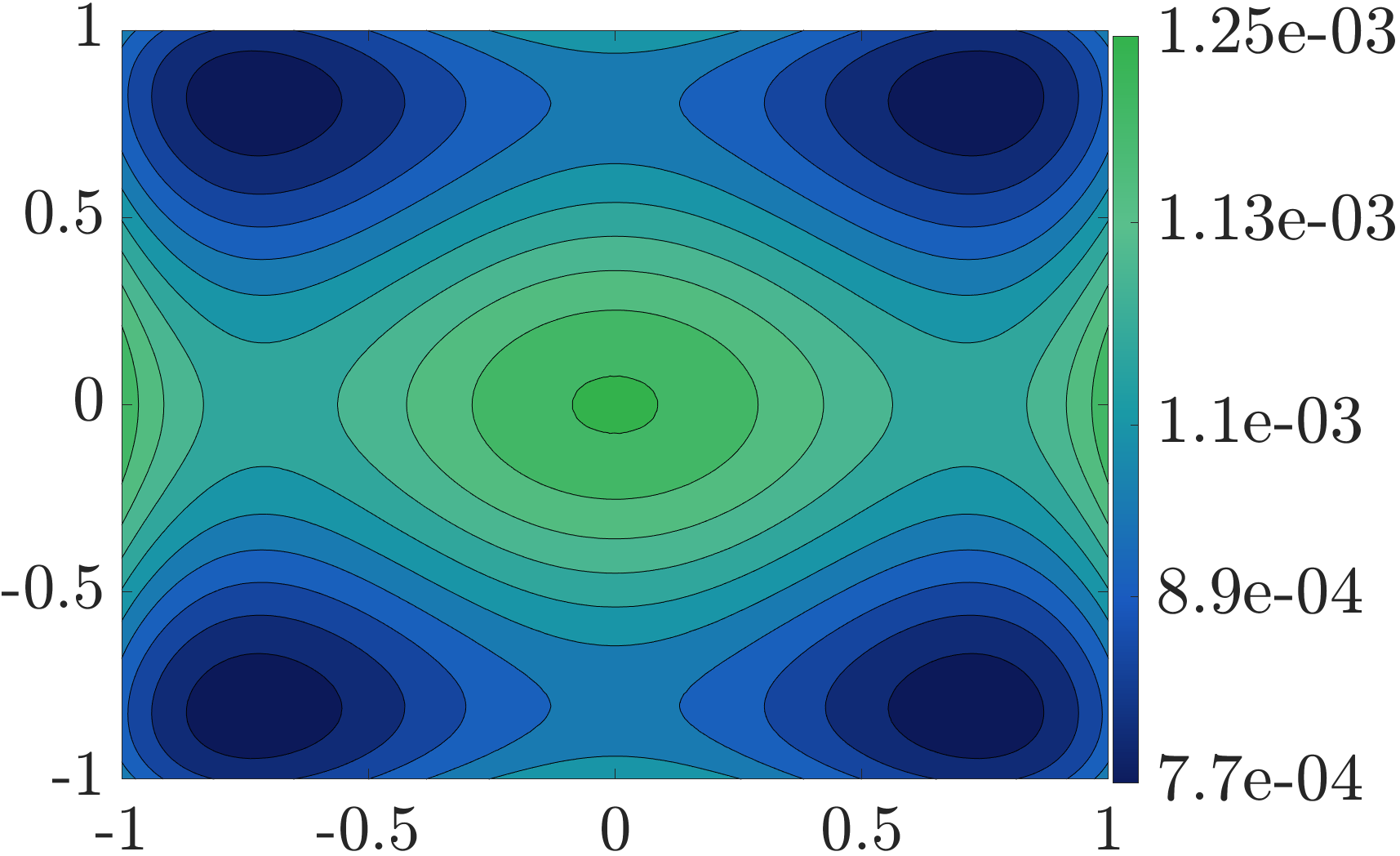}
\end{minipage}

\begin{minipage}{0.333\textwidth} 
\includegraphics[scale= 0.177]{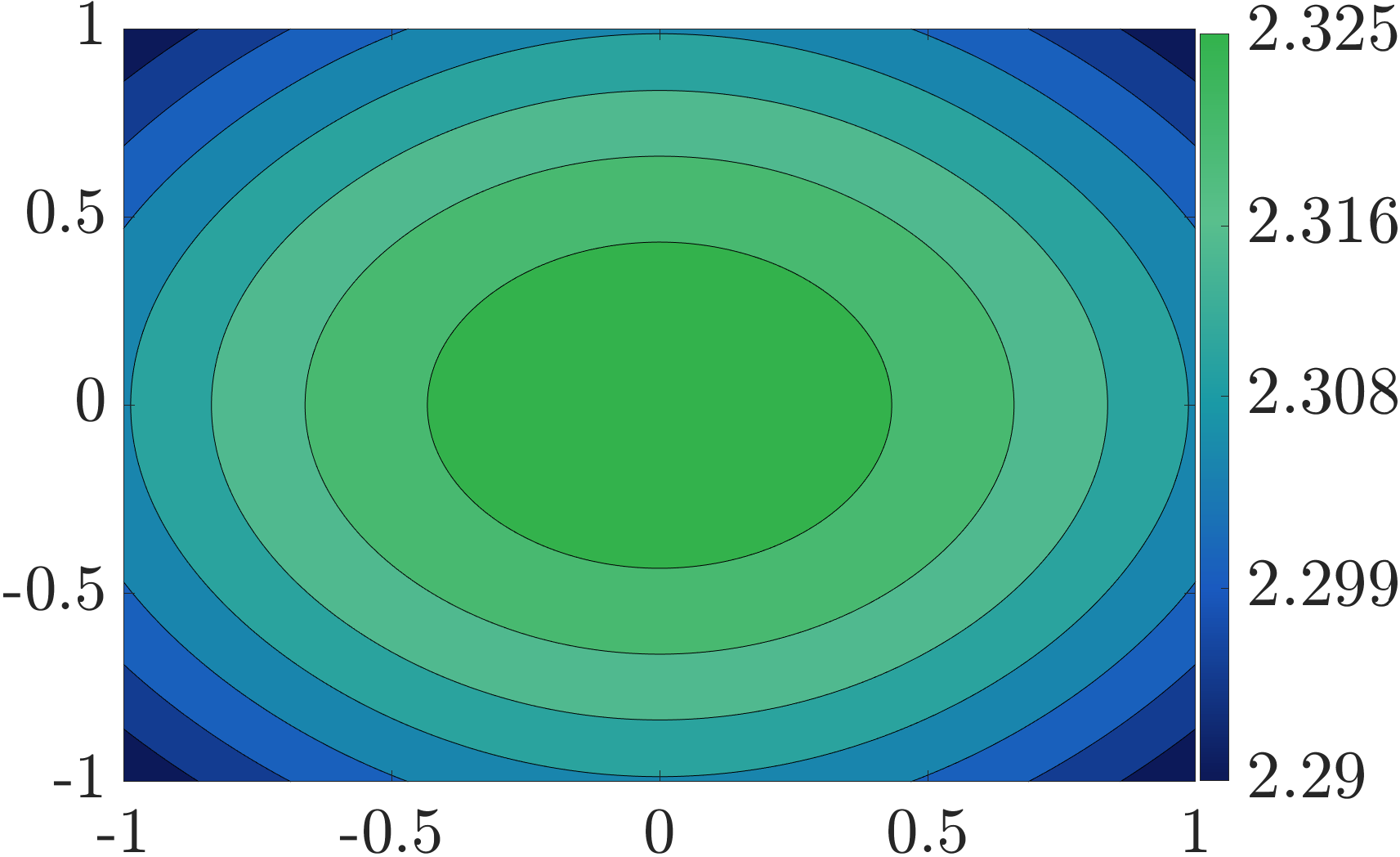}
\end{minipage}%
\begin{minipage}{0.333\textwidth}
\includegraphics[scale= 0.177]{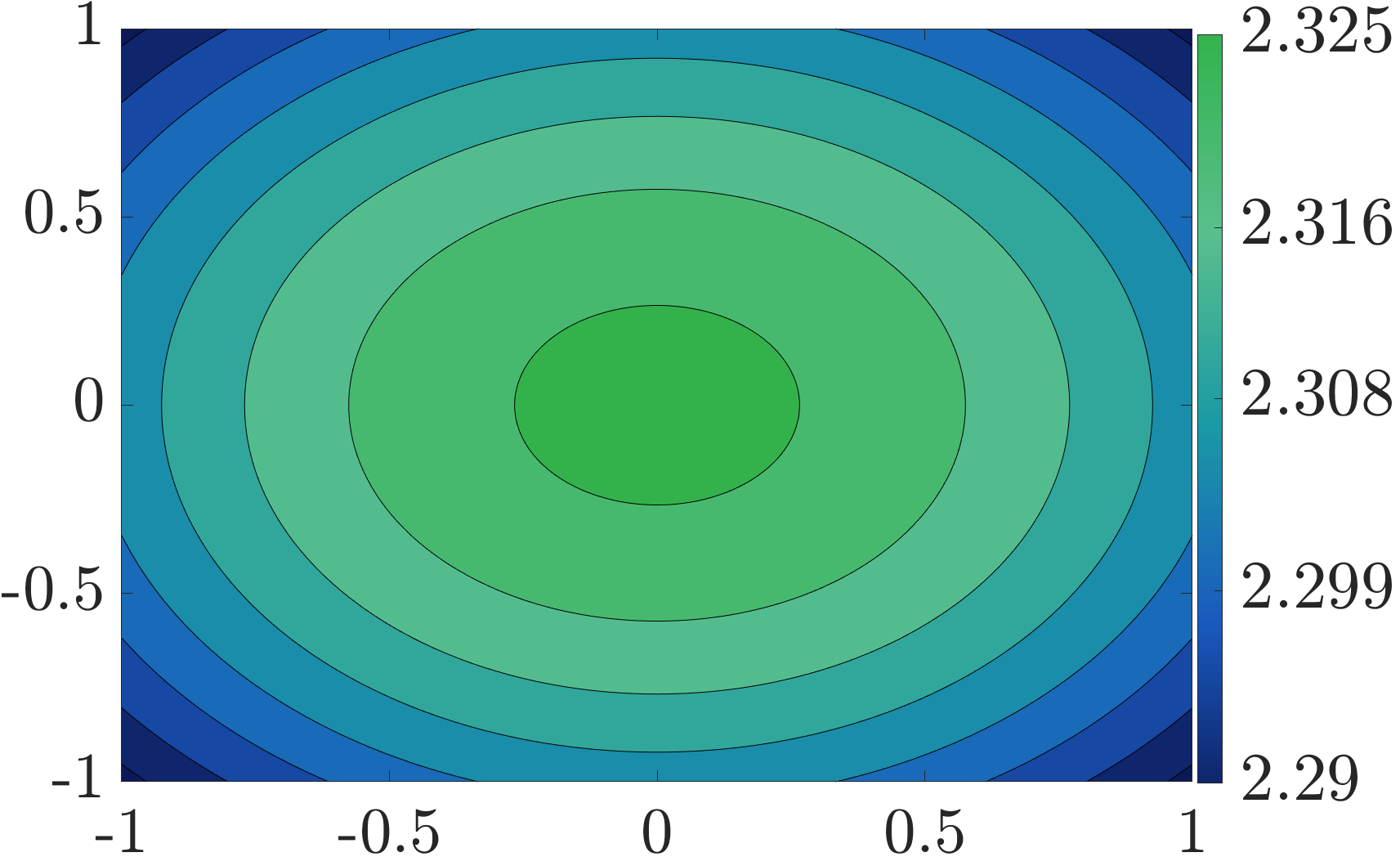}
\end{minipage}%
\begin{minipage}{0.333\textwidth} 
\includegraphics[scale= 0.177]{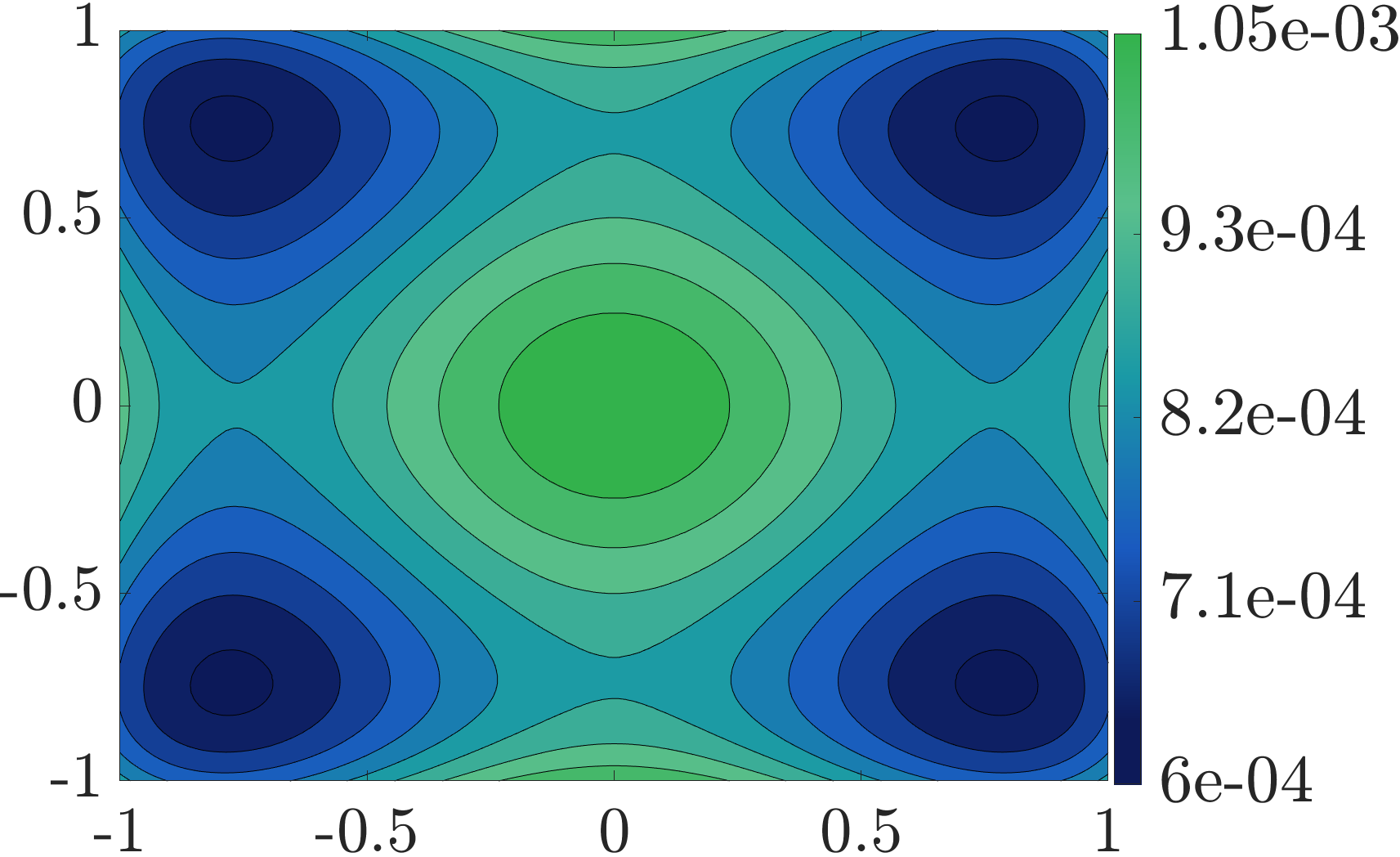}
\end{minipage}

\begin{minipage}{0.333\textwidth} 
\includegraphics[scale= 0.177]{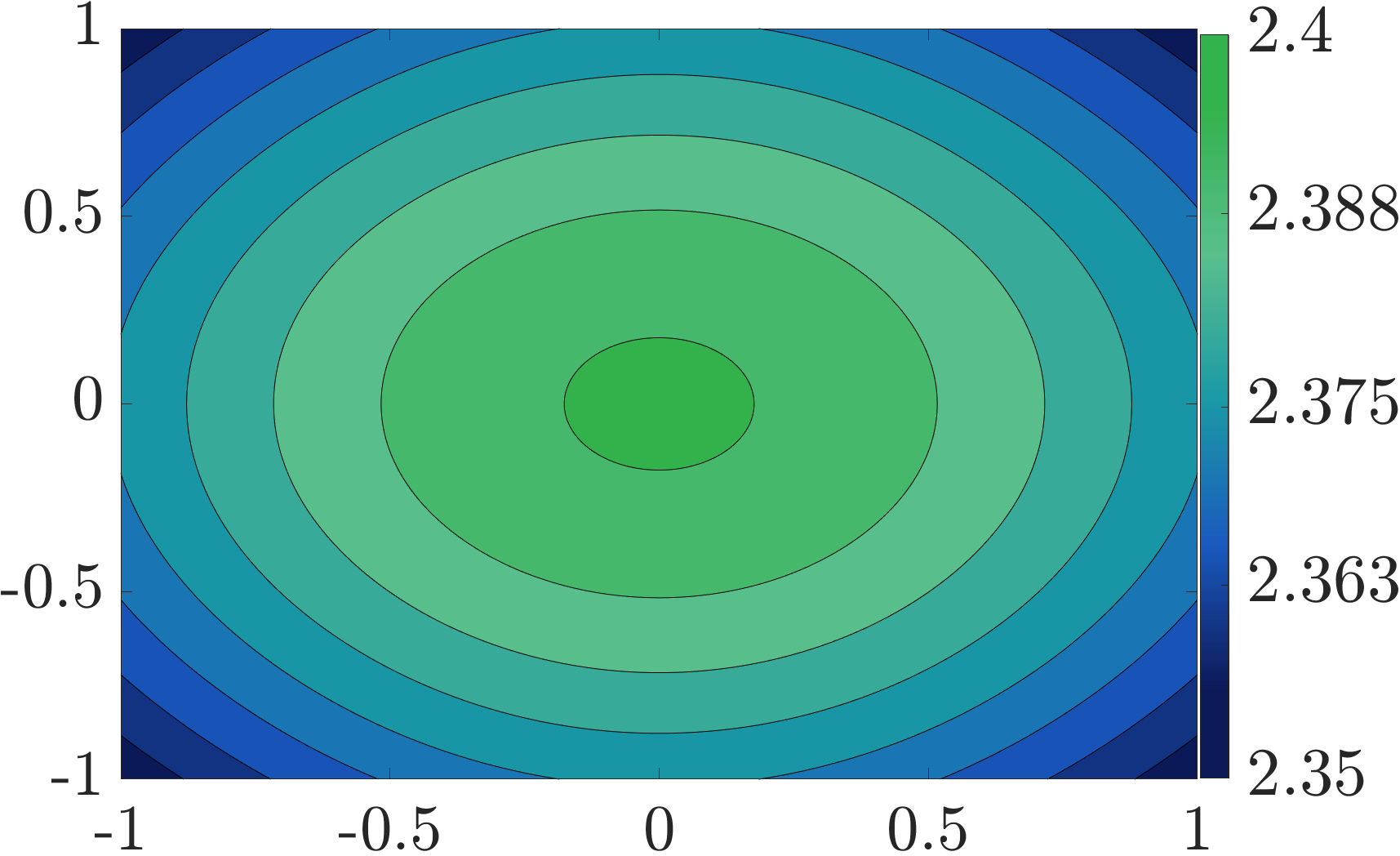}
\end{minipage}%
\begin{minipage}{0.333\textwidth}
\includegraphics[scale= 0.177]{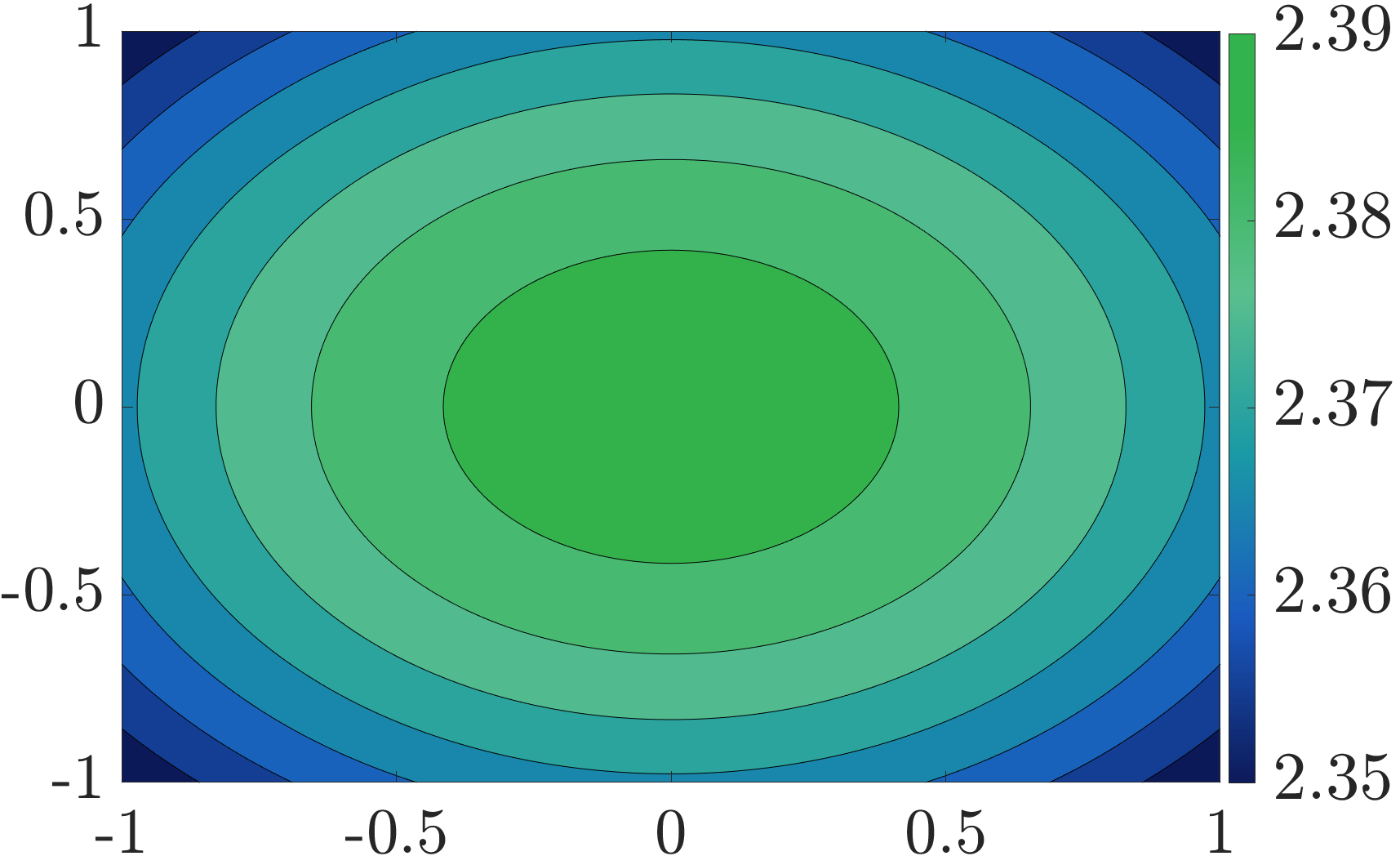}
\end{minipage}%
\begin{minipage}{0.333\textwidth} 
\includegraphics[scale= 0.177]{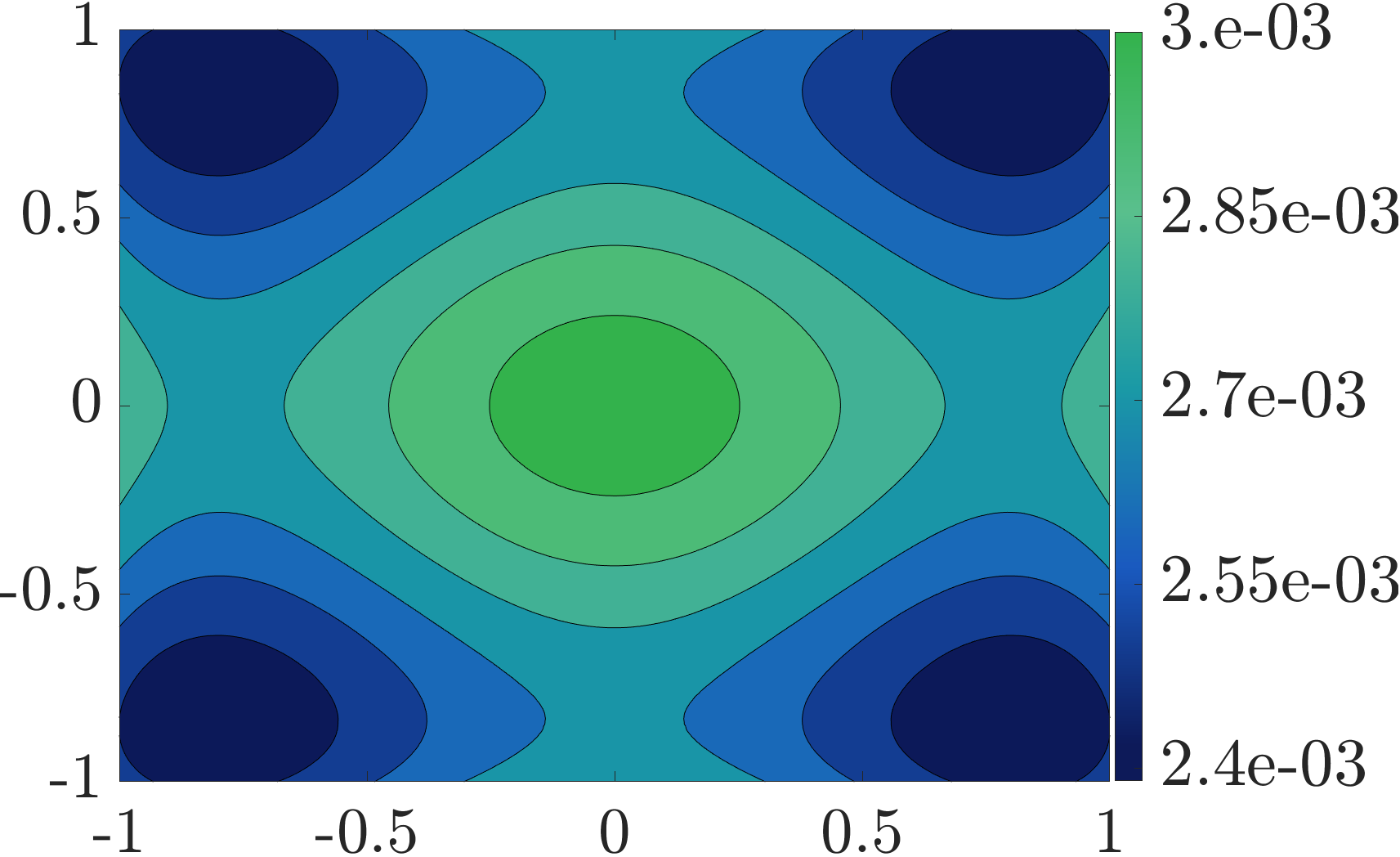}
\end{minipage}
\caption{\footnotesize [55D Semilinear Elliptic] Heatmaps of the reference solution, predicted solution, and the corresponding absolute relative error on two-dimensional slices, with the remaining $58$ dimensions fixed at predefined values. (First row) Dimensions $(27, 32)$. (Second row) Dimensions $(30, 35)$. (Third row) Dimensions $(30, 35)$.}
\label{Semilinear_Compare_ErrorMap}
\end{figure}

We plot the true solution and the predicted solution along selected pairs of dimensions, with resolution $300 \times 300$, in Figure~\ref{Semilinear_Compare_ErrorMap}. From the first and second columns of Figure~\ref{Semilinear_Compare_ErrorMap}, we observe that the two solutions have similar shapes but differ somewhat in magnitude. After the search loop is completed, FEX produces an approximation involving the composite expression $\exp(-(\cdot)^2)$, whose overall shape is similar to that of $\exp(\cos(\cdot))$. This predicted solution is then further fine-tuned during the optimization step to better match the magnitude of the true solution. To further quantify the difference between the two solutions, we also show the \textit{absolute relative pointwise errors} in the last column of Figure~\ref{Semilinear_Compare_ErrorMap}. These results show that, even though FEX does not recover the exact form of the true solution, it is still able to produce a good approximation after the fine-tuning step.

We conclude this section by computing Monte Carlo approximations of the relative $L^2$ error between the true solution and the predicted solution. Since the sampling points used in Monte Carlo integration are randomly generated, the approximation is repeated $50$ times, and the corresponding sample means and standard deviations are reported in Table~\ref{L2_Errors}. These values indicate that the predicted solution provides a reasonably accurate and stable approximation.
% \vspace{-0.3cm}
\subsection{Table of $\mathrm{L}^2$-Errors}
In this section, we present the \textit{relative $L^2$ errors}, defined by $\dfrac{\norm{u-\tilde{u}}_{L^2(\Omega)}}{\norm{u}_{L^2(\Omega)}}$, for the three numerical test cases above. Since these integrals arise in high-dimensional settings, they are estimated using Monte Carlo integration. Specifically, we sample $2000$ interior points to compute the estimates. For completeness, we recall that each error is computed $50$ times to ensure reliability, and the corresponding mean and standard deviation are reported in Table~\ref{L2_Errors}.

\begin{table*}[h!]
\centering
\begin{tabular}{|c|c|c|c|c|}
\hline
\multicolumn{3}{|c|}{} & \text{Mean} & \text{Std}\\
\hline
\multirow{5}{*}{Relateve $\mathrm{L}^2$-Errors} & \multirow{2}{*}{Poisson} &\text{First pool} & $4.17e-07$ & $1.51e-08$ \\
\cline{3-5}
& & Second pool & $3.13e-03$ & $1.24e-04$ \\
\cline{2-5} 
& \multirow{2}{*}{Reaction-Diffusion} & First pool & $7.76e-07$ & $2.20e-08$ \\
\cline{3-5}
& & Second pool & $5.17e-03$ & $1.73e-04$ \\
\cline{2-5} 
& {Semilinear Elliptic}& - & $9.33e-04$ & $3.14e-05$\\ 
\hline
\end{tabular}
\caption{Mean and standard deviation for relative $\mathrm{L}^2$-errors from all test cases.}
\label{L2_Errors}
\end{table*}

% \newpage
\section{Conclusion}
\label{Conclusion}
In this work, we considered an extension of the finite expression method (FEX) in which some explicit operators in the candidate pool are replaced by TransNet-based neural candidates. Through several high-dimensional PDE examples, we examined whether such trainable candidates can be incorporated into the FEX framework without degrading its ability to identify useful solution structures.

Our numerical results show that the proposed extension performs well when the candidate pool contains candidates associated with the main structure of the true solution. When such candidates are excluded, the method can still combine other available candidates to produce a reasonably accurate approximation, although with a visible loss of accuracy. This indicates that the proposed approach offers additional flexibility in the design of the candidate pool, while also highlighting that the success of FEX remains sensitive to the quality and relevance of the available candidates.

% Overall, this work serves as a preliminary numerical investigation of combining TransNet-style training with FEX. A more systematic study of candidate-pool design, operator training, and the balance between expressiveness and interpretability will be pursued in future work.
% \newpage
% \bmsubsection*{Acknowledgments}
% This material is based upon work supported by U.S. Department of Energy, Office of Science, Office of Advanced Scientific Computing Research, Applied Mathematics program under grants DE-SC0025412 and DE-SC0024703. FB would also like to acknowledge the support from U.S. National Science Foundation through grant DMS-2142672. 

% \bmsubsection*{Financial Disclosure}

% None reported.

% \bmsubsection*{Conflicts of Interest}

% The authors declare no conflicts of interest.

\bibliographystyle{plain}
\bibliography{Ref}
\end{document}